\RequirePackage{fix-cm}
\documentclass[envcountsect,numbook,natbib]{svjour3_arxiv4}   

\usepackage{graphicx}
\usepackage{bm}
\graphicspath{{./}}
\usepackage[utf8]{inputenc}
\usepackage{dsfont}

\usepackage{color}
\usepackage[hyphens]{url}
\usepackage{enumerate} 
\usepackage{booktabs}
\usepackage{tabularx}
\usepackage{ltablex}
\usepackage{mathptmx}
\usepackage{amssymb,amsmath,latexsym} 

\usepackage{mathtools}
\mathtoolsset{showonlyrefs=true}
\usepackage{longtable} 
\usepackage{textcomp}
\usepackage{lscape}
\usepackage{bbm}
\usepackage{eurosym}
\usepackage[font=footnotesize,labelfont=bf]{caption}
\usepackage{subcaption}
\usepackage[usenames,dvipsnames,svgnames]{xcolor}

\usepackage[Algorithm]{algorithm}
\usepackage{algorithmicx}
\usepackage{algpseudocode}
\usepackage{listings}

\usepackage{pdfpages}
\usepackage[pdftex]{hyperref}

\usepackage{scalerel}
\usepackage{import}
\usepackage{multicol} 
\usepackage[nice]{nicefrac} 

\usepackage{tikz}
\usetikzlibrary{arrows,positioning}
\usetikzlibrary{shapes.geometric}

\usepackage[12pt]{extsizes}
\usepackage[normalem]{ulem}
\usepackage{geometry}
\usepackage{ifthen}
\usepackage{calrsfs}

\DeclareMathAlphabet{\pazocal}{OMS}{zplm}{m}{n}
\let\mathcal\pazocal

\allowdisplaybreaks

\spnewtheorem{theorem}{Theorem}{\bfseries}{\itshape}
\spnewtheorem{corollary}[theorem]{Corollary}{\bfseries}{\itshape}
\spnewtheorem{lemma}[theorem]{Lemma}{\bfseries}{\itshape}
\spnewtheorem{proposition}[theorem]{Proposition}{\bfseries}{\itshape}
\spnewtheorem{definition}[theorem]{Definition}{\bfseries}{\itshape}
\spnewtheorem{remark}[theorem]{Remark}{\bfseries}{\upshape}
\spnewtheorem{assumption}[theorem]{Assumption}{\bfseries}{\itshape}
\spnewtheorem{algo}[theorem]{Algorithm}{\bfseries}{\itshape}

\counterwithin{theorem}{section}

\definecolor{myred}{rgb}{0.8,0,0}  
\definecolor{mygray}{rgb}{0.5,0.5,0.5}  
\definecolor{mygreen}{rgb}{0,0.7,0}  
\definecolor{mydgreen}{rgb}{0,0.3,0}  

\newcommand{\neu}{\color{blue} }

\newcommand{\mymarginpar}[1]{ \marginpar{{\tiny #1}}}

\renewcommand{\neu}{}

\renewcommand{\mymarginpar}[1]{}

{\vskip\baselineskip\noindent\textbf{Proof of {#1}:}}%
{\hspace*{.1pt}\hspace*{\fill}\BOX\vskip\baselineskip}

\usepackage[version=4]{mhchem}

\def \R{\mathbb{R}}               
\def \1{{\bf 1}}                
\def \0{{\bf 0}}
\def\qed{\hfill$\Box$}
\DeclareMathOperator{\trace}{\mathrm{tr} \,}

\DeclareMathOperator*{\argmin}{arg\,min}
\DeclareMathOperator*{\Cov}{\textrm{Cov}}

\newcommand{\Tc}{T^{\text{LT,in}}}
\newcommand{\mdot}{\dot{m}}
\newcommand{\PG}{P^{\text{G}}}

\newcommand{\AC}{A}
\newcommand{\aC}{a}
\newcommand{\amin}{\underline{a}}
\newcommand{\amax}{\overline{a}}
\newcommand{\aT}{g^{HF}}

\newcommand{\THTin}{T^\text{HT,in}}
\newcommand{\THTout}{T^\text{HT,out}}
\newcommand{\tHTinmax}{\tau^\text{HT,in}_\text{max}}
\newcommand{\tHToutmax}{\tau^\text{HT,out}_\text{max}}
\newcommand{\tHTin}{\tau^{\text{in}}}
\newcommand{\tHTout}{\tau^{\text{out}}}
\newcommand{\PH}{P^{\text{H}}}
\newcommand{\PHfun}{\pi^{\text{H}}}

\newcommand{\lD}{l^{\text{D}}}
\newcommand{\lC}{l^{\text{C}}}

\newcommand{\Tch}{T^{\text{C}}}
\newcommand{\Tdch}{T^{\text{D}}}
\newcommand{\PW}{P^{\text{W}}}

\newcommand{\epsC}{\varepsilon^{C}}
\newcommand{\epsD}{\varepsilon^{D}}

\newcommand{\TSGin}{T^{\text{SG,in}}}
\newcommand{\TSGout}{T^{\text{SG,out}}}
\newcommand{\FSGin}{F^{\text{SG,in}}}
\newcommand{\FSGout}{F^{\text{SG,out}}}

\newcommand{\Hor}{t_{\text{E}}}

\newcommand{\npump}{n_{\text{H}}}
\newcommand{\tl}{t_a}
\newcommand{\tu}{t_b}
\newcommand{\admiss}{\mathcal{U}}
\newcommand{\noise}{Z} 

\newcommand{\kW}{\text{kW}}

\newcommand{\dmin}{d_\text{min}}
\newcommand{\dmax}{d_\text{max}}
\newcommand{\rmin}{r_\text{min}}
\newcommand{\rmax}{r_\text{max}}

\newcommand{\myparagraph}[1]{\paragraph{#1}}
\graphicspath{{./../}}

\providecommand{\keywords}[1]
{
	\small	
	\textbf{\textit{Keywords---}} #1
}

\newcommand{\refAppendix}[3][0]{%
	\ifthenelse{\equal{#1}{0}}{Appendix \ref{#2}}{\cite[Appendix #3]{Pilling2024}}%
}

\newcommand{\eqrefAppendix}[3][0]{%
	\ifthenelse{\equal{#1}{0}}{\eqref{#2}}{#3}%
}

\title{Reinforcement Learning Methods for the Stochastic Optimal Control of an Industrial Power-to-Heat System}
\titlerunning{Reinforcement Learning Methods for the Stochastic Optimal Control} 

\author{Eric Pilling,   Martin Bähr and Ralf Wunderlich}
\institute{Eric Pilling / Ralf Wunderlich \at
	Brandenburg University of Technology Cottbus-Senftenberg, Institute of Mathematics, P.O. Box 101344, 03013 Cottbus, Germany;  
	\email{\texttt{eric.pilling@b-tu.de} / \texttt{ralf.wunderlich@b-tu.de} }  
	\and
	Martin Bähr \at
	German Aerospace Center, Institute of Low-Carbon Industrial Processes,
	Department Simulation and Virtual Design,     Walther-Pauer-Straße 5,  03046 Cottbus, Germany;
	\email{\texttt{Martin.Baehr@dlr.de} }	
}

\date{Version of \today}

\begin{document}

	\maketitle
	\begin{abstract}
		\noindent
		The optimal control of sustainable energy supply systems, including renewable energies and energy storage, takes a central role in the decarbonization of industrial systems. However, the use of fluctuating renewable energies leads to fluctuations in energy generation and requires a suitable control strategy for the complex systems in order to ensure energy supply. In this paper, we consider an electrified power-to-heat system which is designed to supply heat in form of superheated steam for industrial processes. The system consists of a high-temperature heat pump for heat supply, a wind turbine for power generation, a sensible thermal energy storage for storing excess heat and a steam generator for providing steam. If the system's energy demand cannot be covered by electricity from the wind turbine, additional electricity must be purchased from the power grid. For this system, we investigate the cost-optimal operation aiming to minimize the electricity cost from the grid by a suitable system control depending on the available wind power and the amount of stored thermal energy. This is a decision making problem under uncertainties about the future prices for electricity from the grid and the future generation of wind power. The resulting stochastic optimal control problem is treated as finite-horizon Markov decision process for a multi-dimensional controlled state process. We first consider the classical backward recursion technique for solving the associated dynamic programming equation for the value function and compute the optimal decision rule. Since that approach suffers from the curse of dimensionality we also apply reinforcement learning techniques, namely Q-learning, that are able to provide a good approximate  solution to the  optimization problem within reasonable time.
	\end{abstract}
	
	

	\setcounter{tocdepth}{3}
	\tableofcontents
	
	\section{Introduction}
\begin{figure*}[ht] 
		\centering
		\includegraphics[width=0.92\textwidth]{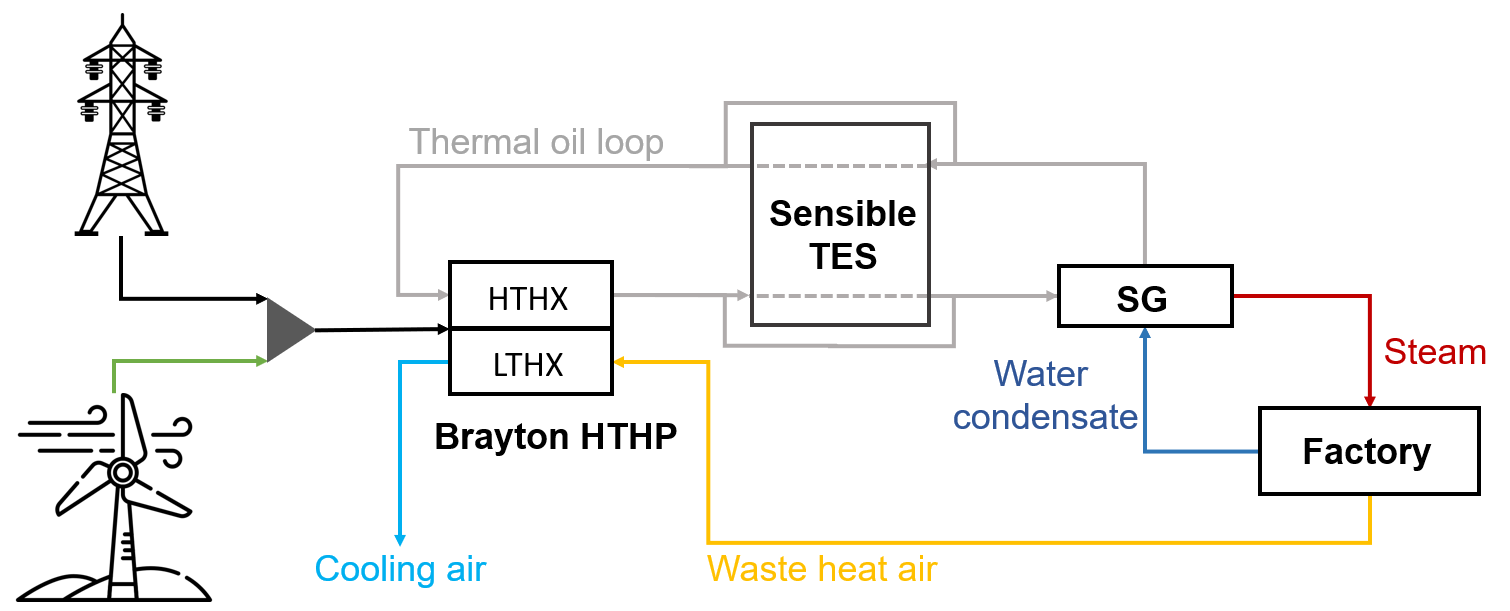}
		\caption{Illustration of the investigated industrial P2H system for electrified steam generation proposed in \cite{Walden2023}. The thermal system consists of an HTHP, a TES and a SG, which are connected via a thermal oil loop. The HTHP uses a waste heat air stream as heat source and is powered by electricity from a WT or the power grid, in order to provide constant heat supply.}
		\label{fig:use_case}
	\end{figure*}

	Nowadays, the supply of process heat for industrial processes by conventional systems leads to high \ce{CO2}-emissions, as these are predominantly based on the combustion of fossil fuels. The electrification of heat generation through the use of novel technologies, such as high-temperature heat pumps (HTHP), is a potential measure to reduce these emissions. In combination with renewable energy sources, the sustainable heat supply for industrial processes is based on complex systems that require realistic modeling as well as cost- and emission-optimized system operation. In particular, electrified energy supply systems face the challenge of determining a cost-optimal operating strategy due to the fluctuating power generation from renewable energies, while ensuring the  required heat demand of the industrial process. In this context, a potential industrial power-to-heat (P2H) system with an HTHP providing heat for a steam generator (SG), see Figure \ref{fig:use_case}, was recently proposed by Walden et al.~\cite{Walden2023}. This P2H system uses the availability of an on-site wind turbine (WT) to generate its own electricity to power the HTHP. This reduces the cost of purchasing electricity from the power grid. These costs can be further reduced by using a thermal energy storage (TES), which serves to balance out the fluctuating generation of renewable energy. The overproduction of electricity can then be stored as thermal energy and used later to supply the system with its own resources.
	
	In \cite{Walden2023}, the cost-optimal operation of this electrified system, with the aim of minimizing the total cost of grid power was treated as a deterministic optimization problem and solved using methods of algebraic nonconvex, nonlinear programming theory. In addition, it was assumed that future wind power generation and electricity prices were already known in advance. However, in real-world scenarios, the problem of optimal management and operation of such systems is a decision making problem under uncertainty, as precise forecasts of future wind energy supply and electricity prices are not possible. Therefore, the problem must be formulated in a stochastic framework.\\ 
	In this context, a typical question needs to be addressed: At what time and at what rate should energy be stored in or withdrawn from the TES to reduce the costs of grid electricity? To address this question, we will treat the cost-optimal management of the underlying industrial P2H system as a stochastic optimal control problem and solve it using Markov decision process (MDP) theory.

	\myparagraph{Literature Review  on  Optimal Management of Industrial P2H Systems} 
	This literature is embedded in numerous studies on optimization problems for energy systems, in particular for electrical and thermal microgrids. However, many of the articles only briefly describe the underlying model and methods. The optimization problems are mainly related to technical aspects, and the control problem is solved with commercial optimization software. The mathematical aspects of the optimal control of energy systems are generally not sufficiently addressed.
	
	The optimal management of \textit{combined heating and power systems} has recently been studied by several authors. In Testi et al. \cite{testi2020stochastic}, an optimal integration of electrically driven heat pumps within a hybrid distributed energy supply system is investigated. There, the authors proposed a multi-objective stochastic optimization methodology to evaluate the integrated optimal sizing and operation of the energy supply systems under uncertainties in climate, space occupancy, energy loads, and fuel costs. In Kuang et al. \cite{kuang2019stochastic}, a stochastic dynamic solution for off-design operation optimization of combined heating and power systems with energy storage is considered. 
     The cost-optimal management of a residential heating system with a geothermal storage under uncertainty is studied in Takam et al. \cite{TakamWunderlich2025}.
    A review on optimal energy management of combined cooling, heating and power microgrids is given in Gu et al. \cite{gu2014modeling}. Further contributions on combined heating, cooling, and power system can be found in  \cite{ehsan2019scenario,zhong2021distributed} and references therein. 
	
	In recent years, \textit{machine learning methods} have also been increasingly used to solve optimization problems.
	For example, Bui et al. \cite{en12091789} and Mohammed et al. \cite{DQNMicrogrid} model a microgrid energy management system with battery storage that is connected to the energy grid and distributed energy sources. To obtain an optimal operation strategy that aims to handle loads, prices and the decision of charging or discharging the battery, Q-learning \cite{Watkins1989} is used. Nakabi and Toivanen \cite{NakabiTahaToivanenPekka} considered a similar application, but used and compared different state-of-the-art reinforcement learning algorithms like Q-learning, deep deterministic policy gradients or proximal policy optimization to achieve the optimal management of their microgrid model. Another application of MDPs is proposed by Yu et al. \cite{DRLSmartHomeEnergyManagement}, who apply reinforcement learning to a home energy management system. In addition to the usage of a battery storage and the connection to the power grid as well as renewable energies, the household utilizes a heating, ventilation and air conditioning system that needs to be operated as cost efficient as possible. Belloni et al. \cite{Belloni2016ASO} use a MDP formulation of a system with a WT and battery storage to obtain the optimal control with dynamic programming. Thermal storage devices combined with a heat pump are used in the papers of Ridder et al. \cite{DERIDDER20112918} and Chenzi et al. \cite{smartgreens21}, which use dynamic programming and Q-learning, respectively.
	
	\myparagraph{Literature Review on Stochastic Optimal Control} 
	The cost-optimal management of energy supply systems under uncertainty can be treated mathematically as a stochastic optimal control problem. There is extensive literature on this theory.
	A considerable part of this literature investigates \textit{dynamic programming} solution techniques. In the continuous-time setting, in which diffusion or jump-diffusion processes form the controlled state process, this leads to the Hamilton-Jacobi-Bellman equation as a necessary optimality condition, see Fleming and Soner \cite{fleming2006control}, Pham \cite{pham2009c}, and Oksendal and Sulem \cite{oksendal2019stochastic}. These nonlinear partial differential equations can usually only be solved using \textit{numerical methods} such as those in Shardin and Wunderlich \cite{shardin2017partially}, Chen and Forsyth \cite{chen2010implications}.	
	For discrete-time models, the theory of \textit{Markov decision processes} provides a solution algorithm based on backward recursion. We refer to Bäuerle and Rieder \cite{Baeuerle2011}, Puterman \cite{Puterman2014}, Hern{\'{a}}ndez-Lerma and Lasserre \cite{HernandezLerma1996} and Powell \cite{Powell2011}. Such MDPs are also obtained by time discretization of continuous-time control problems. 
	
	\noindent For high-dimensional state spaces, the solution of MDPs suffers from the curse of dimensionality. To overcome this problem, powerful numerical methods have been developed in recent years.	
	Examples are the least squares Monte Carlo method introduced in \cite{longstaff2001valuing,tsitsiklis2001regression}, 
	approximate dynamic programming, Q-learning  and related reinforcement learning  methods \cite{Powell2011,Sutton2018,Watkins1989},  
	optimal quantization methods \cite{pages2004optimal} as well as neural network \cite{li2024neural,nielsen2015neural} and deep learning methods \cite{bengio2017deep,hure2021deep}.

	\myparagraph{Our Contribution}
    This article presents a mathematical model for the operation of an industrial P2H system with an HTHP and a TES. It explicitly takes into account the stochasticity of intermittent renewable energy sources such as wind power and the fluctuating market prices for electricity in the power grid. These variables are modeled by suitable stochastic processes that are calibrated to real-world data. Furthermore, the model takes into account that permanent changes in the operating points should be avoided. Therefore, the cost-optimal energy management problem is treated as a discrete-time stochastic optimal control problem, where the controls are kept constant between two discrete points in time. Nevertheless, the dynamics of the state and system variables describing the operation of the system are treated in continuous-time to avoid unnecessary time discretization errors. 
	
	The optimization problem is formulated as an MDP and solved using dynamic programming methods. Since the state of the control problem is three-dimensional, the numerical solution already faces the curse of dimensionality. The problem becomes even more serious when we extend and refine our stylized model to include more details of the HTHP operation. Then the computational effort for the numerical solutions becomes prohibitively high. Therefore, in this paper, we investigate reinforcement learning methods such as Q-learning to find a faster numerical approximate solution. Finally, we present the results of extensive numerical experiments in which we compare the results of the different numerical methods.

	\myparagraph{Paper Organization}
	Section \ref{sec:model_description} is devoted to a thorough mathematical modeling of the considered industrial P2H system. It introduces the state and control variables and additional system variables, as well as the underlying assumptions on the P2H system. In Section \ref{sec:StochasticOptimalControlProblem}, a MDP formulation of the stochastic optimal control problem is derived. This section provides details on the formulation of the stochastic processes for wind speed and the electricity price, state and control constraints, and the cost functions of the optimal control problem. In Section \ref{sec: Backward Dynamic Programming}, the classical approach to solve MDP problems by backward recursion based on the associated Bellman equation is presented. Approximate solutions based on reinforcement learning techniques, in particular Q-learning, are described in Section \ref{sec:MachineLearning}. Finally, Section \ref{sec:NumResults} presents results of numerical experiments in which the optimal control problem is solved using the methods proposed in Section \ref{sec: Backward Dynamic Programming} and \ref{sec:MachineLearning}. 
	An appendix collects proofs and technical results that have been removed from the main text.
	
	\section{Mathematical Modeling of the Industrial P2H System}
	\label{sec:model_description}
	In this section, a mathematical model for the operation of the industrial P2H system is developed. It is treated as a control system with an endogenous state variable that can be influenced by a control variable and exogenous stochastic states. The P2H system is subject to various operational constraints that lead to state and control constraints. For more technical details about the underlying P2H system, we refer the reader to \cite{Walden2023}. In the following, we first briefly introduce the industrial P2H system and then describe its mathematical modeling as a control system.
	
	\subsection{Industrial P2H System}
	\label{sec:ProblemDescr}
	The industrial P2H system based on renewable energy and energy storage shown in Figure \ref{fig:use_case} is designed to supply constant process heat in the form of superheated steam. The system consists of the following four components:  \\	
	\begin{tabular}[t]{l p {0.9\textwidth}}
		(i) &  an on-site \textit{wind turbine} that generates renewable electricity,\\
		(ii) &  a \textit{high-temperature heat pump} for heat supply, which is powered by electricity from the WT or the power grid,\\
		(iii) &  a sensible \textit{thermal energy storage} to store excess energy in times of high wind power production or low grid electricity prices,  and \\
		(iv) &  a \textit{steam generator}  to provide constant process steam.
	\end{tabular}\\[0.5ex]
	The thermal system components (ii)-(iv) are connected via a thermal oil loop. A more detailed system configuration is depicted in Figure \ref{HTHP_scheme}. The HTHP generates high-temperature process heat, which is fed into a thermal fluid loop with thermal oil as the heat transfer fluid (HTF). Via a fluid bypass, the charging factor $\lC \in [0,1]$ determines the proportion of the HTF that is routed through the TES, while the remaining proportion $1-\lC$ passes directly into the SG. The bypass is used to regulate the charging process of the TES. During charging, the hot HTF flows through the cooler TES and heats the storage medium.  
	
	\begin{figure*}[ht]
		\tikzstyle{block} = [draw,rectangle,minimum height=2cm,minimum width=3cm]
		\tikzstyle{block2} = [draw,rectangle,minimum height=2cm,minimum width=2.cm]
		\tikzstyle{block3} = [draw,rectangle,minimum height=2cm,minimum width=4.0cm]
		\tikzstyle{triangle1} = [isosceles triangle,draw,minimum size =0.1cm,isosceles triangle apex angle=60]
		\tikzstyle{triangle2} = [rotate=180,isosceles triangle,draw,minimum size =0.1cm,isosceles triangle apex angle=60]
		\tikzstyle{triangle3} = [rotate=90,isosceles triangle,draw,minimum size =0.1cm,isosceles triangle apex angle=60]
		\tikzstyle{triangle4} = [rotate=270,isosceles triangle,draw,minimum size =0.1cm,isosceles triangle apex angle=60]
		\tikzstyle{ball} = [draw,circle,cross,minimum size=0.5cm]
		\scalebox{0.7}{
			\begin{tikzpicture}[>=latex',every text node part/.style={align=center},cross/.style={path picture={ 
						\draw[black]
						(path picture bounding box.south east) -- (path picture bounding box.north west) (path picture bounding box.south west) -- (path picture bounding box.north east);
				}}]
				\node[block3] at (0,0) (block_HTHP) {{\huge HTHP} \\[8pt] \hspace{0.8cm}{\Large $\PH$}};
				\node[rotate=-90] at (1.7,0) {HTHX};
				\node[rotate=-90] at (-1.7,0) {LTHX};
				\draw[color=black,line width=0.3mm,dashed] (-1.4,1) -- (-1.4,-1);
				\draw[color=black,line width=0.3mm,dashed] (1.4,1) -- (1.4,-1);

				\node[block] at (9,0) (block_S) {{\huge TES}\\[8pt]{\Large $\AC,~~R$}};
				\node[block2] at (18,0) (block_G) {\huge{SG}};
				\node[triangle1] at (-0.315+4.49,-0.64) (block_V1) {};
				\node[triangle2] at (0.315+4.51,-0.64) (block_V2) {};
				\node[triangle3] at (0.315+4.185,-0.97) (block_V2a) {};
				\node[triangle1] at (0.315+12.88,-0.64) (block_B1) {};
				\node[triangle2] at (0.315+13.52,-0.64) (block_B2) {};
				\node[triangle3] at (0.315+13.2,-0.975) (block_B2a) {};    
				\node[triangle1] at (-0.315+13.49,0.64) (block_V3) {};
				\node[triangle2] at (0.315+13.51,0.64) (block_V4) {};
				\node[triangle4] at (0.315+13.185,0.97) (block_V4a) {};
				\node[triangle1] at (0.315+3.875,0.64) (block_B3) {};
				\node[triangle2] at (0.315+4.525,0.64) (block_B4) {};
				\node[triangle4] at (0.315+4.2,0.97) (block_B4a) {};  
				\node at (-1.7,2) (dummy_HTHP1) {}; 
				\node at (-0.2,2) (dummy_HTHP1a) {};     
				\node at (0.8,2) (dummy_HTHP1b) {};  
				\node at (-1.7,-2) (dummy_HTHP2) {};

				\draw[->,color=black,line width=0.4mm] (2.01,-0.64) -- (block_V1) node [pos=0.5, below, color=black] {\Large $ \THTout $};
				\draw[->,color=red,line width=0.4mm] (block_V2) -- (block_S.203) node [pos=0.5, below, color=black] {\Large $\lC\cdot\dot{m}$};
				\draw[->,color=red,line width=0.4mm] (block_S.-23) -- (block_B1) node [pos=0.5, below, color=black] {\Large $\Tch$};
				\draw[->,color=black,line width=0.4mm] (block_B2) -- (17,-0.64) node [pos=0.5, below, color=black] {\Large $\TSGin$};
				\draw[->,color=red,line width=0.4mm] (block_V2a) -- +(0,-1) node [pos=1.5, right, color=black] {\Large ~~~~~$(1-l^\text{C})\cdot\dot{m}$}-- +(9,-1) --+(9,-0.12);
				\draw[->,color=black,line width=0.4mm] (17,0.64) -- (block_V4) node [pos=0.5, below, color=black] {\Large $ \TSGout$};
				\draw[->,color=blue,line width=0.4mm] (block_V3) -- (block_S.23) node [pos=0.5, below, color=black] {\Large $\lD\cdot\dot{m}$};
				\draw[->,color=blue,line width=0.4mm] (block_S.157) -- (block_B4) node [pos=0.5, below, color=black] {\Large $\Tdch$};
				\draw[->,color=black,line width=0.4mm] (block_B3) -- (2.01,0.64) node [pos=0.5, below, color=black] {\Large $ \THTin $};
				\draw[->,color=blue,line width=0.4mm] (block_V4a) -- +(0,1) node [pos=1.5, left, color=black] {\Large $(1-\lD)\cdot\dot{m}$~~~} --+ (-9,1) --+(-9,0.12);
				\draw[->,color=black,line width=0.4mm] (dummy_HTHP1) -- (-1.7,1) node [pos=0, above, color=black] {\Large $\Tc$};    
				\draw[->,color=black,line width=0.4mm] (-1.7,-1) -- (dummy_HTHP2) node [pos=1, below, color=black] {\Large $ T^{\text{LT,out}}$};  
				\draw[-,color=black,line width=0.4mm] (dummy_HTHP1a) -- (-0.2,1.5)  node [pos=0, above, color=black] {\Large ~~$\PW$} -- (0.3,1.5);
				\draw[-,color=black,line width=0.4mm] (dummy_HTHP1b) -- (0.8,1.5)  node [pos=0, above, color=black] {\Large ~$\PG$} -- (0.3,1.5);
				\draw[->,color=black,line width=0.4mm] (0.3,1.5) -- (0.3,1);

			\end{tikzpicture}
		}

		\caption{Detailed flow diagram \cite{Walden2023} of the studied industrial P2H system (cf. Figure \ref{fig:use_case}) with HTHP, TES and SG. The charging and discharging factor $\lC, \lD \in [0,1]$ determine the heat flow to the SG and the HTHP depending on the thermal state of the TES. Exemplarily, the \textit{red} lines indicate the charging mode if $\lC\in (0,1)$, the \textit{blue} lines the discharging mode for $\lD\in (0,1)$. Simultaneous charging and discharging is not allowed. \textit{Charging} mode is characterized by $\lC \in (0,1], ~\lD=0$, \textit{discharging} by $\lD\in (0,1], ~\lC=0$, and  \textit{idle} mode by  $\lC=\lD=0$.}
		\label{HTHP_scheme}
	\end{figure*}
	
	A second bypass from the SG outlet returning to the high-temperature heat exchanger (HTHX) of the HTHP is used to discharge the TES. The discharge factor $\lD \in [0,1]$ specifies the proportion of the HTF that is passed through the TES, such that the remaining part $1-\lD$ enters directly into the HTHP. During the discharging process, the HTF cooled in the SG flows through the warmer TES and lowers the temperature of the storage medium. As the fluid now enters the HTHX at a higher temperature, the HTHP’s electricity consumption is reduced.
	
	In idle operation, characterized by $\lC=\lD=0$, the TES is completely bypassed by the HTF. At the low-temperature heat exchanger (LTHX) of the HTHP, a waste heat air stream from the industrial consumer is used as heat source. We note that the cold air outlet stream at the LTHX is not used for cooling applications in the current configuration. Further, it is assumed that the system components, in particular the HTHP and SG, operate in steady state. This means that the dynamic behavior of the components during operating point changes is neglected.
	
	In our setup, we use the heat flow rates that determine the charging and discharging operations of the TES to describe the control of the P2H system.   These rates have a direct functional relationship with  the HTHX inlet and outlet temperatures, the HTHP's compressor shaft speed and electricity consumption, which is described in more detail below in Subsection \ref{sec:System_Variables_Operational_Constraints}. The HTHP electricity consumption that is not covered by wind energy determines the amount of electricity drawn from the grid, resulting in a direct functional relationship with the running electricity costs that are included in the performance criterion of the optimization problem. More details follow below in Subsection \ref{sec:OperationalCosts}.
	
	The mathematical description of nonlinear component models for HTHP and SG is based on process simulation software that also takes into account the part-load behavior of the heat pump. Based on this, the physical characteristics are then approximated by algebraic surrogate models that appropriately mimic the input-output behavior of the components.  
	For our purposes, it is sufficient to model the state of charge of the TES only by its spatially averaged temperature and neglect the detailed spatial temperature distribution, as for example in \cite{TakamWunderlichPamen2023, TakamWunderlich2025Energies}. This avoids complex calculation of internal heat propagation and facilitates the solution of the optimization problem.
	
	The actual WT power output is typically modeled by a function of the wind speed at rotor height. This dependence is given by the so-called  power curve, which is explained in Subsection \refAppendix{sec:Wind_Turbine_Power}{F.3}.   
	We emphasize that the system and component modeling does not take into account pressure and heat losses as well as friction losses and also no electricity consumption of auxiliary systems such as fluid pumps. The design and dimensioning of the system components was determined by engineering calculations, for which we refer to \cite{Walden2023}.
	
	Recall, the aim is to determine the cost-optimal operation of the P2H system that minimizes the expected total costs over a finite planning horizon $\Hor>0$ from the purchase of grid electricity and the revenues from the sale of WT overproduction, taking into account the uncertainty of the fluctuating wind energy supply and electricity prices. To derive the mathematical formulation of this optimization problem in form of an MDP in Section \ref{sec:StochasticOptimalControlProblem}, we describe in the following the details of state and control variables, additional  system  variables and operational constraints.    
	
	\subsection{Time Discretization}
	\label{sec:Time_Discretization}
	While the state and system variables of the P2H system evolve continuously over time, the control variables available to the controller are typically not changed permanently, but only at discrete points in time and then kept constant until the next time point. This is caused by the fact that operating the HTHP, i.e., changing the HTHX outlet temperature by varying the compressor shaft speed, induces thermal stresses in the heat exchangers during transient operations. For this reason, rapid changes in the operating points should be avoided and limited to a few discrete points in time.
	
	We therefore divide the planning horizon  $[0,\Hor]$ into $N\in \mathbb{N}$ uniformly spaced subintervals of length $\Delta t = {\Hor}/{N}$ and define the time grid points  $t_n = n \Delta t$ for $n =0,\ldots,N$. Let $G:[0,\Hor]\to \R$ be a given continuous-time function. In this work we use the short-hand notation $G_n:=G(t_n)$ for the sampled value at the time grid point $t_n$. The control variables and some related system variables are assumed to be piecewise constant functions on the time grid introduced above and take the values $G(t)=G_n$ for $t\in[t_n,t_{n+1})$ with $n=0,\ldots,N-1$. The remaining system variables and in particular the state variables of the control problem are treated as continuous-time functions governed by certain equations that capture the dynamics of the system. However, the controller only uses the values at the time grid points for the control decisions. 	
	
	The electricity price for trading on the intraday spot market is not quoted continuously over time, but generally only every 15 minutes. For the sake of simplicity we model this price as a continuous-time stochastic process and discretize it accordingly.
	
	\subsection{State and Control Variables}
	\label{sec:state_control}
	\myparagraph{State Variables} 
	For the formulation of the stochastic optimal control problem, the following three variables describe the state of the control system 	
	at time $t \in [0,\Hor]$:\\[0.5ex]
	\hspace*{3em}\begin{tabular}[t]{rlc}
		$R(t)$, &  the average TES temperature &[$^\circ$C],\\
		$W(t)$, &  the wind speed & [m/s], \\
		$S(t)$, &  the electricity price &  [\euro/MWh].
	\end{tabular}\\[0.5ex]	
	These state variables can be divided into endogenous and exogenous quantities. Here, $R$ is the only \textit{endogenous} variable that is subject to the control action, while $W$ and $S$ are exogenous variables and determined outside the model. The storage temperature $R$ changes during charging or discharging operation and is directly related to the associated  heat flow rate, which form the control variable  introduced below.
	
	The \textit{exogenous} states, on the other hand, are stochastic variables, as the wind speed and the electricity price are subject to a certain degree of uncertainty, meaning that future values are not known exactly in advance and are afflicted with considerable forecasting errors. They must therefore be modeled as stochastic processes, where the detailed description is deferred to Subsection \ref{sec:StateDynamics}.
	
	\myparagraph{Control Variable} 
	In our model, we suppose that the operation of P2H system at time $t \in [0,\Hor]$ is controlled by\\[0.5ex]
	\hspace*{3em}\begin{tabular}[t]{rlc}
		$A(t)$, &  the heat flow rate related to the TES& [kW].
	\end{tabular}\\[0.5ex]
	with values in some action or control space $\mathcal{A}\in\R$ which will be specified below.
	In the following, we use the sign convention that a positive heat flow rate corresponds to charging, while negative values of $\AC$ indicate discharging. The idle mode is represented by $\AC(t)=0$. We show in Subsection \ref{sec:System_Variables_Operational_Constraints} that specifying this variable is sufficient to adjust the other system variables describing the HTHP operation accordingly. There it will also be explained how the HTHP electricity demand depends on the control $\AC$, which in turn determines the demand for electrical energy drawn from the grid, when this demand is not fully covered by wind energy, or the overproduction of wind energy that can be fed into the grid. The operational cost associated with the control $\AC$ are derived in Subsection \ref{sec:OperationalCosts}. Note that the controller's choice of the heat flow rate $\AC$ is subject to various constraints, which we explain in detail in Subsection \ref{sec:ControlConstraints}.
	
	\myparagraph{Endogenous State Variable} 
	The control $\AC$ directly determines the dynamics of the only controlled/endogenous state variable in our control system, namely the TES temperature $R$. It results from an energy balance that describes the change in thermal energy in the TES due to the inflow and outflow of energy during charging and discharging.
	Note, that we neglect thermal losses to the environment. Then the change in thermal energy in the time interval $[t_a,t_b]$ with $0\le t_a<t_b\le \Hor$ in the continuous-time setting is given by  $\int_{t_a}^{t_b}\AC(s)\,\mathrm{d} s$. On the other hand, it is equal to $m_\text{s}^{} c_\text{p,s}^{}(R(t_b)-R(t_a))$, where $c_\text{p,s}^{}$ and $m_\text{s}^{}$ denote the specific heat capacity (assumed to be temperature independent) and the mass of the storage medium, respectively. This leads to the following relation for the TES temperature  for $t\in [t_a,t_b]$
	\begin{align}
		\label{storage_dynamic_cont}
		R(t) = R(t_a) +  \frac{1}{m_\text{s}^{} c_\text{p,s}^{}} \int_{t_a}^{t_b}\AC(s)\,\mathrm{d}s.
	\end{align}
	As already mentioned in Subsection \ref{sec:Time_Discretization}, we make the following 
	\begin{assumption}[Piecewise constant control]\label{Ass:ConstantControl}
		The control $\AC$ is kept constant between two consecutive grid points of the time discretization, i.e.
		\begin{align}
			\AC(t) = \AC(t_n)=: \AC_n,\quad t\in[t_n,t_{n+1}),~n=0,\ldots,N-1.
		\end{align}		
	\end{assumption}
	Below in Subsection \ref{sec:System_Variables_Operational_Constraints}, we notice that the assumption of  constant heat flow rates in the TES corresponds to constant oil temperatures at the inlet and outlet of the HTHX within the periods between the time grid points. This avoids rapid changes in the HTHP operation. Only at the time grid points $t_n$ the heat flow rate changes immediately from $\AC_{n-1}$ to $\AC_{n}$, whereby the transient behavior of the HTHP and TES components is neglected.	
	
	Under the above assumption that the heat flow rate $\AC$ is piecewise constant and does not vary within a time period $[t_n,t_{n+1})$, the dynamics \eqref{storage_dynamic_cont} of the TES temperature within such a period simplifies then to 	
	\begin{align}
		R(t) = R_n +  \frac{1}{m_\text{s}^{} c_\text{p,s}^{}} \AC_n (t-t_n), \quad\text{and in particular }\quad R_{n+1} = R_n +  \frac{1}{m_\text{s}^{} c_\text{p,s}^{}} \AC_n \Delta t,
		\label{storage_dynamic}
	\end{align}
	where we recall the notation $R_n=R(t_n)$ and $\Delta t =t_{n+1}-t_n$.		
	
	\subsection{Additional System  Variables and Operational Constraints}	\noindent \label{sec:System_Variables_Operational_Constraints}
The mathematical modeling of the operation of the P2H system, as shown in Figures \ref{fig:use_case} and \ref{HTHP_scheme}, requires the consideration of several additional variables that have not been included in the set of state and control variables. They are referred to as system variables and are subject to certain operational constraints, which are explained in this subsection. These system variables and their dynamics are needed to derive state-dependent control constraints of the control problem that we formulate in Section \ref{sec:StochasticOptimalControlProblem}. However, the system variables are regarded as internal variables whose specific values do not need to be observed by the controller and which are not included in the decision-making process. The latter is based solely on knowledge of the state variables.
	
	To formulate the mathematical model discussed in this work, we make the following simplifying
	\begin{assumption}[Constant HTF mass flow and waste heat temperature]\label{Ass:ConstantMassFlowTcIn}
		The mass flow $\dot{m}$ of the thermal oil stream and the temperature $\Tc$ of the waste heat air stream at the LTHX inlet are constant over the entire period $[0,\Hor]$.
	\end{assumption}
	
	Note that in \cite{Walden2023} the HTF mass flow rate $\dot{m}$ can vary within a certain range, here we assume that $\dot{m}$ is constant. Although this simplification leads to a lower system flexibility, it avoids the introduction of an additional control variable and thus reduces the complexity of the problem as well as the computational effort required to compute the numerical solution.
	
	\subsubsection{Steam Generator}\noindent
	\label{sec:OpConstrSG}
	The constant heat demand of the SG must be satisfied at all times, leading to the following relations between the mass flow $\mdot$ and temperatures at the inlet $\TSGin$ and outlet $\TSGout$ of the SG:
	\begin{align}
		\TSGin = \FSGin(\dot{m})  \quad\text{and}\quad  \TSGout = \FSGout(\dot{m}).
		\label{F3andF4}
	\end{align}	
	The nonlinear functions $\FSGin$ and $\FSGout$ represent surrogate models for the underlying energy balances and are generated using process simulations, see Section \refAppendix{app:SurrogateModels}{B}. 
	According to Assumption \ref{Ass:ConstantMassFlowTcIn}, the mass flow $\mdot$ is constant in our model, so $\TSGin$ and $\TSGout$ are also constant over the entire period $[0,\Hor]$. It is obvious that $\TSGin>\TSGout$, because the SG can simply be understood as a heat exchanger to supply the factory with superheated steam. 
	
	\subsubsection{High-Temperature Heat Pump} \noindent
	\label{sec:HTHP}
	We now describe the relationship between the HTHP, in particular the inlet and outlet temperatures at the HTHX and its electricity demand, and the control variable.
	
	\myparagraph{Piecewise Constant HTHX Inlet and Outlet Temperature}
    According to Assumption \ref{Ass:ConstantControl}, the control variable $\AC$ representing the heat flow rate in the TES is piecewise constant, i.e., $\AC(t)=\AC_n$ in each interval $[t_n,t_{n+1})$ for $n=0,\ldots,N-1$. This property is transferred to the HTHX inlet and outlet temperatures $\THTin$ and $\THTout$, which result from the given system configuration in which the TES is integrated, the fact that $\TSGin$ and $\TSGout$ are constant on $[0,\Hor]$ and the following energy balances. 
	
	In period $[t_n,t_{n+1})$ during charging, we have an inflow of thermal energy from the HTHP to the thermal oil loop with rate $\npump\mdot c_\text{p,f}^{} \THTout(t)$ and an outflow from this loop to the SG with constant rate  $\npump\mdot c_\text{p,f}^{} \TSGin$. Here, $c_\text{p,f}^{}$ denotes the specific heat capacity (assumed to be temperature independent) of the thermal oil and $\npump$ the number of HTHPs operating in parallel in the underlying P2H system, which is not explicitly shown in Figure \ref{HTHP_scheme}). Since we neglect losses to the environment, the difference between the two rates gives the constant inflow rate $\AC_n$ to the TES during the charging process. Recall that $\AC_n>0$ holds in charging mode, thus it follows $\THTout> \TSGin$. 
	
	\noindent Analogously, during discharging, there is an inflow of thermal energy from the SG with constant rate $\npump\mdot c_\text{p,f}^{} \TSGout$ and an outflow to the HTHP with rate $\npump\mdot c_\text{p,f}^{} \THTin(t)$. The difference of the two rates determines the outflow rate from the TES, which is $\AC_n$. Since $\AC_n<0$ during discharging, it follows $\THTin> \TSGout$. 	
	
	\noindent Based on this, we obtain the relations
	\begin{equation}
		\begin{split}
			\AC_n & = \npump \mdot \,c_\text{p,f}^{} (\THTout(t) -\TSGin)> 0, \quad \text{during charging ~and } \\
			\AC_n & = \npump \mdot \,c_\text{p,f}^{} (\TSGout -\THTin(t)) <0, \quad \text{during discharging}
			\label{Control_inlet_outlet0}
		\end{split}	
	\end{equation}    
	from which follows for the HTHX outlet and inlet temperatures
	\begin{align}
		\label{HTHX_inlet_outlet}
		\begin{array}{rllll}
			\THTout(t)&=\THTout_n&=\tHTout(A_n) & \text{with}\quad \tHTout(a)&= \TSGin + \frac{a^+}{\npump \mdot \,c_\text{p,f}^{}},  \quad \text{ and}\\[1ex]
			\THTin(t)&=\THTin_n&=\tHTin(A_n) &\text{with}\quad \tHTin(a)&= \TSGout + \frac{a^-}{\npump \mdot \,c_\text{p,f}^{}},  
		\end{array}		 
	\end{align}
	which are constant in each of the $N$ time periods. Here, $a^{\pm}=\max(\pm a,0)$ denotes the positive and negative parts of $a$. From this, the following unified notation can be derived for the mapping given in \eqref{Control_inlet_outlet0}, which reads
	\begin{align}
		\label{Control_inlet_outlet}
		\AC_n =  \aT(\THTin_n,\THTout_n),	
	\end{align}
	with $\aT(\tHTin,\tHTout) = \npump \mdot \,c_\text{p,f}^{} (\tHTout -\tHTin - (\TSGin-\TSGout))$. Note that $\THTin_n=\TSGout$ during charging, $\THTout_n=\TSGin$ during discharging, and both equations hold true during idle periods. 
	
	\myparagraph{HTHX Outlet Temperature} 
	From \cite{Walden2023} it is known that there is a complex relationship between the HTHX oil outlet temperature $\THTout$ and the HTHX oil inlet temperature $\THTin$, the HTF mass flow $\dot{m}$, the waste heat air temperature $\Tc$ at the LTHX inlet, and the compressor shaft speed $D$, which can be expressed by a surrogate model in terms of a multivariate cubic polynomial $F_1$ as $\THTout(t) = F_1(\THTin(t), \dot{m}, \Tc,D(t))$. 
	As we know from \eqref{HTHX_inlet_outlet} that $\THTin$ and $\THTout$ are stepwise functions, i.e., $\THTout(t)=\THTout_n$ and $\THTin(t)=\THTin_n$ in each interval $[t_n,t_{n+1})$ for $n=0,\ldots,N-1$, and based on the Assumption \ref{Ass:ConstantMassFlowTcIn}, it follows directly that $D$ is also constant between two subsequent discrete time points, so we obtain
	\begin{align}
		\THTout_n = F_1(\THTin_n, \dot{m}, \Tc,D_n).
		\label{Tout}
	\end{align}
	More specifically, the compressor shaft speed $D$ can be varied (only at discrete time grid points) to determine the pressure ratio within the HTHP and thus the thermal oil temperature $\THTout$ at the HTHX outlet.    
	Details on the definition of the polynomial function $F_1$ are provided in Section \refAppendix{app:SurrogateModels}{B}.

	\myparagraph{HTHP Electricity Consumption} 
	Another complex relationship from \cite{Walden2023} holds for the electricity consumption $\PH$ of the HTHP. Again, this is expressed by a surrogate model in the form of a multivariate quadratic polynomial $F_2$ for $n=0,\ldots,N$ as
	\begin{align}
		\PH_n = \npump F_2(\THTin_n, \dot{m}, \Tc, D_n).
		\label{PowerBalance}
	\end{align}
	The factor $\npump$ in \eqref{PowerBalance} as already mentioned above, denoting the number of HTHPs operating in parallel, takes into account that the surrogate model represents only a single HTHP. The detailed definition of $F_2$ can be found in Section \refAppendix{app:SurrogateModels}{B}.
	
	Relation \eqref{PowerBalance} together with Assumption \ref{Ass:ConstantMassFlowTcIn} and the piecewise constant quantities $\THTin_n$ as well as $D_n$ implies that $\PH_n$ is also constant between two subsequent discrete time points. This demand has to be covered by the sum of the power $\PG$ drawn from or fed into the grid and the power $\PW$ generated by the WT, both of which are generally time-varying, meaning that
	\begin{equation}
		\PH_n=\PG(t)+\PW(t),\quad t\in [t_n, t_{n+1}), ~n=0,\ldots,N-1.
	\end{equation}
	For $\PG > 0$, electricity is drawn (purchased) from the grid, for $\PG<0$, electricity is fed (sold) into the grid. Note, that the dependence of the WT power $\PW$ on the wind speed  $W$ is described by the so-called power curve, $\PW=P_{\text{WT}}(W)$. This is a nonlinear function that grows cubically \cite{WANG2019109422} at medium wind speeds until it reaches the rated power. This value is kept constant  for higher speeds and is set to zero above the cut-out speed and below a cut-in speed of the turbine. More details are given in Subsection \refAppendix{sec:Wind_Turbine_Power}{E}. 
	
	\myparagraph{Dependence of Electricity Consumption on the Control} 
	Based on \eqref{Tout} and \eqref{PowerBalance}, the HTHP's electricity consumption $\PH_n$ can be determined at each time point $t_n$ for a given control $\AC$, which determines, by \eqref{HTHX_inlet_outlet}, the temperatures $\THTout_n$ and $\THTin_n$ of the HTF at the outlet and inlet of the HTHX, respectively. 
	Suppose that at time $t_n$ the control is set to be $\AC_n=\aC$, then \eqref{HTHX_inlet_outlet} implies $\THTout_n =\tHTout(a)$ and $\THTin_n = \tHTin(a)$, and recall that $\mdot$ and $\Tc$ are constants. Then, in a first step, the corresponding shaft speed $D_n=d$ is determined by solving \eqref{Tout} for the unknown $d$, i.e., $\tHTout(a)=F_1(d) :=F_1(\tHTin(a), \dot{m}, \Tc,d)$. Since $F_1$ is a cubic polynomial in $d$, the solution is among the real-valued roots of this polynomial. For the given $F_1$ and using the fact that the shaft speed is restricted to values within the interval $[\dmin,\dmax]$, defined by technical conditions \cite{Walden2023}, there is a unique root $d^*=d^*(a)$ in this interval.
	In a second step, the corresponding electricity consumption $\PH_n$ is obtained by substituting $d^*$ into \eqref{PowerBalance} via
	\begin{equation}
		\PH_n=\npump F_2(\tHTin(a),\dot{m}, \Tc, d^*(a)).
		\label{consumption dependency}
	\end{equation}
	To emphasize the dependence of $\PH$ on the control $\AC$, we introduce the function $\PHfun$ that maps $\AC$ to the electricity consumption for $n=0,\ldots,N-1$, i.e., 
	\begin{align}
		\label{PHfun_def}
		\PH_n = \PHfun(A_n), \quad\text{with}\quad \PHfun(a) =\npump F_2(\tHTin(a),\dot{m}, \Tc, d^*(a)).
	\end{align}
    
	\subsubsection{Thermal Energy Storage Operational Modes} \noindent
	\label{sec:OpConstrTES}
	The TES operation can be divided into three operational modes: charging, discharging and idle, where we assume that simultaneous charging and discharging is not allowed. From the previous subsections, the following relations are known for each operational mode:
	\begin{equation}
		\begin{array}{llll}
			\text{Charging:} & \quad \AC_n > 0, & \quad \THTin=\TSGout, &\quad \lC\in(0,1],\\
			\text{Discharging:} &\quad \AC_n < 0, & \quad \THTout_n=\TSGin, &\quad \lD\in(0,1],\\
			\text{Idle:} & \quad \AC_n = 0, & \quad \lC=\lD=0. &
			\label{test}
		\end{array}
	\end{equation}  
	In charging mode, the excess thermal energy is stored in the sensible TES and increases its medium temperature, while during discharging the cooler HTF with temperature $\TSGout$ absorbs heat from the TES, and decreases the temperature of the storage medium. 
    For details on the charging and discharging mode, we refer to Section \refAppendix{app:OperationalModes}{D}. In idle mode, the HTF bypasses the TES completely and obviously implies $\THTout_n=\TSGin$ and $\THTin_n=\TSGout$, which also follows directly from \eqrefAppendix{TSGin}{(D.2)} and \eqrefAppendix{HTHXin}{(D.5)} given in the Supplementary Material. There, we provide further details, including the dependence of the time-varying charging and discharging factors $\lC,\lD$ on the TES temperature $R$ and the chosen control $A_n$.
	
	\section{Stochastic Optimal Control Problem}
	\label{sec:StochasticOptimalControlProblem}		
	In this section, we formulate the stochastic optimization problem using the Markov decision process framework for the cost-optimal energy management of the industrial P2H system introduced above. Most of the MDP theory can be found in the books of Bäuerle and Rieder \cite{Baeuerle2011}, Hernandez and Lerma \cite{HernandezLerma1996}, Powell \cite{Powell2011} and Puterman \cite{Puterman2014}. We would like
	to refer the interested reader to those books for further information about MDPs. The goal is to find the optimal control that minimizes the expected total cost for electricity consumption from the grid over a finite planning horizon, taking into account the uncertainties of future wind energy and electricity prices as well as the revenues from selling overproduction. 
	The derived stochastic control problem consists of the following blocks: state dynamics, state and control constraints, operational cost functions, state and control space, transition operator and performance criterion.
	
	\subsection{State Dynamics} 
	\label{sec:StateDynamics}
	In Subsection \ref{sec:state_control}, we introduced the three state variables $R,W$ and $S$. The dynamics of the endogenous (or controlled) variable $R$ representing the TES temperature is already given in \eqref{storage_dynamic_cont} and \eqref{storage_dynamic}. Here, we focus on the two exogenous states, the wind speed $W$ and the electricity price $S$. Starting with the continuous-time approach, we derive recursions for the state values at the discrete time points $t_n$ for $n=0,\ldots,N$.	
	We decompose $W$ and $S$ into a non-random function that captures the seasonal patterns and a stochastic Ornstein-Uhlenbeck process that is mean-reverting to zero to describe the unpredictable fluctuations.
	
	Throughout this paper, all stochastic processes and random variables are supposed to be defined on a filtered probability space  $(\Omega,\mathcal{F},\mathbb{F},\mathbb{P})$. In particular, that space carries a two-dimensional Brownian motion $B = (B^W,B^S)$, with two independent standard Brownian motions $B^W, B^S$ on $[0,\Hor]$,  which will be used below to drive stochastic differential equations (SDEs) describing the dynamics of $W$ and $S$.
	The filtration $\mathbb{F}$  is assumed to be generated by $B$, that is, $\mathbb{F}=\mathbb{F}^B=(\mathcal{F}^B(t))_{t\in [0, \Hor]}$  with the $\sigma$-algebras $\mathcal{F}^B(t)=\sigma\{B(s), s\le t\}$, augmented by the $\mathbb{P}$-nullsets, so that $\mathbb{F}$ satisfies the usual assumptions. 	
	While, the market price of electricity  $S$ can also take negative values, the wind speed $W$ is always non-negative. We therefore replace $W$ with $\log W$ and assume 
	for  $t\in [0,\Hor]$
	\begin{align}
		\begin{split}
			\log  W(t) &= \mu_W^{}(t) + Y^W(t) \quad \text{ and } \quad
			S(t) = \mu_S^{}(t)  + Y^S(t).
		\end{split}
		\label{ModelDescribtion:OrnsteinUhlenbeckProcesses}
	\end{align}
	Here, the functions $\mu_W^{},\mu_S^{} : [0,\Hor] \to \mathbb{R}$ describe seasonal patterns, and $Y^W$ and $Y^S$ are Ornstein-Uhlenbeck processes defined by SDEs
	\begin{align}
		\begin{split}
			\mathrm{d}Y^W(t) &= - \lambda_W Y^W(t) \, \mathrm{d}t + \sigma_W^{} \, \mathrm{d}B^W(t), \\
			\mathrm{d}Y^S(t) &= -\lambda_S(c_W^{} Y^W(t) + Y^S(t)) \, \mathrm{d}t + \sigma_S^{} \, \mathrm{d}B^S(t),
		\end{split}
		\label{ModelDescription:OUY_t}  
	\end{align}
	with mean reversion speeds $\lambda_W, \lambda_S >0$, diffusion coefficients $\sigma_W^{}, \sigma_S^{} > 0$ and a constant  $c_W^{} \geq 0$. Due to the different natures of wind speed and electricity price, we assume that $\lambda_W \neq \lambda_S$ to simplify our analysis. A positive constant $c_W^{}$ leads to a negative correlation of the wind speed $W$ and the price process $S$, see Lemma \ref{Lemma: E of Y^S, Y^W} and Proposition \ref{prop bivariate dist} below. This is often observed in energy markets.  
	
	The two SDEs above have analytical solutions that allow us to derive the following closed-form expressions for the joint distribution of the pair of random variables $W$ and $S$ as follows, with the proofs provided in Section \refAppendix{Appendix:Section:Proofs}{E.1}.

	\begin{lemma}[Distribution of solutions $Y^W,Y^S$ to the SDEs \eqref{ModelDescription:OUY_t}] \label{Lemma: E of Y^S, Y^W}	
		Let $0\le \tl <\tu\le \Hor$.	Then the solutions of the SDEs \eqref{ModelDescription:OUY_t} on $[\tl, \tu]$ to given initial values 
		$Y^W(\tl)= y_W^{}$ and $Y^S(\tl)= y_S^{}$ with $y_W^{},y_S^{}\in\R$ are  Ornstein-Uhlenbeck processes with
		\begin{align}
			\begin{split}
				Y^W(t) &= y_W^{} e^{-\lambda_W\tau} + \int_{\tl}^t \sigma_W^{} e^{-\lambda_W(t-r)}\, \mathrm{d}B^W(r),  \\
				Y^S(t) &= y_S^{} e^{-\lambda_S\tau} - \lambda_S c_W^{} \int_{\tl}^t e^{-\lambda_S (t-r)} Y^W(r) \mathrm{d} r +  \int_{\tl}^t \sigma_S^{} e^{-\lambda_S (t-r)} \mathrm{d} B^S(r),
			\end{split}
			\label{model_descibtion::closed_form_soulution_OU}
		\end{align}
		for $t\in[\tl,\tu]$ and $\tau = t -\tl$.
		The conditional distribution of the pair  $(Y^W(t), Y^S(t))$ given  $(Y^W(\tl), Y^S(\tl))=(y_W^{},y_S^{})$ is bivariate Gaussian with mean $m_Y^{}(\tau,y_W^{},y_S^{})$ and positive definite covariance matrix $\Sigma_Y(\tau)$ given by 
		\begin{align}
			m_Y^{}(\tau,y_W^{},y_S^{}) = \begin{pmatrix}
				m_{Y^W}(\tau,y_W^{}) \\
				m_{Y^S}(\tau,y_W^{},y_S^{})
			\end{pmatrix}, ~~~~ \Sigma_Y(\tau) = \begin{pmatrix}
				\Sigma_W^2(\tau) & \Sigma_{WS}(\tau) \\
				\Sigma_{WS}(\tau) & \Sigma_S^2(\tau)
			\end{pmatrix},
		\end{align}
		where for $\tau\ge 0$
		\begin{align}
			\begin{split}
				m_{Y^W}(\tau,y_W^{}) &= y_W^{} e^{-\lambda_W \tau},\\
				m_{Y^S}(\tau,y_W^{},y_S^{}) &= y_S^{} e^{-\lambda_S \tau} - \frac{\lambda_S c_W^{}}{\lambda_S - \lambda_W} y_W^{}\big(e^{-\lambda_W\tau} - e^{-\lambda_S\tau}\big),\\ 
				\Sigma_{W}^2(\tau) & = \tfrac{\sigma_W^2}{2\lambda_W}\big(1- e^{-2\lambda_W \tau}\big),\\						
				\Sigma_S^2(\tau) &= \Sigma_{Y^S}^2(\tau) + \tfrac{(\lambda_S c_W^{})^2}{(\lambda_S - \lambda_W)^2} \Big[ \Sigma_{W}^2(\tau) + \tfrac{\sigma_W^2}{\sigma_S^2}\Sigma_{Y^S}^2(\tau) - \tfrac{2 \sigma_W^2}{\lambda_S+ \lambda_W}\big(1- e^{-(\lambda_S+\lambda_W) \tau}\big) \Big],
			\end{split}
			\label{model_describtion:: conditional_mean_variance}
		\end{align}	
		with ~~$\Sigma_{Y^S}^2(\tau)  = \tfrac{\sigma_S^2}{2\lambda_S}\big(1- e^{-2\lambda_S \tau}\big)$, and the covariance
		\begin{align}
			\Sigma_{WS}(\tau)=  -\tfrac{\lambda_S c_W^{}}{\lambda_S -\lambda_W} \Big[\Sigma_{W}^2(\tau) - \tfrac{\sigma_W^2}{\lambda_S+ \lambda_W}\big(1- e^{-(\lambda_S+\lambda_W) \tau}\big) \Big].
		\end{align}
		It holds $\Sigma_{WS}(\tau) \leq 0$ for $c_W^{} \ge 0$ with equality for $c_W=0$.  
	\end{lemma}	
	Combining the above result for $\tl=t_n$ and $ \tu=t_{n+1}=\tl+\Delta t$ with the definitions from \eqref{ModelDescribtion:OrnsteinUhlenbeckProcesses}, and recall the notation $ W_{n}=W(t_n), S_{n}=S(t_n)$ for $n=0,\ldots,N$, we obtain the following result for the joint distribution of $(\log W_{n+1},S_{n+1})$ given $(\log W_{n},S_{n})$. This will be useful for the construction of the transition operator and the corresponding transition kernel for the MDP's state process below in \eqref{model_description::Transition_function}.

	\begin{proposition}[Conditional distribution of $(\log W_{n+1}, S_{n+1})$ given $(\log W_{n}, S_{n})$] 
		\label{prop bivariate dist}
		The conditional distribution of the pair  $(\log W_{n+1}, S_{n+1})$ given  $(\log W_{n}, S_{n})=(\log w,s)$ with $w>0,s\in\R$, is bivariate Gaussian with  mean
		\begin{align}
			m^{WS}_{n+1}(w,s) &= \begin{pmatrix}
				m^{W}_{n+1}(w) \\
				m^{S}_{n+1}(w,s)
			\end{pmatrix},\label{eq:mWS}
		\end{align}
		and the constant  and positive definite covariance matrix $\Sigma=\Sigma_Y(\Delta t)$ and 
		\begin{align}
			\label{mWmS_def}
            \begin{split}
			m^W_{n+1}(w) &= \mu_W^{}(t_{n+1}) +  m_{Y^W}(\Delta t,\log w -\mu_W^{}(t_{n}) ), \\
			m^S_{n+1}(w,s) &= \;\mu_S^{}(t_{n+1}) \,+ m_{Y^S}(\Delta t,\log w -\mu_W^{}(t_{n}), s-\mu_S^{}(t_{n}) ),
            \end{split}
		\end{align}		
		where $m_{Y^W},m_{Y^S}$ and $\Sigma_Y$ are given in Lemma \ref{Lemma: E of Y^S, Y^W}.
	\end{proposition}
	
	Note that the wind speed $W$ follows a log-normal distribution because $\log W$ is Gaussian. The above result on the conditional distribution of the pairs  $(\log W_{n}, S_{n})$ and the fact that the dynamics of the stochastic fluctuations $Y^W, Y^S$ are driven by Brownian motions, i.e., processes with independent increments, can be used to derive a recursion for the discrete-time dynamics of the sequence of these pairs which is driven by a sequence of  independent standard normally distributed random vectors. 
	\begin{corollary}[Discrete-time dynamics of wind speed and energy price]\label{cor_WS-sequence}
		Let the Cholesky decomposition of the symmetric and positive definite covariance matrix $\Sigma$ given in Proposition \ref{prop bivariate dist}	be of the form
		\begin{align}
			\Sigma=AA^\top \quad \text{with}\quad  A = \begin{pmatrix}
				\Sigma_W & 0 \\
				\rho \Sigma_S & \sqrt{1- \rho^2}\Sigma_S
			\end{pmatrix} \quad \text{and}\quad 
			\rho = \frac{\Sigma_{WS}}{\Sigma_W \Sigma_S},
		\end{align}
		where $\rho$ denotes  the associated correlation coefficient. Then there exists a sequence $(\noise_n)_{n=1,\ldots,N}$ of  independent standard normally distributed random vectors $\noise_n=(\noise_n^W,\noise_n^S)^\top \sim \mathcal{N}(0_2,I_2)$ such that 
		\begin{align}
			(\log W_{n+1}, 	S_{n+1}) 
			=  m_{WS}^{}({n+1},W_n,S_n)	
			+ A \noise_{n+1},
		\end{align}
		with $m_{WS}^{}$ given in \eqref{eq:mWS}. Further, it holds
		\begin{align}
			\begin{split}
				W_{n+1} &= \exp\Big( m^W_{n+1}(W_n) + \Sigma_W \noise_{n+1}^W\Big), \\
				S_{n+1} & = m^S_{n+1}(W_n,S_n) + \Sigma_S \Big(\rho \noise_{n+1}^W + \sqrt{1 -\rho^2} \noise_{n+1}^S \Big).
			\end{split}
			\label{disccrete dynamics wind and price}
		\end{align}
		
	\end{corollary}
	\subsection{State Constraints} \noindent
	\label{StateConstraints}
	The operation and technical design of the underlying P2H system restrict the state and control variables. The state-dependent control constraints derived below in Subsection \ref{sec:ControlConstraints} result from the following constraints on the state variables.
	
	Due to the configuration of the P2H system, the TES temperature is bounded, i.e., $R_n \in [r_\text{min},r_\text{max}]$ for all $n$. More precisely, ensuring constant SG inlet and outlet temperatures by fluid bypass regulation implies that the storage temperature $R_n$ cannot be greater than $r_\text{max}=\TSGin$ and cannot fall below $r_\text{min} = \TSGout$. 
	
	In contrast, the exogenous state variables are generally based on our modeling approach in \eqref{ModelDescribtion:OrnsteinUhlenbeckProcesses}. Wind speeds are by nature non-negative and potentially unbounded, implying that $W_n \in (0,\infty)$. Unlike wind speeds, the electricity prices are also allowed to become negative and we have $S_n \in (-\infty,\infty)$. A negative price may occur in the case of overproduction of electricity, while at the same time there is a lower demand in the grid. In addition, in this case, producers are penalized for feeding in additional power, while consumers are rewarded for using electrical energy from the grid.

    \begin{figure}[H]
		\centering
		\hspace*{-0.5cm}\includegraphics[scale=0.4]{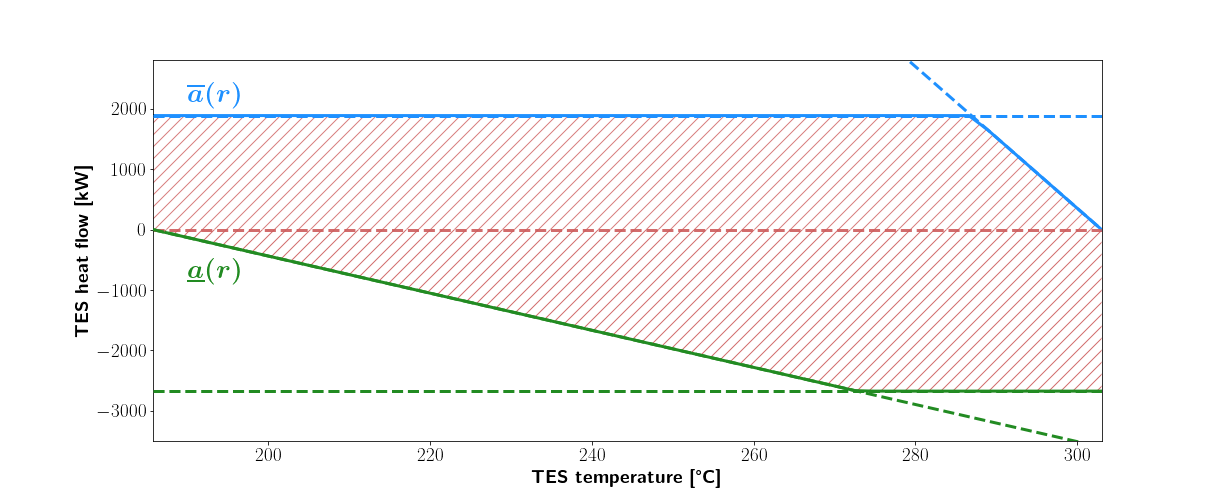}
		\caption{Visualization of the control constraints \eqref{model:describtion action_space_structure} as a function of the TES temperature.			 
			Upper bound  $\amax(r)$ (blue), lower bound  $\amin(r)$ (green),  and the sets  of feasible controls for the  control $\AC_n(x)$ (red). If the TES is  almost full, the upper bound $\amax(r)$ is decreasing and approaches zero to prevent  overheating  during charging. If the TES is not sufficiently full, the decreasing lower limit $\amin(r)$ prevents undercooling during discharging. In both cases, the heat flow is throttled accordingly. The maximum of the positive upper bound $\amax(r)$ and the minimum of the negative lower bound $\amin(r)$ result from the maximal inlet and outlet temperatures of the HTHX.}
		\label{Fig admissible control constraints}
	\end{figure}
    
	\subsection{Control Constraints} \noindent
	\label{sec:ControlConstraints}
	The various state and operational constraints mentioned in the previous subsections imply constraints on the control and lead to state-dependent sets of feasible controls from which the controller can select the actions. 	
	In particular, the control $A_n$ for the period $[t_n,t_{n+1})$ can only be selected such that technical upper limits for the HTHX outlet temperature  $\THTout_n = \tHToutmax$ and the HTHX inlet temperature $\THTin_n = \tHTinmax$ are not exceeded. While the maximal HTHX outlet temperature $\tHToutmax$ is directly related to the maximal compressor shaft speed, the maximal inlet temperature $\tHTinmax$ is set by system constraints of the HTHP. Furthermore,  the controller needs to consider that the time-varying charging and discharging factors $\lC,\lD$, which determine the bypasses, can only take values in $[0,1]$ and must ensure that the TES temperature $R$, which is also time-varying, does not leave the range $[\rmin,\rmax]=[\TSGout,\TSGin]$. Given that the state at the beginning of the period $[t_n,t_{n+1})$ is $X_n=x=(r,w,s)$, the control $A_n$ can be selected from a set of feasible controls  
	\begin{align}
		\mathcal{A}_n(x) = [\amin(r),\amax(r)] \subset \mathcal{A}, 
		\label{model:describtion action_space_structure}
	\end{align}
	where $\amin(r),\amax(r)$ are piecewise linear functions of the TES temperature, which we derive in \eqrefAppendix{upper_aO_contraint}{(G.1)} and \eqrefAppendix{upper_aI_contraint}{(G.2)}. Figure \ref{Fig admissible control constraints} illustrates the derived set of feasible controls and their dependence on $R$ for the system parameters listed in the Supplementary Material \ref{Appendix:TableOfParameters}.

	\subsection{Operational Costs}
	\label{sec:OperationalCosts}
	The operational costs of the system are directly related to the HTHP's electricity consumption $\PH$, which is linked to the HTHX outlet and inlet temperature $\THTout$ and $\THTin$ chosen by the controller, see \eqref{consumption dependency}. Covering the electricity demand depends on the available power output $\PW$ generated by the WT, which in turn is a function of the wind speed $W$. 
	The difference $\PG=\PH-\PW$ must be drawn from the grid at the price $S$ if $\PG>0$. We assume that $\PW$ is free of charge and does not incur any additional costs such as operational and maintenance costs. Consequently, only the consumed grid power $\PG$ must be paid and incurs costs. A negative $\PG$ means an overproduction of WT power that can be sold to the grid for revenue. Here, the selling price is usually lower than the purchase price $S$.

	\myparagraph{Running Cost}
	In our model, we consider the running operational costs $C_n: \mathcal{X} \times \mathcal{A} \to \mathbb{R}$ in each of the periods $n=0,\ldots,N-1$. These are  defined as the expected cumulated costs in the period $[t_n,t_{n+1})$, given that at time $t_n$ the state $X_n=x = (r,w,s)^\top$ and the control $\AC_n=\aC $ is chosen, and read as
	\begin{align}
		C_n(x,\aC) &= \mathbb{E}_{n,x,a}\bigg[\int_{t_n}^{t_{n+1}} \Psi(t,W(t),S(t),\aC) \mathrm{d} t \bigg],
		\label{cost_functional}
	\end{align}
	where $\Psi$ is defined by
	\begin{align}
		\Psi(t,W(t),S(t),a) = S(t)(\PHfun(a) - \PW(t))^+ - \zeta S_\text{sell}(t)(\PHfun(a) - \PW(t))^-,
		\label{Running Cost Psi}
	\end{align}
	with $\PHfun$ introduced in \eqref{PHfun_def}.
	The conditional expectation $\mathbb{E}_{n,x,a}(\cdot) = \mathbb{E}(\cdot| X_n = x,~A_n =a)$ emphasizes the dependence on the current time grid point $t_n$ and the current state $X_n=x$ as well as the action $A_n = a$ selected at this state. Further, we denote by $z^+=\max(z,0)$ and $z^-=\max(-z,0)$ the positive and negative part of $z\in\R$, respectively.
	The functional $\Psi$ in \eqref{Running Cost Psi} is divided into two parts: (i) the costs for buying electricity from the grid at price $S(t)$ and (ii) the revenue for selling excess energy at  a lower price $S_\text{sell}(t)$. Here, $\zeta \in \{0,1\}$ is a user defined model parameter that indicates if selling is allowed or not. If it is not allowed to sell excess energy to generate revenue, we set $\zeta = 0$. In this case, the surplus or overproduction of energy is discarded. For $\zeta = 1$, energy is fed into the grid for a reduced market price $S_\text{sell}$ given by 
	\begin{align}
		S_\text{sell}(t) = S(t) - \eta(t),
	\end{align}
	where $\eta: [0,\Hor] \to \mathbb{R}^+$ is called spread.
	This spread reflects  transaction fees, taxes or the willingness of the grid operator to buy energy only at a certain discount on the market price $S(t)$. 
	A special situation occurs in times of negative market prices $S$, which are often caused by energy overproduction. In this case, buying from the grid leads to a reward, i.e., one gets paid for purchasing energy. Selling, on the other hand, causes additional cost to keep the grid stable due to the abundance of energy. Here, a spread $\eta>0$ leads to a further reduction of the selling price, which results in higher costs for feeding energy into the grid and therefore makes selling less attractive. For more details on the computation of the running costs $C_n(x,a)$ and in particular the conditional expectation in \eqref{cost_functional}, see Subsection \refAppendix{Appedix:Contruction_Cost}{H.1}.
	
	\myparagraph{Terminal Cost} 
	At the end of the planning period, a terminal cost function $G_N: \mathcal{X} \to \mathbb{R}$ can be used to evaluate the terminal state of the system, in particular the amount of thermal energy stored in the TES, in monetary terms. A typical example are \textit{penalty and liquidation payments} that are applied if the TES medium temperature is below or above a certain user-defined critical value $r_\text{crit} \in [r_\text{min}, r_\text{max}]$. Suppose that the terminal state is $X_N=x =(r,w,s)^\top$, then the terminal cost is defined by
	\begin{align}
		G_N(x) = \begin{cases}
			g_\text{Pen}(r)s_\text{Pen}, & r < r_\text{crit}, \\
			g_\text{Liq}(r)s_\text{Liq}, & r \geq r_\text{crit}, 
		\end{cases}
		\label{terminal_penalty}
	\end{align}
	where the functions $g_\text{Pen}: \mathcal{X} \to \mathbb{R}^+$ and $g_\text{Liq}: \mathcal{X} \to \mathbb{R}^-$ describe the amount of thermal energy required to adjust the TES temperature from $r$ to $r_\text{crit}$. In the case of penalization, $g_\text{Pen}$ units of thermal energy must be fed in, while for liquidation, $g_\text{Liq}$ units are withdrawn.
	If the TES is not sufficiently filled, i.e., $r < r_\text{crit}$, a penalty is applied for energy consumption at a fixed price $s_\text{Pen} \geq 0$, depending on the respective temperature difference. In the opposite case, excess energy in the TES is liquidated, which means that energy is sold at the fixed price $s_\text{Liq} \geq 0$, which generates a revenue that appears as a negative terminal cost $G_N$. It should be noted that this definition of the terminal cost includes a worthless expiration of the TES by setting $s_\text{Pen} = s_\text{Liq} = 0$.
	
	\subsection{State and Action Space} \noindent
	Summarizing all the information from above, the state process $X_n \in \mathcal{X}$ of the P2H system at time $t_n$ is described by $X =(X_n)_{n=0,\ldots,N}$ with
	\begin{align}
		X = (R, W, S),
	\end{align}
	where the state space $\mathcal{X} \subset \mathbb{R}^3$ is defined as \hfill 
	\begin{align}
		\mathcal{X} =  [r_\text{min},r_\text{max}]  \times (0,\infty) \times (-\infty, \infty),
		\label{state space definition}
	\end{align}
	with the boundaries according to $r_\text{min} = \TSGout$ and $r_\text{max} = \TSGin$, resulting from system-related and technical constraints as well as model assumptions. 
	The control process $A =(A_n)_{n=0,\ldots,N-1}$ at a given state $X_n$ is specified by  $A_n = u_n(X_n) \in \mathcal{A}$ with decision rules
	\begin{align}
		u_n : \mathcal{X} \to \mathcal{A} , \quad x \mapsto u_n(x) \in \mathcal{A}_n(x),\quad  n=0,\ldots,N-1.
	\end{align}
	The  sequence  $u = (u_n)_{n=0,\ldots,N-1}$ of decision rules is called a policy.
	Moreover, the system at a state $X_n = x$ can be controlled by choosing the action
	$u_n(x) = a $.
	The set of feasible controls $\mathcal{A}_n(x) \subset \mathcal{A}$ in state $x \in \mathcal{X}$ at time $t_n$ is based on the derived control constraints \eqref{model:describtion action_space_structure} in Subsection \ref{sec:ControlConstraints} and reads for $x=(r,w,s)$ as
	\begin{align}
		\mathcal{A}_n(x) = [\amin(r),\amax(r)]. 
	\end{align}

	\subsection{Transition Operator} \noindent
	\label{subsec: transition operrator}
	The transition from one state to another, within the feasible set $\mathcal{X}$, is mathematically described by the transition operator. For state $X_n $, action $A_n= u_n(X_n)$ at time point $t_n$ and a random disturbance $\noise_{n+1} = (\noise^W_{n+1},\noise^S_{n+1}) \sim \mathcal{N}(0_2,\textrm{I}_2)$, the state dynamics of the system is defined by the transition operator as
	\begin{align}
		X_{n+1} = \mathcal{T}_n(X_n,A_n, \noise_{n+1}).		\label{model_description::Transition_function}
	\end{align}
	In this context, according to \eqref{storage_dynamic}, the endogenous state dynamics for the TES temperature is given by $R_{n+1} := g^R(X_n,A_n)$ with			
	\begin{align}
		g^R(x,a) = r + \frac{1}{m_\text{s}^{} c_\text{S}^{}} a \Delta t, \quad  \text{for} \quad x = (r,w,s).
	\end{align}	
	The wind speed $W$ as exogenous, stochastic state is modeled as an exponential discrete-time Ornstein-Uhlenbeck process, see \eqref{disccrete dynamics wind and price}. Its dynamic reads as
	\begin{align}
		W_{n+1} := g^W_{n+1}(X_n,\noise_{n+1}) \quad\text{with}\quad 
		g^W_{n+1}(x,z) = \exp(m^W_{n+1}(w) + \Sigma_W z^W),
	\end{align}
	for $x = (r,w,s)$ and $z=(z^W,z^S)$. The electricity price $S$, the second exogenous and stochastic state, is described using \eqref{disccrete dynamics wind and price} with
	\begin{align}
		S_{n+1} := g^S_{n+1}(X_n,\noise_{n+1}) \quad\text{with}\quad
		g^S_{n+1}(x,z) =  m^S_{n+1}(w,s) + \Sigma_S (\rho z^W + \sqrt{1 -\rho^2} z^S).
	\end{align}    
	Putting all together, the transition operator \eqref{model_description::Transition_function} is given by
	\begin{align}
		\mathcal{T}_n(x,a,z) = \begin{pmatrix}
			g^R(x,a),~
			g^W_{n+1}(x,z),~
			g^S_{n+1}(x,z)
		\end{pmatrix}.
	\end{align}
	For the discrete-time system we will consider a filtered probability space with the  filtration $\mathbb{F}=(\mathcal{F}_n)_{n=0,\dots,N}$ where the  $\sigma$-algebras $\mathcal{F}_n = \sigma(\noise_1,\dots,\noise_n)$ are  generated by the independent random variables $\noise_1,\dots, \noise_n$,  and $ \mathcal{F}_0 = \{ \emptyset, \Omega \}$ is the trivial $\sigma$-algebra.

	\subsection{Performance Criterion and Optimization Problem} \noindent
	\label{sec:optproblem}
	The combination of the discrete-time system and the formulation of the corresponding cost functional allows to determine the operational performance criterion. In this context, the optimal control of the system is related to a cost-optimal policy such that the expected aggregated running costs \eqref{cost_functional} for operating the P2H system and the terminal costs \eqref{terminal_penalty} for an initial state $X_0 = x \in \mathcal{X}$ are minimized.  A policy $u = (u_n)_{n=0,\ldots,N-1}$ is a sequence of decision rules $u_n: \mathcal{X} \to \mathcal{A}$ that maps a given state $x \in \mathcal{X}$ to an admissible action $a \in \mathcal{A}_n(x)$. At each point of time, the associated objective function or performance criterion $J^u_n:\mathcal{X}\to\mathbb{R}$ is given by
	\begin{align}
		J_n^u(x) = \mathbb{E} \bigg[ \sum_{i=n}^{N-1} C_i(X_i, u_i(X_i)) + G_N(X_N)~ \bigg|~ X_n =x \bigg],
		\label{performance criterion}
	\end{align}
	with the running and terminal cost defined in Subsection \ref{sec:OperationalCosts}. We denote by $\admiss$ as the set of all admissible policies $u$ such that the objective function \eqref{performance criterion} is well-defined and $\AC_n = u_n(X_n)$ satisfies the control constraints for all $n=0,\ldots,N-1$.
	An admissible policy $u\in  \admiss $ is called optimal if  
	\begin{align}
		J_n^{u^*}(x)= V_n(x) = \underset{u \in \admiss }{\inf \,} J_n^u(x).
	\end{align}
	The function $V_n$ is called \textit{value function} and describes the minimal expected aggregated running costs. Therefore, finding the value function is equivalent of finding the optimal policy $u^*$.
    
	\myparagraph{Recursion Property and Bellman Equation} 
	In practice, it is not tractable to minimize over the space of all admissible policies $\admiss$. We will see that the performance criterion \eqref{performance criterion} satisfies a recursion property from which we are able to deduce an alternative way for obtaining the value function.  
	The following theorems can be found in Bäuerle and Rieder \cite[p. 21-23]{Baeuerle2011}. Theorem \ref{thm1} states that the objective function \eqref{performance criterion} fulfils a recursion property.
	\begin{theorem}[Recursion property]
		Let $u = (u_n)_{n=0,\ldots,N-1}$ be a fixed policy. Then for every $n = 0, \ldots, N-1$ and $x \in \mathcal{X}$ the objective function $J^u_n(x)$ satisfies
		\begin{align}
			\begin{split}
				&J_N^u(x) = G_N(x), \\
				& J_n^u(x) = C_n(x, a) + \mathbb{E}_{n,x,a} [J_{n+1}^u(X_{n+1})] .
			\end{split}
		\end{align}
		\label{thm1}
	\end{theorem}
	Using Theorem \ref{thm1}, we obtain that the value function $V_n(x)$ is the solution of the well-known \textit{Bellman equation}.
	\begin{theorem}[Bellman equation]
		For every $x \in \mathcal{X}$, the value function $V_n(x)$ for $n = 0, \ldots, N$ satisfies the Bellman equation
		\begin{align}
			\begin{split}
				& V_N(x) = G_N(x), \\
				&V_n(x) = \underset{a \in \mathcal{A}_n(x)}{\inf} \bigg\{ C_n(x,a) + \mathbb{E}_{n,x,a} [V_{n+1}(X_{n+1}) ]\bigg\}.
			\end{split}
			\label{bellman}
		\end{align}
		\label{thm2}
	\end{theorem}
	The Bellman equation reduces the problem of finding an optimal policy $u^* \in \admiss $ to a recursion in which only the optimal actions for each time point $n$ must be found.  
	
	\myparagraph{(Optimal) State Action Function and Properties} 
	Many applications and algorithms use an alternative performance criterion for $J^u_n$, which is often called the state action function and is denoted by $Q_n^u$ \cite{Bertsekas1996}. It is defined for all $(x,a) \in \mathcal{X} \times \mathcal{A}_n(x)$  as 
	\begin{align}
		Q_n^u(x,a) = \mathbb{E}_{n,x,a} \bigg[ \sum_{i=n}^{N-1} C_n(X_i, u_i(X_i)) + G_N(X_N) \bigg].
		\label{state action performance criterion}
	\end{align}
	Note, that the expectation is conditional to $X_n =x, u_n(x) = a$, and given a policy $u$ with deterministic decision rules $u_n$ it holds that
	\begin{align}
		J_n^u(x) = Q^u_n(x,u_n(x)),
	\end{align}
	for all $n =0,\ldots,N$. As the name suggests, the function $Q^u_n$ assigns a value to each feasible state action pair $(x,a)$, instead of a single value for each state $x$ as is the case with the performance criterion $J^u_n$. By Theorem \ref{thm1}, we get that $Q_n^u$, satisfies the recursion 
	\begin{align}
			Q_n^u(x,a) &= C_n(x,a) + \mathbb{E}_{n,x,a} [J_{n+1}^u(X_{n+1}) ]\nonumber\\
			&= C_n(x,a) + \mathbb{E}_{n,x,a} [Q_{n+1}^u(X_{n+1},u_{n+1}(X_{n+1}))].
	\end{align}
	The optimal state action function $Q_n^*(x,a)$ for $(x,a) \in \mathcal{X} \times \mathcal{A}_n(x)$ is given by
	\begin{align}
		Q^*_n(x,a) = \underset{u \in \admiss }{\inf \,} Q^u_n(x,a).
	\end{align}
	It further relates to the original value function $V_n(x)$ by
	\begin{align}
		V_n(x)  =   \underset{a \in \mathcal{A}_n(x)}{\inf \,} Q_n^*(x,a),
		\label{mdpModel:: valuefunction_stateactionfunction_relation}
	\end{align}
	with optimal policy $u^*  = (u_n^*)_{n=0,\ldots,N-1}$ and solves the Q-version of the Bellman equation \cite{annurev-statistics-031219-041220} below.
	\begin{theorem}[Q-version Bellman equation]
		For every $(x,a) \in \mathcal{X} \times \mathcal{A}_n(x)$, the optimal state action function $Q^*_n(x,a)$ for $n = 0,\ldots,N$ satisfies
		\begin{align}
			\begin{split}
				&Q_N^*(x,a) = G_N(x), \\
				& Q_n^*(x,a) = C_n(x,a) +  \mathbb{E}_{n,x,a} \bigg[\inf\limits_{a^\prime \in \mathcal{A}_{n+1}(X_{n+1})}   Q_{n+1}^*(X_{n+1},a^\prime) \bigg].
			\end{split}
			\label{QBellman}
		\end{align}
		\label{thm5}
	\end{theorem}
	The main difference between the Bellman equation \eqref{bellman} and the Q-version \eqref{QBellman} is that the expectation and minimization are interchanged. This offers certain computational advantages both in the calculation of the expected value and in the minimization, which is now a deterministic problem.
	
	\section{Backward Dynamic Programming}
	\label{sec: Backward Dynamic Programming}
	In this section, we introduce an algorithm that solves the stochastic optimal control problem. This method is inspired by the Bellman equation \eqref{bellman} and is known as backward dynamic programming (BDP). We will also discuss issues that naturally arise when using BDP and propose techniques on how to address them.
    
	\subsection{Backward Recursion Algorithm} \noindent
	Note that given a state $X_n = x$ and action $\AC_n=a$, the next state $X_{n+1}$ is obtained by the transition operator \eqref{model_description::Transition_function}, which depends on the random disturbance $\noise_{n+1}$. Furthermore, the randomness in the subsequent state is only induced by the disturbance allowing to take the expectation with respect to the random variable $\noise_{n+1}$, instead of $X_{n+1}$ leading to
	\begin{align}
		\mathbb{E}_{n,x,a} [V_{n+1}(X_{n+1}))] = \mathbb{E}_{} [V_{n+1}(\mathcal{T}_n(x,a,\noise_{n+1}))].
		\label{TransitionOperator_Expectation}
	\end{align}
	Algorithm \ref{backwardDP} summarizes the backward recursion procedure of BDP, which will be used to solve the MDP derived in Section \ref{sec:StochasticOptimalControlProblem}.
	
	\begin{algorithm}[ht!]
		\caption{Backward Dynamic Programming} \label{backwardDP}
		\begin{algorithmic}
			\State{\textbf{1:}} Initialize $V_N(x) = G_N(x)$ for all $x \in \mathcal{X}$ and set $n = N-1$
			\State{\textbf{2:}} Compute for all $x \in \mathcal{X}$ the value function at time $n$
			\begin{align*}
				V_n(x) = \underset{a \in \mathcal{A}_n(x)}{\inf \,} \bigg\{ C_n(x,a) + \mathbb{E}_{} [V_{n+1}(\mathcal{T}_n(x,a,\noise_{n+1}))] \bigg\}
			\end{align*}
			and the associated optimal decision rule by
			\begin{align*}
				u_n^*(x) \in \underset{ a \in \mathcal{A}_n(x)}{\argmin} \bigg\{ C_n(x,a) + \mathbb{E}_{} [V_{n+1}(\mathcal{T}_n(x,a,\noise_{n+1}))] \bigg\}
			\end{align*}
			\State{\textbf{3:}} If $n > 0$ set $ n = n-1$ and go to step 1 else stop the algorithm.
		\end{algorithmic}
	\end{algorithm}
    \vspace{-0.5cm}
	\subsection{Approximate Solution of the Bellman Equation} \noindent
	BDP may face some issues. For instance, if the state and action space are large or high-dimensional, it may suffer from the curse of dimensionality, or in the case of continuous spaces, the optimization problem in step 2 of Algorithm \ref{backwardDP} may be difficult to solve. Another practical issue is the calculation of the expected value for all $x \in \mathcal{X}$, which may be computationally intractable due to the unavailability of closed-form expressions.
	Since it may be difficult to solve the Bellman equation exactly, we need to make certain simplifications in order to solve the issues described above.
	\myparagraph{State and Action Space Discretization} 
	Firstly, the state space $\mathcal{X} \subset \mathbb{R}^3$ given in \eqref{state space definition} is discretized into distinct grid points. The value function is then calculated and saved for the given reference grid points.
	For the discretization, the seasonalities $\mu_W^{}$ and $\mu_S^{}$ of the Ornstein-Uhlenbeck processes \eqref{ModelDescribtion:OrnsteinUhlenbeckProcesses}, are used to construct time-varying sets of grid points. The advantage of introducing this time dependence is that the value function is only calculated for regions of interests, i.e., subsets of $\mathcal{X}$ that are more likely to appear at certain times. This leads to the family of discretized state spaces  
	\begin{align}
		\widetilde{\mathcal{X}}_n = \{r_1,\ldots,r_{n_\text{R}} \} \times 
		\{w_{1,n},\ldots,w_{n_\text{W},n} \} \times \{s_{1,n},\ldots,s_{n_\text{S},n} \}, \quad  n=0,\ldots,N-1,
		\label{NumRes:: discrete_spaces}
	\end{align}
	with $n_R^{},n_W^{},n_S^{} \in \mathbb{N}$. The specific choice of grid points used in our calculations can be found in Section \refAppendix{Appendix: State Discretization}{I}. Note, that the discretization points for the TES temperature is the same at every time point and only those of wind speed and electricity price change.
	The discrete structure of the state space $\widetilde{\mathcal{X}}_n$ allows us to solve the Bellman equation for each $x \in \widetilde{\mathcal{X}}_n$ separately.
	To solve the problem of minimization over the action space $\mathcal{A}_n(x) \subset \mathbb{R}^2$ in BDP, the action space $\mathcal{A}_n(x)$ is also discretized into grid points $\widetilde{\mathcal{A}}_n(x) = \{a_1,\ldots,a_{n_A} \}$ for given $n_A \in \mathbb{N}$. Hence, the minimization consists of calculating a value for each action and picking the action with the smallest value. 
	
	\myparagraph{Approximation of the Expected Value} 
	Given a state $X_n = x$ and action $\AC_n=a$, the conditional expected value in the Bellman equation with respect to the next state $X_{n+1}$ is given by an unconditional expectation with respect to the random disturbance $\noise_{n+1}$, see \eqref{TransitionOperator_Expectation}. Let $\mathcal{Z} = \{z_1,\ldots,z_L\}$ be a set of value that $\noise_{n+1}$ can take and denote $\widehat{\noise}_{n+1}$ as the discrete random variable taking values in $\mathcal{Z}$. Further let $p_l = \mathbb{P}(\widehat{\noise}_{n+1} = z_l)$ be the corresponding probability that $\widehat{\noise}_{n+1}$ takes value $z_l$, then the expected value can be approximated as weighted sum
	\begin{align}
		\mathbb{E}_{} [V_{n+1}(\mathcal{T}_n(x,a,\noise_{n+1}))]  \approx  \mathbb{E}_{} [V_{n+1}(\mathcal{T}_n(x,a,\widehat{\noise}_{n+1}))] = \sum_{l=1}^L p_l^{} V_{n+1}(\mathcal{T}_n(x,a,z_l)).
	\end{align}
	The set $\mathcal{Z}$ is called a \textit{quantizer} of $\noise_{n+1}$ and defines a partition on $\mathbb{R}^2$ into $L$ subsets $\mathcal{C}(z_l),~l=1,\ldots,L$, where each point $z_l$ is uniquely assigned to a subset. As a consequence the probability $p_l$ corresponds to the probability that $\noise_{n+1}$ takes values in $\mathcal{C}(z_l)$, more precisely
	\begin{align}
		p_l= \mathbb{P}(\noise_{n+1} \in \mathcal{C}(z_l)).
	\end{align}
	The calculation of these probabilities often requires solving high-dimensional integrals over the subsets $\mathcal{C}(z_l)$ with respect to the density of $\noise_{n+1}$. Below, we will explain how to numerically obtain these probabilities. But before we come to this, we want to point out a practical problem that is caused by the quantizer $\mathcal{Z}$.
	
	\myparagraph{Interpolation and Extrapolation} 
	In general, the states $x_l = \mathcal{T}_n(x,a,z_l)$ do not coincide with the grid points in $\widetilde{\mathcal{X}}_{n+1}$ for which the value function $V_{n+1}$ is calculated and saved. If the points are allocated in between existing values, an interpolation can be used to determine $V_{n+1}(x_l)$, otherwise this value must be calculated by extrapolation. In this paper, we use linear interpolation if $x_l = \mathcal{T}_n(x,a,z_l)$ is in between existing point and extrapolate with the value of the nearest neighbor in the set $\widetilde{\mathcal{X}}_{n+1}$. The corresponding extrapolation errors are usually larger than those resulting from interpolation. In any case, the value $V_{n+1}(x_l)$ is weighted by the probability $p_l$, and if this is small, the corresponding error introduced will be scaled down by $p_l$. For this reason, the probabilities $p_l$ can be used to reduce and control the errors in the calculation of the value function. We also note that an appropriate choice of the discretized state space $\widetilde{\mathcal{X}}_{n+1}$ helps to mitigate extrapolation. Details on the construction of the discretization that takes this fact into account can be found in Section \refAppendix{Appendix: State Discretization}{I}.
	
	\begin{remark}
		Due to the approximation of the expected value and action space discretization, the value function obtained is an approximation of $V$ and will be denoted by $\widetilde{V}$. Therefore, the control corresponding to $\widetilde{V}$ is an approximation of the optimal control. The calculation of the value function $\widetilde{V}_n$ for $x \in \widetilde{\mathcal{X}}_n$ in the BDP Algorithm \ref{backwardDP} reduces to 
		\begin{align}
			\widetilde{V}_n(x) = \underset{a \in \widetilde{\mathcal{A}}_n(x)}{\min \,} \bigg\{ C_n(x,a) + \mathbb{E}[\widetilde{V}_{n+1}(\mathcal{T}_n(x,a,\widehat{\noise}_{n+1}))] \bigg \}.
		\end{align}
	\end{remark}
	
	\subsection{Optimal Quantizer for the Expected Value}  \noindent
	The choice of the quantizer $\mathcal{Z}$ is important in order to obtain a good approximation of the expected value in \eqref{TransitionOperator_Expectation}. An approach to obtain them was proposed by Pagès \cite{PagèsPrintems,Pagès}, in which a so-called optimal quantizer $\mathcal{Z}^* = \{z_1^*,\ldots,z_L^* \}$
	is selected. In the following, we will briefly discuss the optimality of this quantizer as well as some theoretical definitions and results. For more information, we refer the reader to the work of Pagès \cite{PagèsPrintems,Pagès}.

	\myparagraph{Optimal Quantizer of $\noise$}  
	For a square integrable random variable $\noise$ in $\mathbb{R}^2$ with probability density of $\noise$ by $f_{\noise}$, we denote the set $\mathcal{Z} = \{z_1,\ldots,z_L\}\subset \mathbb{R}^{2}$ and Borel-measurable function $q: \mathbb{R}^2 \to \mathcal{Z}$. The random vector $q(\noise)$ is called a $L$-quantization of $Z$ and $\mathcal{Z}$ is called $L$-quantizer. The aim is to find a $L$-quantization $q$ such that the quadratic distortion $D_L^{\noise}$ given by
	\begin{align}
		D_L^{\noise}(\mathcal{Z}) = \mathbb{E}(\parallel \noise - q(\noise)\parallel_2^2)
	\end{align}
	is minimized. It can be shown that the so-called Voronoi $L$-quantization defined by
	\begin{align}
		q_{\text{Vor}}(z) = \sum_{l=1}^L z_l \mathbbm{1}_{\mathcal{C}(z_l)}(z)
	\end{align}
	is minimizing $D_L^{\noise}$, where $\mathcal{C}(z_l),l=1,\ldots,L$ are Voronoi-cells with
	\begin{align}
		\mathcal{C}(z_l) \subset \{ y \in \mathbb{R}^d |~\parallel z_l-y\parallel_2 \leq \parallel z_i-y\parallel_2,~i=1,\ldots,L\}.
	\end{align}
	Note, that the Voronoi-cells form a partition of $\mathbb{R}^2$. These quantizations can be understood as the nearest neighbor projection of $\noise$ onto the set $\mathcal{Z}$. We denote the Voronoi $L$-quantization of $\noise$ by $\widehat{\noise} = q_{\text{Vor}}(\noise)$. Moreover, the probability that $\widehat{\noise}$ takes the value $z_l$ is given by
	\begin{align}
		p_l = \mathbb{P}(\widehat{\noise} = z_l) = \mathbb{P}(\noise \in \mathcal{C}(z_l)) = \int\limits_{\mathcal{C}(z_l)} f_{\noise}(z)\, \mathrm{d}z.
	\end{align}
	For the Voronoi $L$-quantization $\widehat{\noise}$ the quadratic distortion can be written as
	\begin{align}
		D_L^{\noise}(\mathcal{Z}) = \mathbb{E}( \parallel\noise - \widehat{\noise}\parallel_2^2) = \sum_{l=1}^L\mathbb{E}(\mathbbm{1}_{\mathcal{C}(z_l)}(Z)\parallel\noise-z_l\parallel_2^2) = \int\limits_{\mathbb{R}^2} f_{\noise}(z)\underset{1\leq l \leq L}{\min\,} \parallel z - z_l\parallel_2^2\,\,\mathrm{d} z. 
	\end{align}
	Now for $\widehat{\noise}$ the mapping $\mathcal{Z} \mapsto D_L^X(\mathcal{Z})$ is continuous and yields a minimum $\mathcal{Z}^* = \{z_1^*,\ldots,z_L^*\}$ having distinct components \cite{PagèsPrintems}. This set is called an optimal $L$-quantizer and satisfies
	\begin{align}
		D_L^{\noise}(\mathcal{Z}^*) = \underset{\mathcal{Z} \subset \mathbb{R}^{2}}{\min \,} D_L^{\noise}(\mathcal{Z}).
	\end{align}
	
	The existence of an optimal quadratic $L$-quantizer and the convergence are proven in \cite{Pages2018}. In addition, Zardor's theorem \cite{Zador1982} provides a prescribed level of accuracy for the number of quantization points $L$. Apart from an upper bound on the quadratic distortion with an error rate of order $L^{-1/d}$, where $d$ is the dimension of the random variable (in our case $d = 2$), this theorem also establishes asymptotic convergence to zero as the number of quantization points $L$ goes to infinity.
	
\begin{figure}[H]
		\centering
		\includegraphics[scale=0.41]{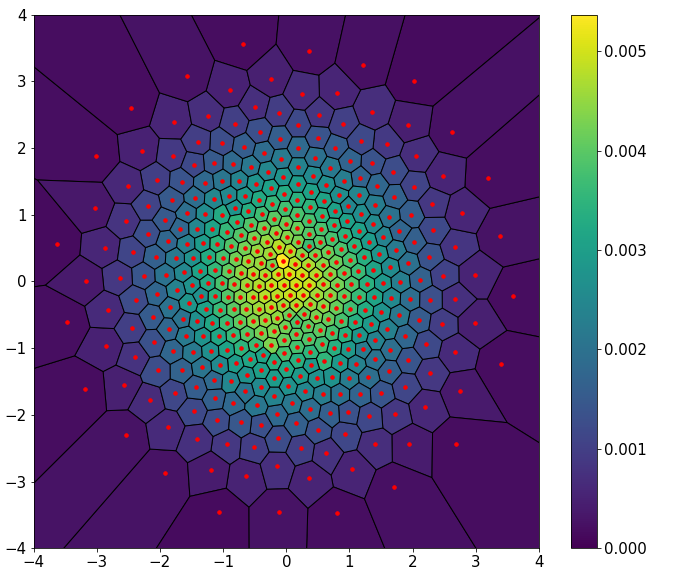}
		\caption{An optimal quadratic $200$-quantizer (red dots) with Voronoi-cells for a standard bivariate Gaussian random variable, taken from \url{www.quantize.maths-fi.com}. The color of the Voronoi-cells indicates their corresponding probability mass.}
		\label{2D_opt_quantizer_gaussian}
	\end{figure}

	\myparagraph{Calculation of Quantizers and Probabilities}
	Numerical methods such as the Competitive Learning Vector Quantization (CLVQ) or (randomized) Llyods algorithm are often used to compute optimal quantizers, see \cite{Pages2015,PagèsPrintems,montes:tel-03112849}. For standard normally distributed random variables in $\mathbb{R}^d$, pre-calculated optimal quadratic $L$-quantizers for different $L,d \in \mathbb{N}$, with their corresponding probability mass of the Voronoi-cells are available on \url{www.quantize.maths-fi.com}. Due to the accessibility of the high precision precalculated optimal quantizers $\mathcal{Z} = \{z_1^*,\ldots, z_L^* \}$ and probabilities $p_l,~l=1,\ldots,L$ of $\noise_{n+1}$, we will use them in this paper. The optimal quadratic $200$-quantizer $\mathcal{Z}^*$ of $\noise$ with its respective Voronoi-cells and corresponding probability mass is shown in Figure \ref{2D_opt_quantizer_gaussian}. 
	
	\myparagraph{Application to the Bellman Equation}  
	When applying the optimal quantization to \eqref{bellman}, it is obvious to replace the continuous random variable $\noise$ by its optimal quantization $\widehat{\noise}$, associated with the optimal quantizer $\mathcal{Z}^*$ in order to obtain a reasonable approximation of the expected value. However, $g(\widehat{\noise})$ is in general not an optimal quantizer for $g(\noise)$, when applying the nonlinear transformation $g(z) = V_{n+1}(\mathcal{T}_n(x,a,z))$. If $g$ is bounded and continuous, then the convergence of $\mathbb{E}_{}[\widehat{\noise}]$ to $\mathbb{E}_{}[\noise]$ as the number of quantization points $L$ grows to infinity, see \cite{Pages2018}, implies the convergence $\mathbb{E}_{}[g(\widehat{\noise})] \to \mathbb{E}_{}[g(\noise)]$. If we make further assumptions on $g$, the convergence in the sense
	\begin{align}
		\underset{L \to \infty}{\lim} L^{\alpha_g} | \mathbb{E}[g(\noise)] -\mathbb{E}[g(\widehat{\noise})] | \leq  C_{g,\noise}, 
	\end{align}
	can be proven for different classes of functions $g$ and precise convergence rates $\alpha_g>0$ can be formulated, for more details see \cite{Lemaire2020}. The constant $C_{g,\noise}$ depends on the properties of $g$ and the disturbance $Z$. In particular, if $g$ is a Lipschitz continuous function with Lipschitz coefficient $L_g$, we obtain that
	\begin{align}
		| \mathbb{E}[g(\noise)] -\mathbb{E}[g(\widehat{\noise})] | \leq L_g \mathbb{E}[\parallel \noise -\widehat{\noise} \parallel_1 ] \leq L_g \mathbb{E}[\parallel \noise -\widehat{\noise} \parallel_2] = L_g \sqrt{D^\noise_L (\mathcal{Z}^*)}.
	\end{align}

	\section{Reinforcement Learning Techniques} \noindent
	\label{sec:MachineLearning}
	This section presents reinforcement learning algorithms that can tackle some of the problems mentioned in the context of BDP.
	We first introduce a quite general class of algorithms called temporal difference (TD) learning methods and then study Q-learning as a special case. Their aim is to approximate the value function in an appropriate parameter space and to construct the optimal policy with respect to this approximation. These methods rely on gradient descent to update the parameters with information obtained by generating samples of the controlled state process. In practice, the information for these methods does not have to come from an explicit model, as in our case. Instead, it can also be provided as data from an observed real world process, which is why these algorithms are often referred to as model-free. 
	
	\subsection{Temporal Difference Learning} \noindent
	\label{subsec:TD_Learning}
	In the following, we use a function approximation to approximate $V_n(x)$ for all $x \in \mathcal{X}$ and $n = 0, \ldots, N-1$. Let $\theta_n \in \mathbb{R}^p$ be a parameter vector that describes an approximation $\overline{V}_n$ of the exact value function $V_n$ in terms of $p \in \mathbb{N}$ parameters $\theta_n^1,\cdots,\theta^p_n$
	\begin{align}
		\overline{V}_n: \mathcal{X} \times \mathbb{R}^p \to \mathbb{R},
	\end{align}
	such that $V_n(x) \approx \overline{V}_n(x,\theta_n)$. Let us give some typical examples of $\overline{V}_n(x,\theta_n)$ below.
	
	\myparagraph{Linear Function Approximation}  In this class of functions \cite{Konidaris_Osentoski_Thomas_2011,Sutton2018}, $\overline{V}_n(x,\theta_n)$ is represented by a linear combination of basis or ansatz functions $\phi_i,~i=1,\ldots,p$ with $\phi_i: \mathcal{X} \to \mathbb{R}$ as
	\begin{align}
		\overline{V}_n(x,\theta_n) = \sum_{i=1}^p \theta_n^i \phi_i(x).
		\label{mdp_theory:: linear approximation}
	\end{align}
	The parameter vector $\theta_n \in \mathbb{R}^p$  corresponds to the coefficients of the linear combination. Polynomial ansatz functions, Fourier basis functions or radial basis functions (RBF) are examples of function classes that are used for the linear function approximation $\overline{V}_n$. 
	
	\myparagraph{Feedforward Neural Networks}  
    Feedforward neural networks (FNN) are simple artificial neural networks and are a popular choice for nonlinear function approximators \cite{DeVore_Hanin_Petrova_2021,ELFWING20183}.  Essentially, they consist of affine-linear maps and nonlinear activation functions. Let $d_0 = \mathbb{R}^d$ and $d_L \in \mathbb{R}$ denote the input and output dimension of the FNN, then  $\overline{V}_n(x,\theta_n)$ is represented by the recursion
	\begin{align}
		\overline{V}_n(x,\theta_n) = A_L \rho(A_{L-1} \rho(\cdots \rho(A_1 x + b_1) \cdots) + b_{L-1}) + b_L),
	\end{align}
	where $L$ is the number of layers, $A_l \in \mathbb{R}^{d_l \times d_{l-1}}$ and $b_l \in \mathbb{R}^{d_l},~ l=1,\ldots,L,$ are weights and biases for each layer with width $d_l \in \mathbb{N}$ and $\rho: \mathbb{R} \to \mathbb{R}$ is a nonlinear activation function that is applied component-wise. The parameter vector $\theta_n$ is the collection of all matrices $A_l$ and vectors $b_l$. Examples for activation functions are the  \textit{sigmoid function} $\rho(x) = \frac{1}{1 + e^{-x}}$ or the \textit{rectified linear unit} (ReLU) $\rho(x) = \max\{x,0\}$. Nowadays, FNNs are frequently used because it is known that they fulfill the universal approximation property \cite{Hornik1989}, i.e., any continuous function can be approximated arbitrarily well.
	
	\myparagraph{TD-Learning Loss Functional} 
	The corresponding parameter update for TD-learning can be derived by minimizing a loss functional given by the expected squared distance 
	\begin{align}
		\mathcal{L} (\theta_n) = \frac{1}{2}\overline{\mathbb{E}}  [(V_n(X_n) - \overline{V}_n(X_n,\theta_n))^2],\quad n = 0,\ldots,N-1.
		\label{TDfunctional}
	\end{align}
	The underlying distribution in the expectation of loss \eqref{TDfunctional} is called the state distribution and is used to sample states $X_n$. Normally, this distribution is chosen such that it reflects the importance of certain states that $X_n$ can take, i.e., states that are of interest for the controller or are likely to appear. A natural choice would be the distribution of $X_n$ or the so-called steady-state distribution \cite{Sutton2018}. For a given policy $u$, it describes the likelihood of $X_n$ taking a specific state for a given initial state. Sampling from the steady-state distribution is realized by creating trajectories starting from an initial state, while following the policy $u$. Here, the initial state is sampled from a predetermined distribution, for example a uniform distribution over the state space $\mathcal{X}$.
	Minimizing the loss with respect to $\theta_n$ can be achieved by gradient descent, which leads to an iterative update rule 
	\begin{align}
		\theta_n^{k+1} = \theta_n^k -  \alpha_n^k ~ \nabla_{\theta_n} \mathcal{L}(\theta_n^k).
		\label{gradientUpdate}
	\end{align}
	Here $k$ denotes the current iteration of the parameters and $\alpha_n^k > 0$ is the step size or learning rate.
	Note, that by interchanging the gradient with the expectation, we formally obtain
	\begin{align}
		\begin{split}
			\nabla_{\theta_n} \mathcal{L}(\theta_n) &=  \overline{\mathbb{E}} [(V_n(X_n) - \overline{V}_n(X_n,\theta_n) ) \nabla_{\theta_n} \overline{V}_n(X_n,\theta_n)].
		\end{split}
	\end{align}
	Thus, we can obtain samples for the gradient of the loss $\mathcal{L}(\theta_n)$ by using samples of $X_n$ according to the state distribution in \eqref{TDfunctional}. To achieve a good and unbiased estimator for the expected value of the gradient, multiple realizations with batch size $M \in \mathbb{N}$ are used and averaged.
	The update rule \eqref{gradientUpdate} is therefore replaced by
	\begin{align}
		\theta_n^{k+1} = \theta_n^k - \alpha_n^k \frac{1}{M} \sum_{j=1}^M  \delta_n^j \nabla_{\theta_n} \overline{V}_n(x_{n}^j,\theta_n^k),
		\label{SGD_update}
	\end{align}
	with $\delta_n^j = V_n(x_n^j) - \overline{V}_n(x_n^j,\theta_n^k)$. The iterative gradient update rule \eqref{SGD_update} is a special case of the Robbins-Monro algorithm \cite{Robbins1951} and is referred to as stochastic gradient descent (SGD) \cite{Garrigos2023}.
	By applying \eqref{bellman}, we get
	\begin{align}
		\delta_n^j =  \underset{a \in \mathcal{A}_n(x_n^j)}{\inf} \big\{ C_n(x_n^j,a) + \mathbb{E}_{n,x_n^j,a} [V_{n+1}(X_{n+1}) ]\big\}  - \overline{V}_n(x_n^j,\theta_n^k).
		\label{TD-1}
	\end{align}
	There are some problems that need to be addressed before a parameter update can be performed.
	
	\myparagraph{Practical and Computational Issues} 
	Firstly, since $V_{n+1}$ is unknown, we need to replace it with an approximation. One way to do this is to use a Monte-Carlo estimation of the performance criterion. Here, multiple trajectories starting from $X_n = x_n^j$ and thus multiple realizations of the performance criterion are obtained and averaged. However, the optimal policy is also unknown as it requires knowledge about the value function. Thus, the best choice for the policy is the (sub)optimal policy induced by the value function approximation. 
	
	A more common approach is \textit{bootstrapping}, which avoids following trajectories with a possibly suboptimal policy. In doing so, the value function $V_{n+1}$ is replaced by its corresponding parameterization $\overline{V}_{n+1}$. This brings some computational advantages, but at the expense that convergence towards the value function can not always be shown, which will be discussed below.
	
	Another issue arises from the replacement of the value function $V_{n+1}$ and the minimization in \eqref{TD-1} which is performed with respect to this approximation. If $u$ is the policy induced by these value function approximations, the optimal control obtained by the associated decision rule is given by $a_n^j = u_n(x_n^j)$. Note, however, that this control may not be optimal in the sense that the minimum over $\mathcal{A}_n(x_n^j)$ is attained and therefore degrades the approximation. 
	
	Last but not least, as for BDP, there are several ways to calculate the expectation in \eqref{TD-1}. For BDP, a quantization approach is used, which could as well be used here. Nevertheless, a more common approach is to use a Monte-Carlo simulation and use samples of $X_{n+1}$.

	\myparagraph{TD-Learning Update} 
	In practice, the TD-learning methods use bootstrapping to replace the  $V_{n+1}$ and one-sample Monte-Carlo estimates instead of extensive calculations of the expectation. This is mostly motivated by the fact that calculations as well as simulation of the state process and evaluation of the optimal policy is time-consuming and therefore computationally intensive. Bootstrapping also offers the advantage to update the parameters $\theta_n$ immediately after observing samples of $X_n$, where Monte-Carlo estimates of the performance criterion must wait until the trajectories end. Parameter updating is performed using samples
	\begin{align}
		(x_n^j,a_n^j,x_{n+1}^j)_{j=1,\ldots,M},\quad n = 0, \ldots, N-1,
		\label{trajectories}
	\end{align}
	with $a_n^j = u_n(x_n^j)$ and $x_{n+1}^j = \mathcal{T}_n(x_n^j,a_n^j,z_{n+1}^j)$, where $z_{n+1}^j$ is a realization of $\noise_{n+1}$.
	This results in the following TD-learning update for the parameters at time $n$
	\begin{align}
		\begin{split}
			\theta_n^{k+1} = \theta_n^k -  \alpha_n^k \frac{1}{M} \sum_{j=1}^M  \delta_n^j  \nabla_{\theta_n} \overline{V}_n(x_n^j,\theta_n^k),
		\end{split}
		\label{TD-learning}
	\end{align}
	with  temporal difference
	\begin{align}
		\delta_n^{j} =  C_n(x_n^j,a_n^j) + \overline{V}_{n+1}(x_{n+1}^j,\theta_{n+1}^k) - \overline{V}_n(x_n^j,\theta_n^k).
		\label{TDoriginal}
	\end{align}
	The scalar $\delta_n$ can be interpreted as the change in information when moving from state $x_n$ to $x_{n+1}$. 
	
	\subsection{Q-Learning} \noindent 
	\label{section::Q-learning}
	Q-learning was first proposed by Watkins in 1989 \cite{Watkins1989} and is essentially a special class of TD-learning. The starting point here is the state action function $Q^u_n(x,a)$ given in \eqref{state action performance criterion} instead of the objective function $J^u_n(x)$ from \eqref{performance criterion}.
	Analogous to the derivation of the TD-learning method, a parameter vector $\theta_n \in \mathbb{R}^p$ is used for all $n = 0,\ldots,N-1$ to describe the parameterization
	\begin{align}
		\overline{Q}_n(x,a,\theta_n): \mathcal{X} \times \mathcal{A} \times \mathbb{R}^p \to \mathbb{R},
	\end{align}
	in order to approximate $Q^*_n(x,a)$. Due to the relation in \eqref{mdpModel:: valuefunction_stateactionfunction_relation}, an approximate value function can also be derived by
	\begin{align}
		\overline{V}_n(x,\theta_n) = \underset{a \in \mathcal{A}_n(x)}{\min \,} \overline{Q}_n(x,a,\theta_n).
		\label{approximate_function_relation}
	\end{align}
	Again, the aim is to minimize the following loss functional
	\begin{align}
		\mathcal{L}^Q(\theta_n) = \frac{1}{2}\overline{\mathbb{E}}  [(Q_n^*(X_n, A_n) - \overline{Q}_n(X_n, A_n,\theta_n))^2].
		\label{Q_learning_functional}
	\end{align}
	As for the TD-learning loss function \eqref{TDfunctional}, a sample state distribution is used. Note, that $u_n(X_n) =A_n$ is a random variable, which emphasizes the choice of an action distribution to obtain samples for $A_n$. For a given state $x$, this distribution can be directly defined by a probability measure on $\mathcal{A}_n(x)$ or by a selection policy $u^S = (u^S_n)_{n=0,\ldots,N-1}$, with $A_n = u_n^S(X_n)$. The relation \eqref{approximate_function_relation} motivates to choose $u^S$ such that the optimal action is sampled frequently. A natural choice is a \textit{greedy} selecting policy
	\begin{align}
		u_n^S(x) = \underset{a\in \mathcal{A}_n(x)}{\min\,} \overline{Q}_n(x,a,\theta_n),
		\label{natural_greedy_selection}
	\end{align}
	with respect to the current approximation. Analogous to TD-learning using \eqref{QBellman} in loss \eqref{Q_learning_functional} in combination with SGD, this results in an iterative Q-learning update for the parameters $\theta_n$
	\begin{align}
		\begin{split}
			\theta_n^{k+1} = \theta_n^k - \alpha_n^k \frac{1}{M} \sum_{j=1}^M  \delta_n^j ~\nabla_{\theta_n} \overline{Q}_n(x_n^j, a_n^j,\theta_n^k) ,
		\end{split}
		\label{Qlearningupdate}
	\end{align}
	with iteration counter $k$, step size $\alpha_n^k > 0$ and temporal difference
	\begin{align}
		\delta_n^j = C_n(x_n^j,a_n^j) + \mathbb{E}_{n,x_n^j,a_n^j} \bigg[\inf\limits_{a^\prime \in \mathcal{A}_{n+1}(X_{n+1})} Q_{n+1}^*(X_{n+1}, a^\prime)\bigg]  - \overline{Q}_n(x_n^j, a_n^j,\theta_n^k).
		\label{Q-1}
	\end{align}
	\myparagraph{Q-Learning Update} 
	As above, we need to tackle similar problems to define a parameter update, such as replacing the unknown $Q_{n+1}^*$ and evaluating the expected value in \eqref{Q-1}. The most common variant of Q-learning uses bootstrapping to replace $Q_{n+1}^*$ and one-sample estimates for the expected value, which results in the temporal difference
	\begin{align}
		\delta_n^j = C_n(x_n^j,a_n^j) + \inf\limits_{a^\prime \in \mathcal{A}_{n+1}(x_{n+1}^j)} \overline{Q}_{n+1}(x_{n+1}^j, a^\prime,\theta_{n+1}^k)  - \overline{Q}_n(x_n^j, a_n^j,\theta_n^k).
	\end{align}
	
	\begin{remark}\label{rem_convergence}
		To ensure convergence for iterative stochastic approximation methods \cite{Bertsekas1996}, like SGD, the step sizes $\alpha_n^k$ for all $n = 0, \ldots, N-1$ need to satisfy the Robbins-Monro conditions \cite{Robbins1951}
		\begin{align}
			\sum_{k=0}^\infty \alpha_n^k = \infty \quad\text{and}\quad \sum_{k=0}^\infty (\alpha_n^k)^2 < \infty.
			\label{RobbinsMonroCond}
		\end{align}
		However, as a consequence of bootstrapping in TD-learning and Q-learning, the resulting gradient estimates may differ from those of the original underlying loss function \eqref{TDfunctional}. These kind of methods are different from SGD and are called semi-gradient methods and require separate convergence results as well as additional assumptions. It should be noted that the main convergence analysis for TD-learning and Q-learning is based on MDPs with infinite time horizon. Convergence results for linear function approximation \eqref{mdp_theory:: linear approximation} are provided in \cite{Melo2008}, \cite{Melo2007}. Nonlinear function approximators such as neural networks require additional assumptions and techniques like projection \cite{NEURIPS2019_98baeb82} or linearization \cite{pmlr-v119-xu20c} to guarantee convergence. Although these convergence results hold for infinite-horizon MDPs, they can still be applied to the finite-horizon case by augmenting the state with time as an additional state variable. Provided that we can formulate an equivalent problem with the augmented state, the aforementioned convergence results can be applied.
		
		Without augmentation, convergence results for finite-horizon TD-learning and Q-learning are derived in \cite{Asis2019} using linear and nonlinear function approximations. These results are based on the recursive properties of the value functions and therefore require less restrictive assumptions than the infinite-horizon setting.
	\end{remark}
    
	\myparagraph{$\varepsilon$-Greedy Selection Policies} 
	Given a state  $X_n = x$, the selection policy $u^S$ is used to create samples of actions that ensure reasonable exploration of the action space $\mathcal{A}_n(x)$ and the state space of $X_{n+1}$. In the following, we discuss a commonly used class of selection policies that differ from the greedy policy \eqref{natural_greedy_selection}. The greedy policy has one major disadvantage that can lead to poor approximations. If this greedy strategy is fully exploited, all decisions are based on the current function approximation, which itself is biased by approximation errors. This can lead to suboptimal action samples in the sense that the optimal action may not be sampled frequently. As a result, the approximation of $\overline{Q}$ will be poor as well as the value function approximation associated with it. We can compensate this by sampling actions from a uniform distribution on $\mathcal{A}_n(x)$. This approach has the advantage that all actions are selected with equal probability and an overall better approximation can be obtained for all actions. A drawback, however, is that the optimal action may again not get sampled very often. In the literature, this problem of choosing an appropriate $u^S$ is known as the exploration-exploitation dilemma. Here, we want to exploit the optimal policy induced by the function approximation as much as possible while still sampling a reasonable number of other actions to not miss out on more optimal actions.
	In practice, a combination of random and greedy policies is used, known as \textit{$\varepsilon$-greedy policy} with exploration rate $\varepsilon \in [0,1]$. Here, an action is either drawn from an uniform distribution on $\mathcal{A}_n(x_n)$ with probability $\varepsilon$ or with probability $1-\varepsilon$ selected greedily as in \eqref{natural_greedy_selection}. The exploration rate could be given by a simple linearly decaying scheme with $\varepsilon^k = \frac{\varepsilon_0}{k}$
	for $\varepsilon_0 \in [0,1]$, where $k$ denotes the iteration counter. 
	\begin{remark}
		In Powell \cite{Powell2011}, examples of exploration rates and step sizes are discussed. It is also mentioned that choosing an appropriate step size that satisfies the Robbins-Monro conditions \eqref{RobbinsMonroCond} is hard in practice. It may happen that the step sizes decrease too quickly, so the parameters converge to a non-optimal solution. Hence, it is suggested to use small constant step sizes $\alpha_n^k = \alpha_0 >0$ as it has been empirically observed that these work well in applications, although the second condition in \eqref{RobbinsMonroCond} is violated. 
	\end{remark}
	
	\subsection{Experience Replay} 
	\label{ML::improving_Q_Learning}
	Q-learning faces the problem that creating trajectories and samples can be a time-consuming task, especially if the dynamics of the system are difficult to simulate or the time horizon is large. This makes it intractable to use generated samples only once and then throw them away when the parameter update is completed.
	We will use a technique called \textit{experience replay} \cite{LinLong-Ji,Mnih2013} that solves the problem of wasting generated samples. This uses a so-called replay memory or replay buffer $\mathcal{R}$ with size $\mathcal{N}_R \in \mathbb{N}$ to store samples $(x_n,a_n,x_{n+1}) \in \mathcal{R}$ for $n=0,\ldots,N-1$ and replay (reuse) them as needed in batches $(x_n^j,a_n^j,x_{n+1}^j),j=1,\ldots,M$ to update the parameters. Sampling past experience, for example, from a uniform distribution on $\mathcal{R}$, also helps to overcome the exploration-exploitation dilemma. Here, the samples obtained in the early phase of the algorithm are used repeatedly, making the choice of the exploration rate less important for the performance of the algorithm. We summarize the Q-learning method with replay buffer in Algorithm \ref{QlearningAlg}.
	\begin{algorithm}
		\caption{Q-Learning with Replay Buffer} 
		\begin{algorithmic}
			\State{\textbf{1:}} Initialize $(\theta_n^0)_{n=0}^{N-1}$; 
			set the maximum number of iterations $k^{\max}>0$ and $k=0$, \\batch size $M$; choose a selection policy $u^S =(u_n^S)_{n=0}^{N-1}$  
			\State{\textbf{2:}} Set $n=0$; choose the initial states $x_0 \in \mathcal{X}$ 
			\While {$n < N$}
			\State Select an action according to $a_n = u_n^S(x_n)$, observe $x_{n+1}$ and store $(x_n,a_n,x_{n+1})$ in  $\mathcal{R}$.
			\State Sample batch $(x_n^j,a_n^j,x_{n+1}^j)_{j=1}^M$ from the replay buffer $\mathcal{R}$
            \For{$j = 1,\dots, M$} calculate $\delta_n^j$
			\If{$n<N-1$}  
			\begin{align*}
				\delta_n^j = C_n(x_n^j,a_n^j) +  \min\limits_{ a \in \mathcal{A}_{n+1}(x_{n+1}^j)} \overline{Q}_{n+1}(x_{n+1}^j, a,\theta_{n+1}) - \overline{Q}_n(x_n^j, a_n^j,\theta_n^k)
			\end{align*}
			\Else{} set \quad ~~$\delta_n^j = C_{N-1}(x_{N-1}^j,a_{N-1}^j) +  G_N(x_N^j) - \overline{Q}_{N-1}(x_{N-1}^j, a_{N-1}^j,\theta_{N-1}^k)$	
			\EndIf
            \EndFor
			\State Choose $\alpha_n^k \in [0,1]$ and update parameters
            \begin{align*}
                \theta_{n}^{k+1} = \theta^k_n - \alpha_n^k \frac{1}{M} \sum_{j=1}^M  \delta_n^j  \nabla_{\theta_n} \overline{Q}_n(x_n^j, a_n^j,\theta_n^k) 
            \end{align*}
			
			\State set $n = n+1$
			\EndWhile 
			\State{\textbf{3:}} Set $k = k+1$; if $k =k^{\max}$ go to step 4; else go to step 2
			\State{\textbf{4:}}  Obtain optimal control for $x \in \mathcal{X}$ and $n = 0,\ldots,N-1$: $u^*_n(x)  \in \underset{a \in \mathcal{A}_n(x)}{\argmin \,} \overline{Q}_n(x, a,\theta_n^{k^{\max}})$
		\end{algorithmic}
		\label{QlearningAlg}
	\end{algorithm}
	Calculating the minimum over all actions might be hard to accomplish if the action space is continuous, as in our case. In our numerical experiments, the action space $\mathcal{A}_n(x)$ is discretized as for BDP, which was explained in Section \ref{sec: Backward Dynamic Programming}.
	
	\section{Numerical Results} \noindent
	\label{sec:NumResults}
	In this section, the numerical results obtained by Algorithms \ref{backwardDP} and \ref{QlearningAlg}, i.e., BDP and Q-learning, are presented and compared with each other. More precisely, we compare the results in terms of accuracy and computational effort of the computed solutions (value function and trajectory). A time horizon of $\Hor=120$ hours is selected for the numerical simulation to find the value function and optimal operation for the industrial P2H system during a working week (5 days). The system and algorithm parameters can be found in Tables \ref{ParametersModel} and \ref{NumericalResultsParameters}, respectively. To keep the results of the proposed algorithms comparable, the same action space discretization is used for both. The evaluation of the minimization over the discretized action space is done by calculating all action values and selecting the action with the minimal value. Penalty costs are applied at terminal time if the TES temperature is below a certain threshold value. In the experiments, the storage must be at least half-full, i.e., the critical value is set to $r_\text{crit} = (r_\text{max}-r_\text{min})/2$. Furthermore, falling below $r_\text{crit}$ is penalized with a penalization price $s_\text{Pen} = 90$ \euro/MWh. We do not reward the liquidation of the TES energy and therefore set $s_\text{Liq} = 0$ \euro/MWh. Selling excess energy into the grid is also not allowed and we set $\zeta = 0$. For the calculation of the value function, both methods made use of a time-varying state space discretization, see Section \refAppendix{Appendix: State Discretization}{I}. For visualization and convenience, we display the value function for fixed grid points selected from the set $[\rmin, \rmax] \times [2,23.5] \times [15,55]$. This set is chosen such that it is a subset of $\mathcal{X}_n$ for all $n=1,\ldots,120$. In order to visualize the results with respect to the three-dimensional state space we will fix some state variables to specific values and plot against the remaining variables. In the following, these fixed values are chosen as the corresponding centers of the $r,w$ and $s$-axis.
	
	\subsection{Backward Dynamic Programming} \noindent 
	
\begin{figure}[t]
		\centering
         \subfloat[\centering]{
         \begin{minipage}{\textwidth}
            \begin{center}
                \hspace*{-0.6cm} \includegraphics[scale = 0.44]{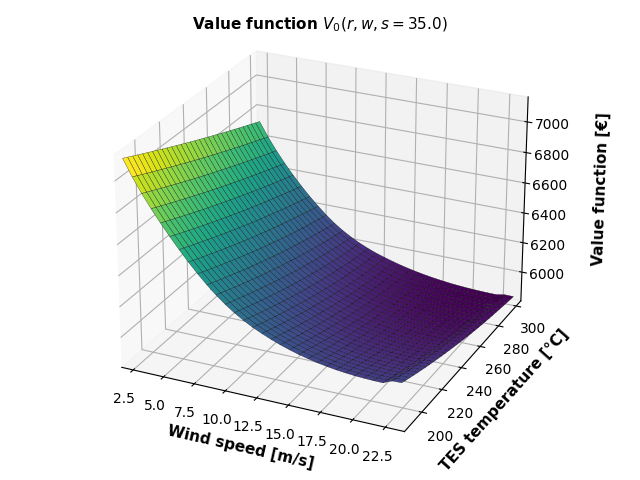}
			  \hspace*{-0.3cm} \includegraphics[scale = 0.44]{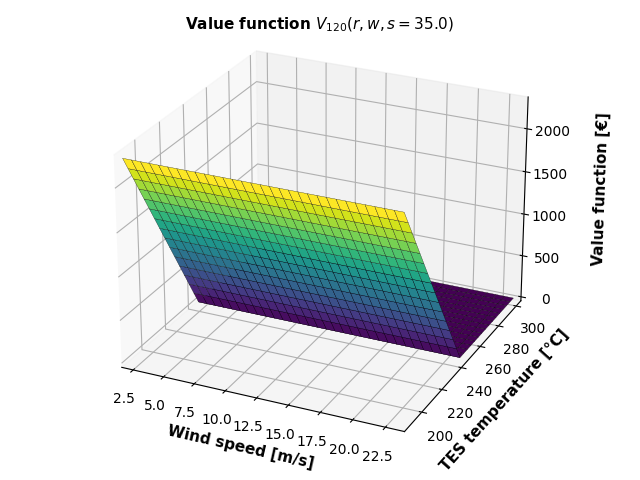}
            \end{center}
        \end{minipage}
        \label{NumRes::BDP_W_R_ValueFct}
         }
  
        \subfloat[\centering]{        
        \begin{minipage}{\textwidth}
            \begin{center}
                \hspace*{-0.6cm} \includegraphics[scale = 0.44]{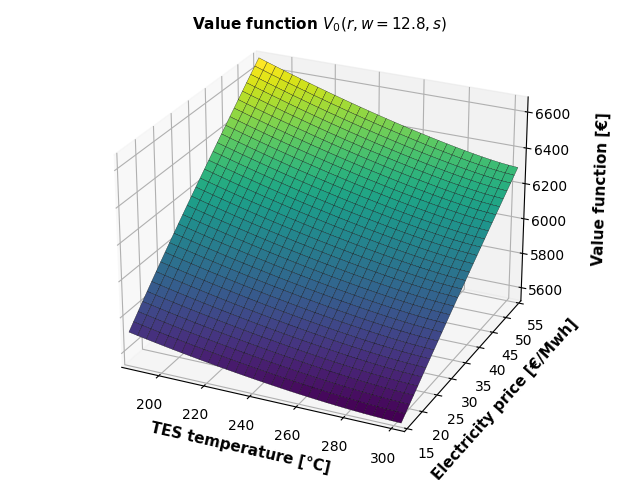}
			  \hspace*{-0.3cm} \includegraphics[scale = 0.44]{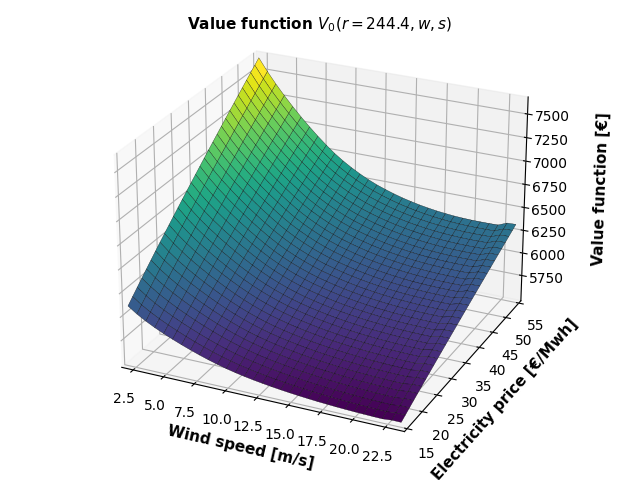}
            \end{center}
        \end{minipage} 
        \label{NumRes::BDP_WS_RS_ValueFct}}
		\caption{BDP: \textbf{(a)} Value function at initial time $n=0$ (left) and at terminal time $n=120$ (right) in terms of storage temperature and wind speed. \textbf{(b)} Value function at initial time $n=0$ depending on the storage temperature and electricity price (left) as well as on wind speed and electricity price (right).}
	\end{figure}

	Let us first discuss the results of the value function computed by BDP, which are presented in Figure \ref{NumRes::BDP_W_R_ValueFct}. The graphics on the left and right show the value function for the initial time step $n=0$ and the terminal costs $n=120$, respectively, as a function of TES temperature and wind speed. Obviously, low wind speed and low TES temperature lead to higher expected costs. While the dependence on the storage temperature is almost linear, the wind speed has a significant nonlinear influence on the value function, whereby the latter effect are induced by the WT power curve model. In addition, the dependence of the value function on storage temperature and electricity price as well as on wind speed and electricity price for the initial time is depicted in Figure \ref{NumRes::BDP_WS_RS_ValueFct} on the left and right, respectively. The electricity price affects the expected costs mostly linearly with respect to the TES temperature, with higher prices leading to higher values. Again, higher storage temperatures reduce costs and compensate for expensive grid electricity. The relation between wind speed and electricity price is also almost linear if we consider changes with respect to the prices. However, in terms of wind speed, the values again reflect the nonlinear power curve model.
	
\begin{figure}[H]
		\centering
        \subfloat[\centering BDP: Visualization of system dynamics and control for one trajectory.]{
        \begin{minipage}{\textwidth}
            \begin{center}
         \hspace*{-0.8cm} \includegraphics[scale = 0.46]{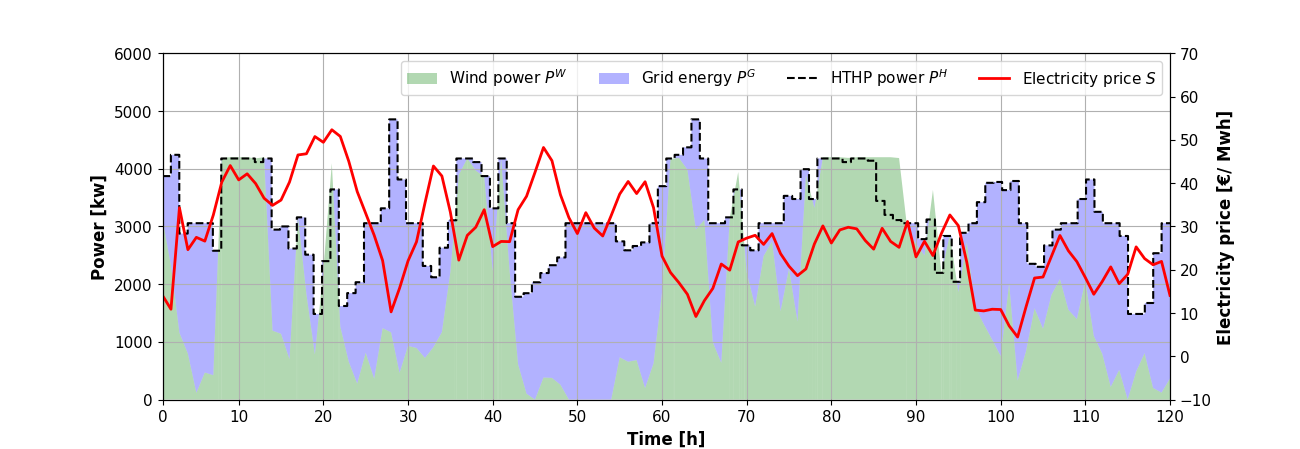} \\
            \hspace*{-0.8cm} \includegraphics[scale =0.46]{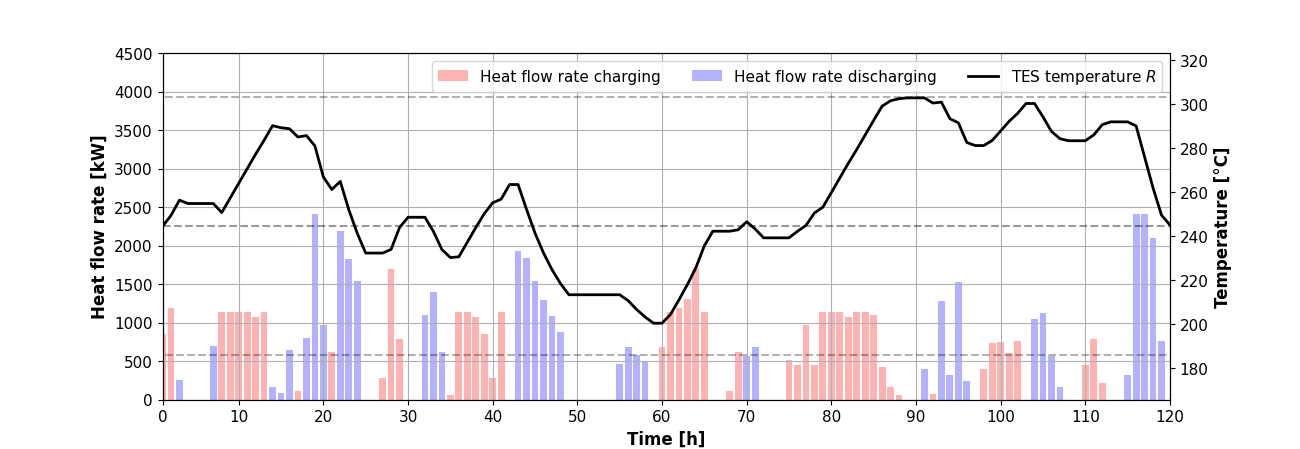}
            \end{center}
        \end{minipage}
        \label{NumRes::BDP_Trajectory}
        }

        \subfloat[\centering Q-learning: Visualization of system dynamics and control for one trajectory.]{
        \begin{minipage}{\textwidth}
            \begin{center}
         \hspace*{-0.8cm} \includegraphics[scale = 0.46]{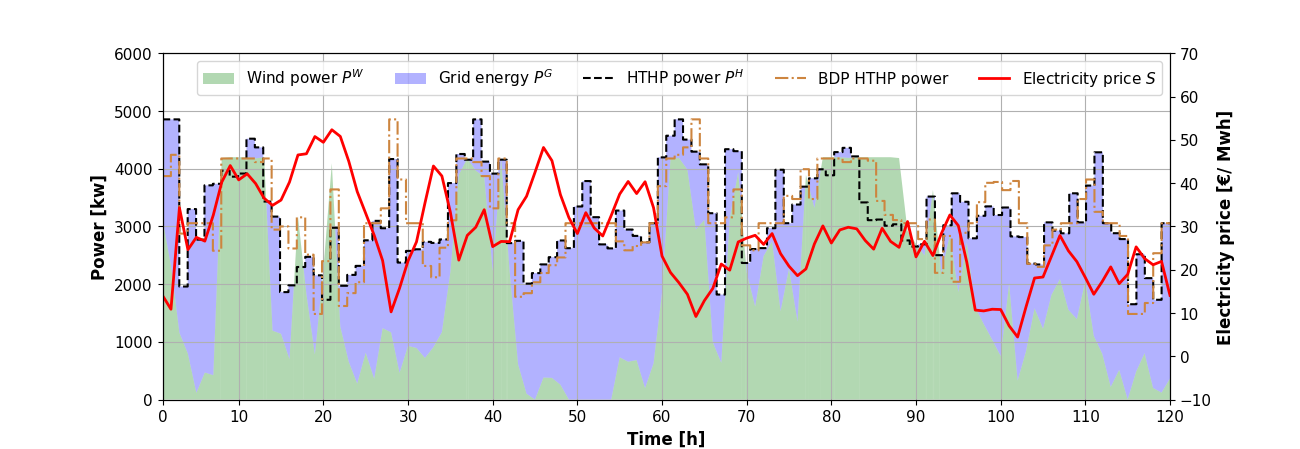} \\
            \hspace*{-0.8cm} \includegraphics[scale = 0.46]{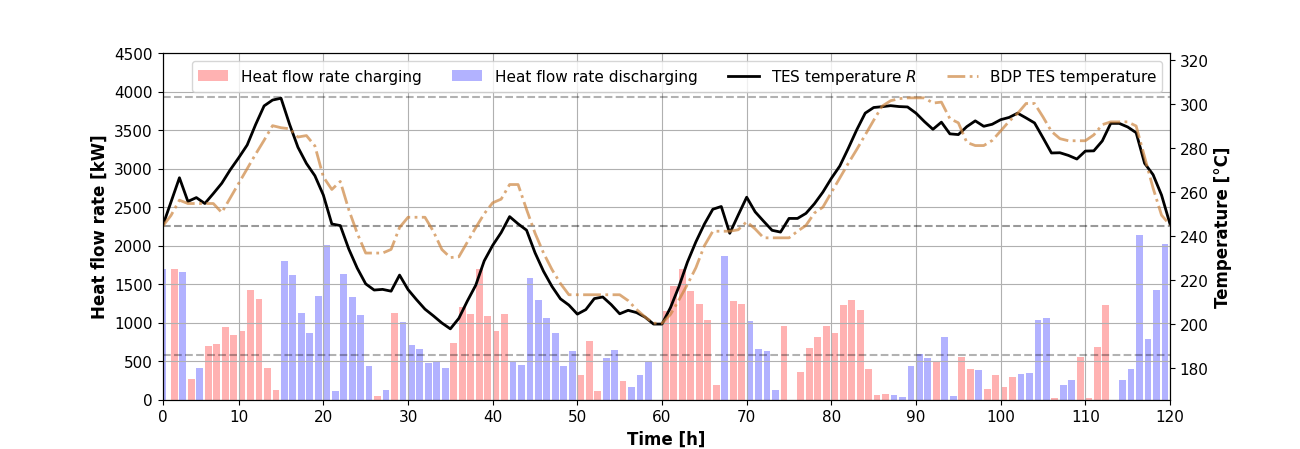}
            \end{center}
        \end{minipage}
        \label{NumRes::Q_Learning_Trajectory}
        }     
		\caption{The upper plot in \textbf{(a)} and \textbf{(b)} shows the electricity consumption (dotted black) for operation the HTHP, with generated wind energy (green) and consumed grid power (blue) stacked, as well as the electricity price (red). The respective lower plots visualize the average TES temperature (black) and the transferred heat flow rate during charging (red) and discharging (blue). For a better comparison of both methods, we include the HTHP electricity consumption and TES temperature (brown) from the BDP solution in \textbf{(a)} into \textbf{(b)}.}
	\end{figure}

	Second, we analyze a trajectory, starting from the initial state $(R_0,W_0,S_0) = (244.4,4,37)$ with $R_0 = r_\text{crit}$, for the control obtained by BDP as shown in Figure \ref{NumRes::BDP_Trajectory}. As expected, the control aims to charge the TES during periods of high wind power production and/or low prices. Charging that solely uses wind energy can be identified when the HTHP's electricity consumption (dashed black line) is covered by the available wind energy (green area), while additional power from the grid (blue area) is used to charge the TES when the price falls to small values.    
	
	For instance in hours $60$ to $65$, a combination of both scenarios can be observed. Conversely, when electricity prices are high and wind energy is not enough to cover the nominal HTHP power consumption, the TES is discharged, see e.g. hours 43 to 48. 
	Approximately four hours before the time horizon $\Hor$ is reached, the storage temperature is reduced to match the desired threshold $r_\text{crit}$.
    
	For further details on the optimal control $\AC$, see Figure \ref{BDP Control}, where the optimal decision rule is plotted as a function of time and TES temperature. At each time step $t_n = n \Delta t$, the decision rule is calculated for the corresponding TES temperature level (left y-axis) and the values of the seasonalities $\mu_W(t_n)$ and $\mu_S(t_n)$ at that time step. The seasonality functions $\mu_W^{}$ (dashed green line) and $\mu_S^{}$ (dotted red line) with their corresponding scales are also presented (right y-axis). The colored red and blue areas correspond to the charging and discharging mode of the system with the respective heat flow rate. The white areas represent the idle mode of the system, i.e., no charging or discharging is operated.     
	
	It can be seen that the control obtained by BDP captures the functional structure of the seasonalities. This means that when prices are globally low ($\mu_S^{}$ takes a global minimum), for example at the hours 25, 50 or 75, charging is the preferred action. Conversely, discharging is performed when prices are high and wind speeds are low, see hours 20 or 45. For At hours 15, 35 or 60, we can observe that when prices and wind speeds are locally high, the optimal action is to wait and operate in idle mode. Another situation where waiting is optimal, appears when the storage temperature is higher than the penalty temperature $r_\text{crit}$ (dashed gray line) and prices are at the global minimum of $\mu_S^{}$, see e.g. hour 45 to 55. However, when a maximum price is reached, the P2H system operates in discharge mode to compensate for the high electricity prices. Furthermore, when observing the seasonal patterns, it can be seen that the controls for charging, discharging and idle mode are usually centered around the local and global extrema of $\mu_S$. In particular, the electricity price appears to have a greater influence on the control than the wind speed (wind energy), as it has a major impact on the operational costs of the P2H system.    

    \newpage
    
	\begin{figure}[H]
		\centering
        \subfloat[\centering Visualization of the value function in terms of storage temperature and wind speed.]{
        \begin{minipage}[ht]{\textwidth}
        \begin{center}
            \begin{multicols}{2}
				\includegraphics[scale = 0.44]{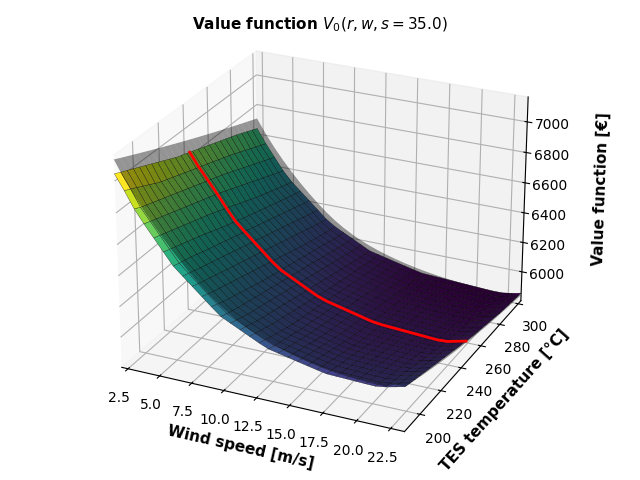} \includegraphics[scale = 0.44]{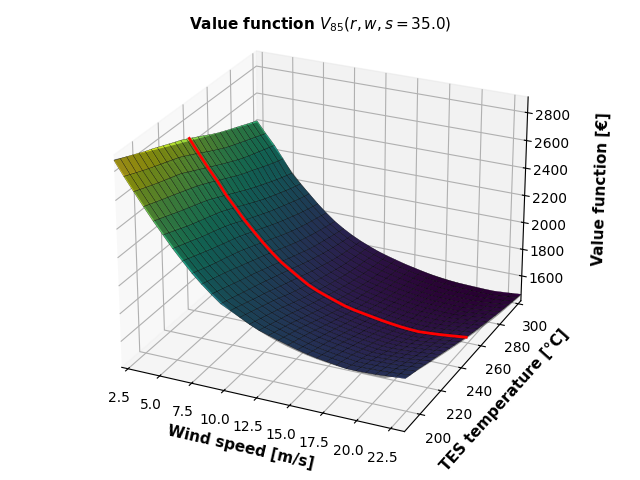}\\ 
				\includegraphics[scale = 0.44]{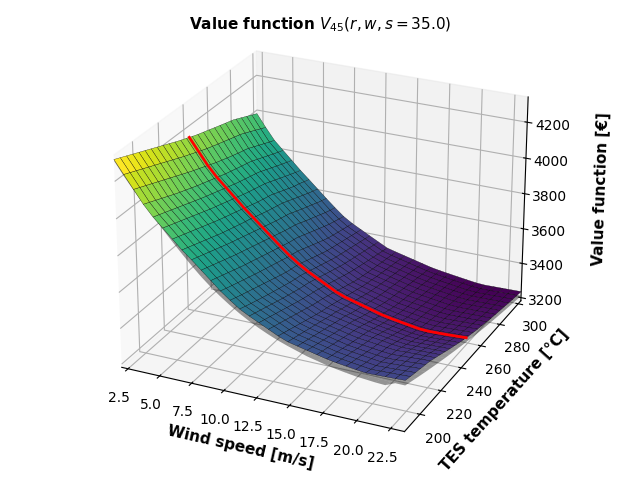}
				\includegraphics[scale = 0.44]{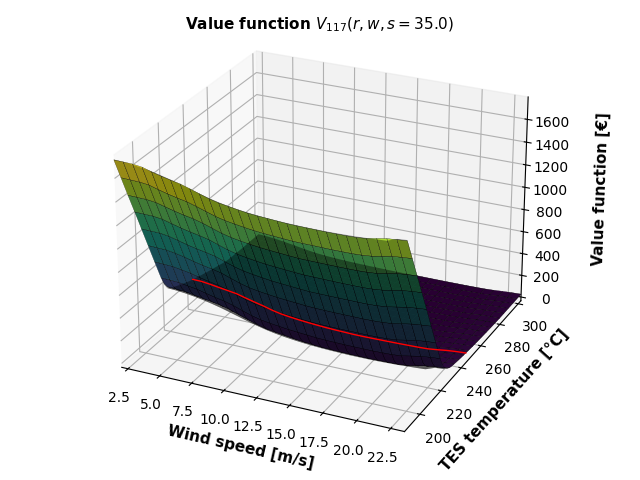}
			\end{multicols}
        \end{center}
		\end{minipage}
        \label{NumRes::Q_Solution_RW}
        }\\
		\subfloat[\centering Visualization of the value function in terms of electricity price and wind speed.]{
        \begin{minipage}[ht]{\textwidth}
        \begin{center}
			\begin{multicols}{2}
				\includegraphics[scale = 0.44]{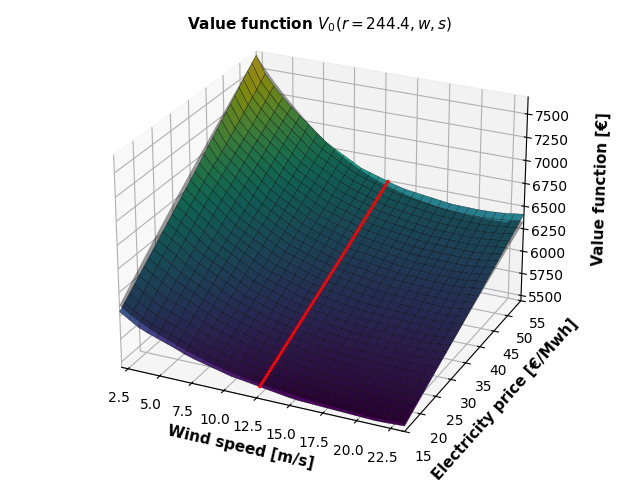} \includegraphics[scale = 0.44]{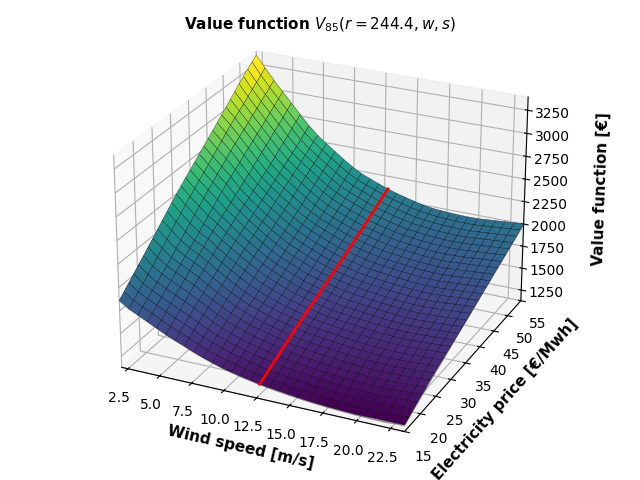}\\ 
				\includegraphics[scale = 0.44]{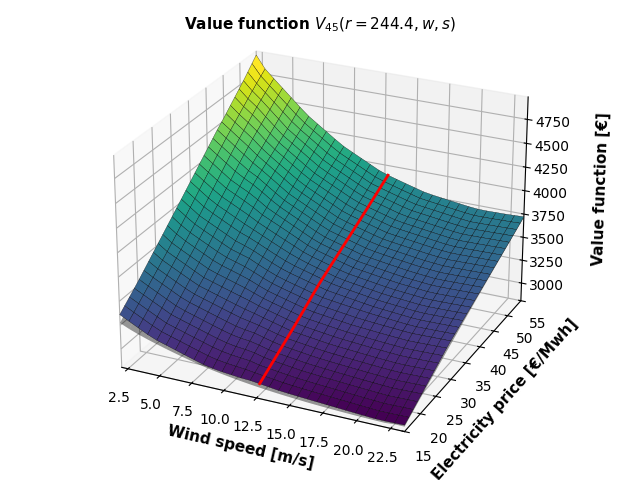}
				\includegraphics[scale = 0.44]{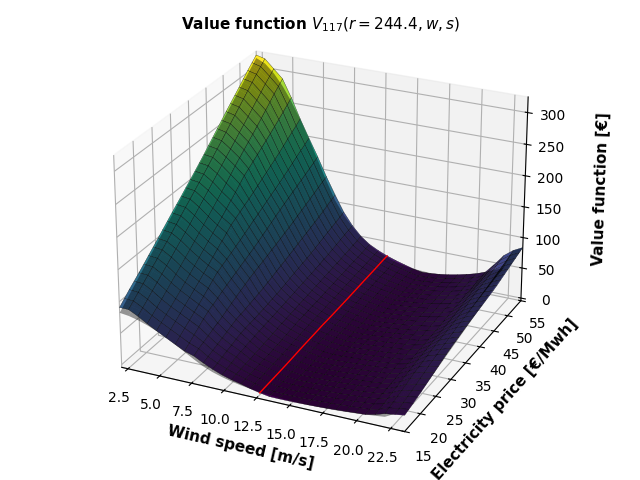}
			\end{multicols}
            \end{center}
		\end{minipage}
        \label{NumRes::Q_Solution_WS}
        }
        
		\caption{Q-learning: Visualization of the value function at times $n= 0,45,85,117$ depending on different state variables. The plot includes the BDP solution (gray) as a reference for comparison. The cross-sections (red) for each of the four value function plots are also shown in the Figure \ref{NumRes::Cross-Sections} for better visualization.}
\label{NumRes::Q_Solution}
	\end{figure}

	\subsection{Q-Learning} 
\begin{figure}[H]
		\hspace{-0.5cm} \includegraphics[scale = 0.48]{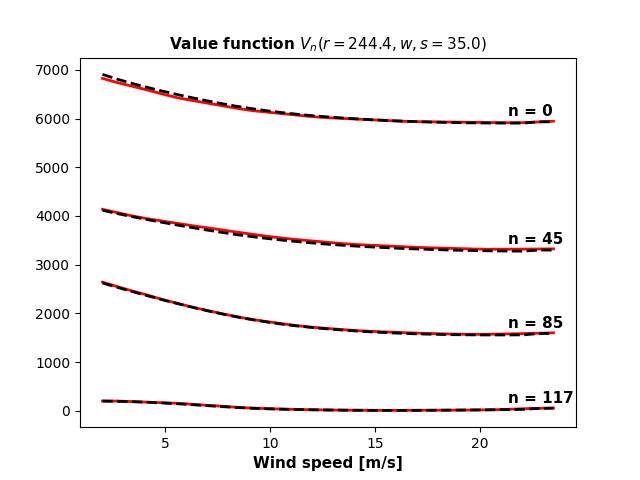}
		 \includegraphics[scale = 0.48]{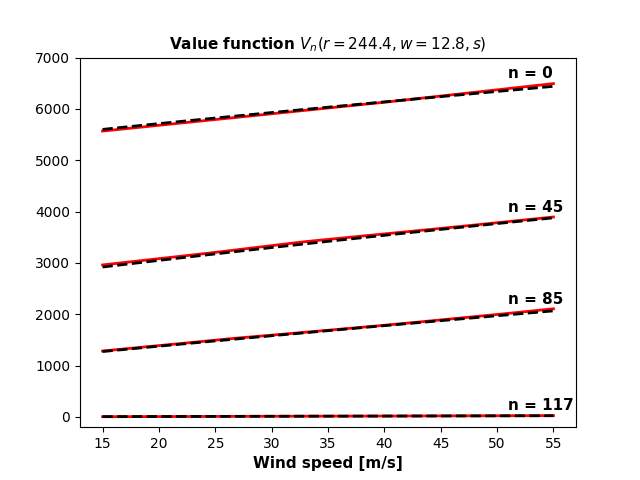}
		\caption{Visualization of the cross-sections from the value functions in Figure \ref{NumRes::Q_Solution_RW} (left) and Figure \ref{NumRes::Q_Solution_WS} (right). The black dashed curves show the value function approximation from BDP, which is compared with the value function from Q-learning given in red.}
		\label{NumRes::Cross-Sections}
	\end{figure}

	We now compare the results of Q-learning using replay buffer, summarized in Algorithm \ref{QlearningAlg}, and BDP based on the computed value function and the corresponding control. Thus, BDP serves as a benchmark, as we can expect accurate results due to the high computational effort. For the parameterization of $\overline{Q}(x,a,{\theta_n})$ a two-layer neural network is used with 128 neurons for each layer and ReLU activation functions. Since the parameterization is defined globally on the state space, no state discretization is required as for BDP. The state distribution of the initial states $X_0 = x_0$ is chosen as a uniform random distribution on the discretized state space $\widetilde{\mathcal{X}}_0$. The Figures \ref{NumRes::Q_Solution_RW} and \ref{NumRes::Q_Solution_WS} show the value function of both methods depending on TES temperature and wind speed as well as on wind speed and electricity price, respectively, for the hours 0, 45, 85 and 117. More precisely, in each value function plot of Q-learning, the reference solution of BDP is visualized as gray shaded graph. In addition, the red lines represent a cross-section of the value function, meaning that the value function is fixed in two variables and is visualized in order to provide a more detailed view at the error between both solutions. For a better comparison, the cross-sections (from Figure \ref{NumRes::Q_Solution}) are also depicted in Figure \ref{NumRes::Cross-Sections}. Compared to BDP, Q-learning is able to capture the same shape of the value function. Especially, for $n$ near the terminal time, both function approximations differ only slightly. The further we move back in time from the terminal time, the more differences become visible.
    
\begin{figure}[H]
    \vspace*{-1cm}
		\raggedleft 
           \subfloat[\centering DP: Visualization of the decision rule $u_n$.]{\begin{minipage}{\textwidth}
            \begin{center}
         \hspace*{-0.8cm} \includegraphics[scale = 0.36]{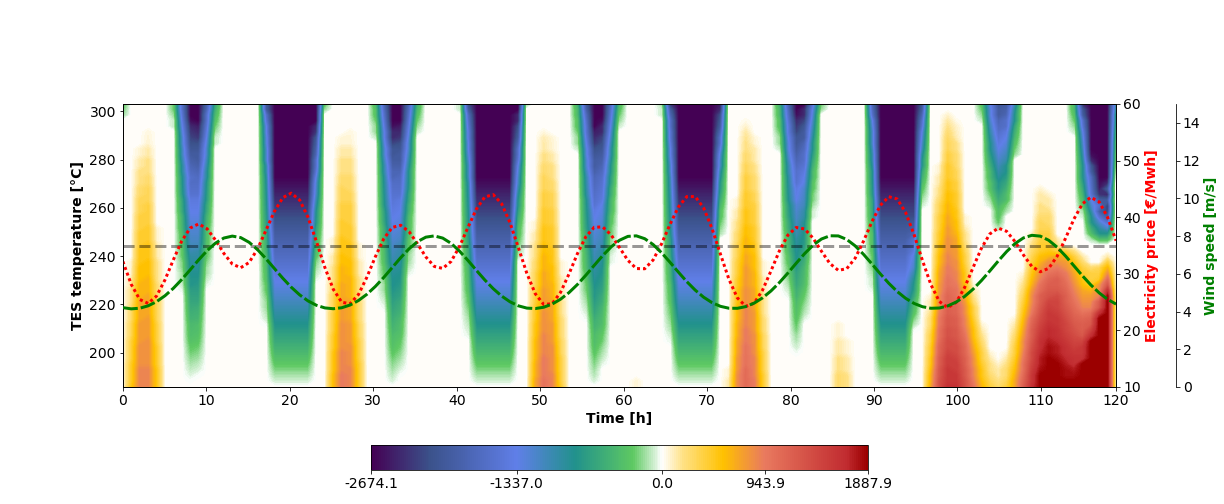}
            \end{center}
        \end{minipage} 
        \label{BDP Control}}
        \vspace*{-1cm}
                \subfloat[\centering Q-learning: Visualization of the decision rule $u_n$.]{
                \begin{minipage}{\textwidth}
            \begin{center}
         \hspace*{-0.8cm} \includegraphics[scale = 0.36]{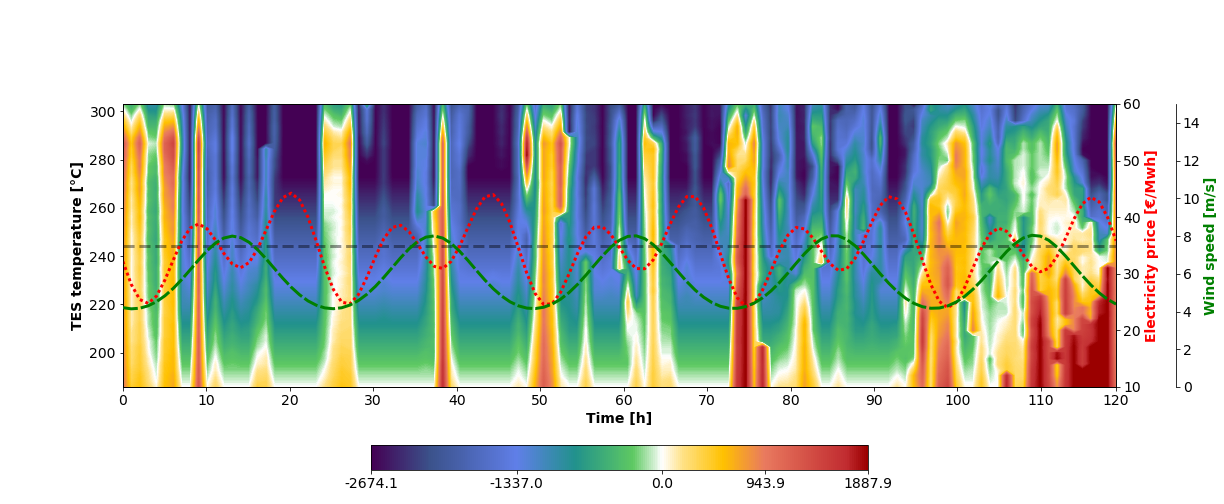} \
            \end{center}
        \end{minipage}
        \label{Q_Learning Control}
                }
		\caption{
			Visualization of the optimal decision rule $u_n$ with respect to the TES temperature together with threshold temperature $r_\text{crit}$ (dashed gray line). At each time point $t_n = n \Delta t$, the decisions calculated for the TES temperature is given in terms of the values of the seasonalities $\mu_W^{}(t_n)$ (dashed green line) and $\mu_S^{}(t_n)$ (dotted red line), i.e. $u_n(r,w=\mu_W(t_n),s=\mu_S(t_n)) = \aC$. The red color represents the HTHX outlet temperature, the blue color the HTHX inlet temperature, while white areas specify idle mode of the system. 
			The corresponding opacity indicates the absolute temperature while charging or discharging.}
           
	\end{figure}

	Figure \ref{NumRes::Q_Learning_Trajectory} also confirms that the approximation of the value function by Q-learning is similar to that of BDP. The difference in the control is mainly reflected in the charging and discharging intensity of the TES and thus affects the HTHP's electricity consumption. Apart from this, the controller aims to charge the TES during times of high wind energy availability or/and low prices and to discharge vice versa. Overall, the comparison shows that the Q-learning control is qualitatively similar to that of BDP. However, an exception is the control in Figure \ref{NumRes::Q_Learning_Trajectory}, which almost does not contain any waiting periods. Instead, the charging or discharging periods are generally longer compared to the BDP control, which can be seen in hours 15 to 25. Furthermore, Figure \ref{Q_Learning Control} provides a more detailed look at the optimal decision rule obtained with Q-learning. Obviously, the charging mode is performed in regions where price seasonality has a global minimum. However, almost every time a peak occurs, the high prices are compensated by discharging the TES. Again, the role of wind speed appears to be less important than the influence of the electricity price on the control itself. Even though seasonality is taken into account, the overall structure and sequence are not captured  as well as with the BDP.

	\myparagraph{Computational Time} 
	In addition to the qualitative comparison of the numerical results with BDP and Q-learning, the computational effort required to compute the numerical solution is also of practical importance. Here, the computational time serves as an indicator of how well the methods are able to deal with the curse of dimensionality. All computations are performed on a compute server with 320 GB RAM running two Intel Xeon Gold 6136 processors, each with 12 cores and 24 threads. BDP calculates and saves the value function for the grid points of the discretized state spaces $\widetilde{\mathcal X}_n$. To speed up the computations of the value function in each time step, we will use all available cores and calculate its values for different grid points in parallel. In total, BDP requires 36 hours computational time, which corresponds to approximately 18 minutes per time step. In contrast, Q-learning only takes around 8 hours on a single core to compute the approximate solution. This means a time saving of a factor of 4 for the problem considered with three-dimensional states and one-dimensional actions. For stochastic optimal control problems in higher dimensions, we can expect much greater savings.
	
	\section{Conclusion} \noindent
	This work presents a mathematical model for the cost-optimal operation of an industrial P2H system. Apart from providing a modeling approach for the stochastic processes that takes into account correlation between wind speed and electricity price, we also calibrated the associated parameters with real-world data (see Supplementary Material \ref{SM:calibration}). The resulting discrete-time stochastic optimal control problem is formulated as an MDP and solved using the classical dynamic programming approach as well as modern reinforcement learning techniques, namely Q-learning. 
	
	A comparison of the numerical results shows that both methods can achieve similar approximations of the value function and yields reliable cost-optimal decision rules. Although the results of Q-learning differ in some aspects, it offers a faster and computationally more efficient solution for complex control problems. This is especially useful for problems with high-dimensional state and control spaces, where the dynamic programming approach will fail due to the curse of dimensionality.
	
	By dropping the assumptions of constant mass flow $\mdot$ and waste heat temperature $\Tc$, we can extend our model and make it more general. In this case, it is necessary to introduce $\Tc$ as an additional state and $\mdot$ as an additional control variable, which increases the dimension of the control problem. Even though the classical backward recursion of dynamic programming might become intractable for this extended model due to the curse of dimensionality, we are confident that Q-learning still offers an efficient solution. However, as the dimension of the action space grows, it becomes infeasible to calculate the minimal action by discretization. Appropriate gradient descent methods could be used, at the cost that this may slow down the algorithm. Reinforcement learning algorithms such as policy gradients or actor-critic methods as in \cite{Sutton2018} might be more suitable for dealing with large (continuous) action spaces, as they already include a way to handle this minimization step. 

	\clearpage
	\appendix
	\section*{Appendix}
	\section{Nomenclature}
    \vspace{-2cm}
    	\begin{table}[H]   
		\begin{minipage}{0.6\textwidth}
			\begin{tabular}{p{1cm}p{6.5cm}}
				\multicolumn{2}{l}{\textbf{Acronyms}}\\[1ex]
				HTHP&High-temperature heat pump\\
				{HTF}&{Heat transfer fluid}\\
				{TES}&{Thermal energy storage}\\
				{LTHX}&{Low-temperature heat exchanger}\\
				{HTHX}&{High-temperature heat exchanger}\\
				{WT}&{Wind turbine}\\
				{SG}&{Steam generator}\\
				{P2H}&{Power-to-Heat}\\
				{MDP}&{Markov desicion process}\\
				{SDE}&{Stochastic differential equation}\\
				{BDP}&{Backward dynamic programming}\\
				{SGD}&{Stochastic gradient descent}\\[2ex]

				\multicolumn{2}{l}{\textbf{Latin symbols}}\\[1ex]
				{$\Tc$}&{LTHX inlet temperature}\\
				{$\TSGin$}&{SG inlet temperature}\\
				{$\TSGout$}&{SG outlet temperature}\\
				{$\THTout$}&{HTHX outlet temperature}\\
				{$\THTin$}&{HTHX inlet temperature}\\
				{$\Tch, \Tdch$}&{TES outlet temp. charging/discharging}\\
				{$\lC, \lD$}&{Charging/discharging factor}\\
				{$\mdot$}&{Mass flow}\\
				{$m_\text{s}$}&{Mass of TES}\\
                {$c_\text{p,s},c_\text{p,f}$}&{Heat capacity of TES/thermal oil} \\
                {$d$}&{Rotational speed}\\
                {$\PG$}&{Electrical power of grid}\\
				{$\PH$}&{Electrical power of HTHP}\\
				{$\PW$}&{Electrical power of WT}\\
                {$\AC$}&{Heat flow rate}\\
				{$B$}&{Brownian motion}\\
			\end{tabular}
		\end{minipage}
		\hspace*{-1cm}\begin{minipage}{0.4\textwidth}
			\begin{tabular}{p{1cm}p{5cm}}
				& \\[2ex]
				\multicolumn{2}{l}{\textbf{}}\\[1ex]
				& \\
                {$X$}&{State process}\\{$Y$}&{Ornstein-Uhlenbeck process}\\
				{$R$}&{Storage temperature}\\
				{$S$}&{Electricity grid price}\\
				{$W$}&{Wind speed}\\
                {$\noise$}&{Random disturbance}\\
                {$\widehat{\noise}$}&{Quantizer for $\noise$}\\
				{$t$}&{Time}\\
				{$\Hor$}&{Terminal time horizon}\\
                {$\Delta t$}&{Time step size}\\
				{$n$}&{Time point}\\
				{$N$}&{Total number of time points}\\
                {$\npump$}&{Number of HTHPs running in parallel}\\
				{$\mathcal{X}$}&{State space}\\
                {$x$} & State variable\\
				{$\mathcal{A}$}&{Action space}\\
                {$\mathcal{T}$}&{Transition operator}\\
                {$a$} &{Control/action variable}\\
                {$C_n$} &{Running cost}\\
                {$G$} &{Terminal cost}\\
                {$J_n$} &{Performance criterion}\\
                {$V_n$} &{Value function}\\
                {$Q_n$} &{State action function}\\
                {$\theta$} &{Parameter vector}\\[2ex]
				
				\multicolumn{2}{l}{\textbf{Greek symbols}}\\[1ex]
				
				{$\mu$}&{Mean reversion level}\\
				{$\lambda$}&{Mean reversion speed}\\
				{$\sigma$}&{Volatility}\\
				$\eta$ &  Spread\\
				{}&{}\\
				& \\ 
				&  
			\end{tabular}
		\end{minipage}
	\end{table}

    \newpage
		
	\section{Model Parameters} 
	\label{Appendix:TableOfParameters}
	\begin{minipage}{\textwidth}
    The parameters of the model and the methods used to solve the MDP in our numerical experiments, i.e., BDP and Q-learning, are listed in Tables \ref{ParametersModel} and \ref{NumericalResultsParameters}. 
	    \begin{table}[H]
		\caption{Parameters used in the description and application of the industrial P2H system, the stochastic processes and the numerical methods.}
		\centering
		\subfloat[\centering Parameters used in the overall industrial P2H system and the stochastic processes.]{
			\label{ParametersModel}
			\begin{tabular}{p{1.4cm}p{7.3cm}p{1.2cm}p{2.1cm}}
				\midrule
				\multicolumn{2}{l}{\textbf{Parameters of the P2H System}} & Value & Unit \hfill\\
				\midrule
				$\mdot$& Thermal oil massflow &  6 & [kg/s]\\
				$\dmin$& Minimal compressor shaft speed &  0.8 & [-]\\
				$\dmax$& Maximal compressor shaft speed &  1.53 & [-]\\
				$m_\text{s}^{}$& Storage mass &  600000& [kg]\\
				$c_\text{p,s}^{}$& Storage specific heat capacity &  1.025 & [kJ/kgK]\\
				$c_\text{p,f}^{}$& Thermal oil specific heat capacity &  2.314 & [kJ/kgK]\\
				$\Tc$& Waste heat air temperature &  80 &[$^\circ$C]\\
				$\TSGout$& SG outlet temperature  &  185.8 & [°C]\\
				$\TSGin$& SG inlet temperature  &  303.0 & [°C]\\
				$\tHTinmax$& Maximum HTHX oil inlet temperature  & 250 &[$^\circ$C]\\
				$\npump$& Number of HTHPs &  3 & [-]\\
				$\epsilon^\text{C}$& Charging efficiency &  0.9 & [-]\\
				$\epsilon^\text{D}$& Discharging efficiency &  0.9 & [-]\\
				\midrule
				\multicolumn{2}{l}{\textbf{Estimated Parameters for the Ornstein-Uhlenbeck Processes}} &  Value & Unit\\ 
				\midrule
				${\lambda}_W$ & Estimate for mean reversion speed of $Y^W$ &0.1702& [1/$s$]\\
				${\sigma}_W$ & Estimate for volatility of $Y^W$ &0.2486& [m/$\sqrt{\text{s}^3}$]\\
				${c}_W^{}$ & Estimate for correlation/conversion constant &0.5483& [(\euro  s)/(MWh\,m)]\\
				${\lambda}_S$ & Estimate for mean reversion speed of $Y^S$&0.2534& [1/$s$]\\
				${\sigma}_S$ &Estimate for volatility of $Y^S$ &0.1072& [\euro /(MWh\,$\sqrt{\text{s}}$)]\\
				\midrule
				\\
			\end{tabular}
		} \hfill
        \vspace{-0.3cm}
		\subfloat[\centering Parameters for the numerical solution of the MDP with backward dynamic programming and Q-learning.]{
			\label{NumericalResultsParameters}
			\begin{tabular}{p{1.4cm}p{7.3cm}p{1.2cm}p{2.1cm}}
				\midrule
				\multicolumn{2}{l}{\textbf{Parameters for BDP and Q-Learning}} &  Value & Unit \\
				\midrule
				$n_R^{},n_W^{},n_S^{}$ & Number of grid-points in each state direction &$51$ & [-]\\
				$n_A^{}$ & Number of grid-points for control $\AC$&$31$ & [-]\\
				$L$ & Number of quantization points & $400$ & [-]\\
				$k^\text{Ref}$ & $\sigma$-rule parameter for $\mathcal{X}^\text{Ref}_n $&$3$ & [-]\\
				$k^\text{Ext}$ & $\sigma$-rule parameter for $\mathcal{X}_n $&$4$ & [-]\\
				$k^\text{max}$ & Number of iteration for Q-learning & $50000$& [-] \\
				$M$ & Batch size in Q-learning & $128$ & [-]\\
				$\mathcal{N}_R$ & Size of replay buffer per time step & $20000$ & [-]\\
				$\varepsilon_0$ & Exploration rate linearly decaying to zero & $1$ & [-]\\
				$\alpha$ & Learning rate & $0.001$& [-]\\
				\midrule
				\multicolumn{2}{l}{\textbf{Other}} &  Value & Unit\\
				\midrule
				$r_\text{crit}$ & Critical terminal temperature & $244.4$& [°C] \\
				$s_\text{Pen}$ & Penalization price & $90$ & [\euro/Mwh]\\
				$s_\text{Liq}$ & Liquidation price & $0$ & [\euro/Mwh]\\
                $\zeta$ & Selling model parameter & $0$ & [-]\\
                $t_\text{E}$ & Time horizon & $120$ & [h]\\
                $N$ & Time horizon partition & $120$ & [-]\\
                $\Delta t$ & Time step size & $1$ & [h]\\
				\midrule
			\end{tabular}
            \vspace{-3cm}
		}
	\end{table}
	\end{minipage}

    \newpage

    \section{Surrogate Models}\noindent
	\label{app:SurrogateModels}
	As described in \cite{Walden2023}, the part-load behavior of the HTHP and the SG is represented by polynomial surrogate models, which are created using process simulation software. 
    
    The HTHP surrogate models are given by the multivariate polynomials $F_1$ and $F_2$. Here, $F_1$ describes the HTHX outlet temperature $\THTout$ and $F_2$ the consumed electrical power $\PH$ at a given HTHX inlet temperature $\THTin= \tau^{H}$, mass flow $\dot m$, LTHX air inlet temperature $\Tc = \tau^{L}$ and compressor shaft speed $d$ by
	\begin{align}
		\begin{split}
			\THTout&  = F_1(\tau^{H}, \dot{m}, \tau^{L},d)
			=95.9612+0.93433 \tau^{H}-0.327753\dot{m}+0.0146542 \tau^{L}\\
			&\hspace{0.3cm} -271.354 d+0.00104853 (\tau^{H})^2+0.0211819 \tau^{H}\dot{m}-0.706122 \tau^{H} d\\
			&\hspace{0.3cm}-0.00388073\dot{m} \tau^{L}+0.0595068 \tau^{L} d-29.4801\dot{m} d+1.04924\dot{m}^2 +562.428 d^2\\
			&\hspace{0.3cm}-0.000716825 (\tau^{H})^2 d -2.18172\dot{m} d^2-151.476 d^3+0.0229386 \tau^{H}\dot{m} d\\
            &\hspace{0.3cm}+0.881391\dot{m}^2 d +0.203578 \tau^{H} d^2-0.0405702\dot{m}^3 -0.00148575 \tau^{H}\dot{m}^2,\\[2ex]
			\PH & = F_2(\tau^{H}, \dot{m}, \tau^{L}, d)  
			= 127.87+2.06342 \tau^{H}+2.55723\dot{m}+0.756419 \tau^{L} \\
			&\hspace{0.3cm}-1164.84 d-1.3829 \tau^{L} d-0.0168942 \tau^{H}\dot{m}-2.60579 \tau^{H} d-0.540713\dot{m}^2 \\
            &\hspace{0.3cm}+13.3204\dot{m} d+1556.66 d^2,
		\end{split}
        \label{Appendix:F1_F2}
	\end{align}
	with $\tau^{H}\in[177,250], \dot{m}\in[5,16], \tau^{L} \in[60,100]$ and $d\in[0.8,1.53]$. It should be noted that the HTHP's surrogate models are derived on the basis of a single HTHP. Consequently, for the investigated P2H system, the electricity consumption $\PH$ (cf. \eqref{PowerBalance}) and the mass flows in the thermal oil loop must be multiplied by the number of parallel operating HTHPs $\npump$. 
	
	The surrogate models for SG inlet and outlet temperatures as functions of mass flow read as
	\begin{equation}
		\TSGin = \FSGin(\dot{m}) = 201.92+\frac{1819.32}{\npump\dot{m}},\quad \TSGout = \FSGout(\dot{m}) = 196.3-\frac{188.4}{\npump\dot{m}}.\label{Appendix:SG_functions}
	\end{equation}
	As mentioned, the latter equations take into account the fact that the mass flows $\mdot$ at the HTHX outlet of the parallel running HTHP's are merged in the thermal oil loop.
	
	\section{Details on Thermal Energy Storage Operational Modes} \noindent
	\label{app:OperationalModes}
	The TES operational modes introduced in Subsection \ref{sec:OpConstrTES} are described in more detail below.
    
	\myparagraph{Charging Mode} 
	During charging, the TES charging outlet temperature $\Tch$ is given by a simplified model defined by a weighted average as 
	\begin{align}
		\Tch(t) = (1 -\epsC) \THTout_n + \epsC R(t),
		\label{Tch}
	\end{align}
	where a constant charging efficiency $\epsC\in[0,1]$ is assumed. An efficiency $\epsC<1$ indicates that the TES does not completely absorb the thermal energy emitted by the thermal oil. The two oil streams that are routed through and bypass the TES with a certain proportion using $\lC$ are mixed and enter the inlet of the SG. We assume an ideal heat transfer so that the mixing temperature is reached immediately, which leads to the following relationship for the constant SG inlet temperature
	\begin{align}
		\TSGin = (1-\lC(t)) \THTout_n + \lC(t) \Tch(t).
		\label{TSGin}
	\end{align}
	Due to $\THTout_n > \TSGin$ and $\lC\in(0,1]$ it follows directly $\Tch(t)\le \TSGin$, so that the charging outlet temperature never exceeds the SG inlet temperature. By charging the TES, its temperature $R$ increases linearly over time (see \eqref{storage_dynamic} and \eqref{Control_inlet_outlet0}) and based on \eqref{Tch}, $\Tch$ also changes with time. For this reason, the charging factor must be adjusted over time to ensure a constant SG inlet temperature. Combining of \eqref{HTHX_inlet_outlet} with \eqref{Tch} and \eqref{TSGin} leads to the following relation for the charging factor
	\begin{align}
		\lC(t) = \frac{\tHTout(\AC_n) - \TSGin}{\epsC (\tHTout(\AC_n) - R(t))}, 
		\label{charge_factor}
	\end{align}
	which must hold in each period $[t_n,t_{n+1})$ and describes the dependence of $\lC$ on the chosen control $\AC_n$ and the TES temperature $R$.
	
	\myparagraph{Discharging Mode}
	Analogously to \eqref{Tch}, the discharging outlet temperature $\Tdch$ is given by 
	\begin{align}
		\Tdch(t) = (1 -\epsD) \TSGout + \epsD R(t),
		\label{Tdch}
	\end{align}
	where the discharging efficiency $\epsD$ is assumed to be a constant in $[0,1]$. For $\epsD<1$, the heat transfer from the TES to the thermal oil is not perfect, as the oil leaves the TES with temperature $\Tdch<R$. Discharging is determined by $\lD$, which describes the proportion of the HTF passing through the TES before entering the HTHX inlet, meaning that 
	\begin{align}
		\THTin_n =  (1 - \lD(t)) \TSGout + \lD(t) \Tdch(t).
		\label{HTHXin}
	\end{align}
	As above, the discharging factor must also be continuously adjusted in each period to ensure the constant HTHX inlet temperature $\THTin_n$. 
	Substituting \eqref{Tdch} into \eqref{HTHXin} and using \eqref{HTHX_inlet_outlet}, yields for $t\in [t_n,t_{n+1})$
	\begin{align}
		\lD(t) = \frac{\tHTin(\AC_n) - \TSGout}{\epsD ( R(t) - \TSGout)}.
		\label{discharge_factor}
	\end{align}
	Note that the time-varying charging and discharging factors \eqref{charge_factor} and \eqref{discharge_factor} depending on $R$ and chosen $\AC_n$ are used to derive the control conditions in Section \ref{app:ControlConstraints}.

	\section{Details on Control Constraints}\noindent
	\label{app:ControlConstraints}
	Here, we provide the derivation of the control constraints given in Subsection \ref{sec:ControlConstraints}.
    
	\myparagraph{Upper Bound} 
	A first upper bound for values of the control process $\AC$ is obtained from a natural upper bound of the HTHX outlet temperature $\THTout_n$. This follows from relation \eqref{Tout} and its monotonicity with respect to the compressor shaft speed. Replacing $D_n$ by the maximum shaft speed $\dmax$ of the HTHP and knowing that $\THTin_n= \TSGout$, during the charging process, it holds $\THTout_n \leq \tHToutmax:= F_1(\TSGout,\mdot, \Tc , \dmax)$ and relation \eqref{Control_inlet_outlet} implies 
	\begin{align}
		\AC_n \leq \amax^1:= \aT(\TSGout, \tHToutmax).
		\label{aO_constraint_1.1}
	\end{align}	
	
	A second upper bound follows from the fact that the charging factor $\lC$  is restricted to values in $[0,1]$. During charging operation, the average TES temperature increases, i.e., $R_n \leq R(t)$ for $t \in [t_n,t_{n+1}]$. Since $\THTout_n$ and $\TSGin$ are constants, \eqref{charge_factor} implies that the charging factor $\lC$ also increases. Thus, to ensure that $\lC\in[0,1]$ on $[t_n,t_{n+1}]$, it is sufficient to require $\lC(t_{n+1}) \leq 1$. This leads to the following upper bound for $\AC$ by inserting \eqref{storage_dynamic}-\eqref{HTHX_inlet_outlet} into \eqref{charge_factor}
	\begin{align}
		\AC_n \leq   \amax^2(R_n) \quad\text{with}\quad 
		\amax^2(r)=\aT\bigg(\TSGout, \TSGin + \frac{\epsC (\TSGin - r)}{1- \epsC(1-\zeta)}\bigg), \quad r \in [\rmin, \rmax], 
		\label{aO_constraint_2}
	\end{align}
	where $\zeta = \npump\Delta t\tfrac{\mdot c_{F}^{}}{m_{S}^{} c_{S}^{}}$, which depends on the TES temperature $R_n$. In addition, the condition $R(t) \leq r_\text{max} =  \TSGin$ on $[0,\Hor]$ must also be fulfilled. During charging, $R$ increases, so it is sufficient to require $R_{n+1} = R(t_{n+1}) \leq \TSGin$.  However, if we substitute $\amax^2$ into $\eqref{storage_dynamic}$, then this condition is satisfied and no additional constraint is applied. In summary, we obtain the following state-dependent upper bound of the control 
	\begin{align}
		\amax(R_n) = \min\{\amax^1,\amax^2(R_n) \}. 
		\label{upper_aO_contraint}
	\end{align}
	
	\myparagraph{Lower Bound}
	A first lower bound results from technical reasons, namely that the HTF temperature entering the heat pump must be bounded from above, by $\tHTinmax$. During discharging it holds $\THTout_n= \TSGin$ so that we obtain from relation \eqref{Control_inlet_outlet} the following:
	\begin{align}
		\AC_n \geq \amin^1:= \aT(\tHTinmax, \TSGin).
		\label{aO_constraint_1.2}
	\end{align}	
	
	A second lower bound follows, similar to the upper bound $\amax^2$, from the fact that the discharging factor $\lD$ only takes values in $[0,1]$. During discharging, the TES temperature decreases with time, i.e., $R(t_{n+1}) \leq R(t)$ for $t \in [t_n,t_{n+1}]$. Since $\THTin_n$ and $\TSGout$ are constants, \eqref{discharge_factor} implies that the discharging factor $\lD$ increases. To ensure that $\lD\in[0,1]$ on $[t_n,t_{n+1}]$, it is therefore sufficient to require $\lD(t_{n+1}) \leq 1$. Inserting \eqref{storage_dynamic}-\eqref{HTHX_inlet_outlet} into \eqref{discharge_factor}, the lower bound depending on $R_n$ is obtained as follows
	\begin{align}
		\label{aI_constraint_2}
		\AC_n \geq  \amin^2(R_n) \quad\text{with}\quad 	
		\amin^2(r) = \aT\bigg(\TSGout + \frac{\epsD(r -\TSGout)}{1 + \zeta \epsD},\TSGin\bigg), \quad r \in [\rmin, \rmax].
	\end{align} 
	
	Finally, the condition $R(t) \geq r_\text{min} =  \TSGout$ must be fulfilled on $t \in [0,\Hor]$. In discharging mode, $R$ is decreasing, so it is sufficient to require  $R_{n+1} = R(t_{n+1}) \geq \TSGout$.  Inserting $\amin^2$ into \eqref{storage_dynamic} implies this condition. To summarize, the following state-dependent lower bound of the control can be specified by
	\begin{align}
		\amin(R_n) = \max \{ \amin ^1,\amin ^2(R_n)\}.
		\label{upper_aI_contraint}
	\end{align}
	
	\section{Proofs}
	\label{Appendix:Section:Proofs}
	\subsection{Proof of Lemma \ref{Lemma: E of Y^S, Y^W}} 
	\label{Paragraph_Solutions_Processes}
	\myparagraph{Solution, Expected Value and Variance for $Y^W(t)$}
	First we consider the Ornstein-Uhlenbeck process $Y^W(t)$ and find a solution on $[\tl, \tu]$ with initial condition $Y^W(\tl)= y_W^{} \in \mathbb{R}$. Applying Itô’s 
	formula with $F(t,Y^W(t)) = e^{\lambda_W t} Y^W(t)$ yields
	\begin{align}
		\mathrm{d} F(t,Y^W(t)) = \lambda_W e^{\lambda_Wt} Y^W(t) \mathrm{d}t + e^{\lambda_Wt} \mathrm{d} Y^W(t) = e^{\lambda_W t} \sigma_W^{} \mathrm{d} B^W(t).
	\end{align}
	Using the integral representation and multiplying both sides with $e^{-\lambda_W t}$ leads to the solution
	\begin{align}
		Y^W(t) = e^{-\lambda_W (t-\tl)} y_W^{} +  I_W(t),
		\label{solution_Y_W}
	\end{align}
	with $I_W(t) = \int_{\tl}^t \sigma_W^{} e^{-\lambda_W(t-r)}\mathrm{d} B^W(r)$.
	The stochastic integral appearing in \eqref{solution_Y_W} is a mean-zero martingale. Using this fact, the conditional expected value of $Y^W(t)$ given $Y^W(\tl)= y_W^{}$ reads 
	\begin{align}
		m_{Y^W}(\tau,y_W^{}) =  \mathbb{E}[Y^W(t)] = y_W^{} e^{-\lambda_W \tau},
	\end{align}
	with $\tau = t-\tl$.
	The variance follows from the Itô isometry, leading to 
	\begin{align}
		\Sigma_{Y^W}^2(\tau) = \textrm{Var}(Y^W(t)) & = \mathbb{E}[ I_W(t)^2 ]  
		= \mathbb{E}\bigg[  \int_{\tl}^t \sigma_W^2 e^{-2 \lambda_W(t-r)}\mathrm{d} r  \bigg]= \frac{\sigma_W^2}{2 \lambda_W} ( 1- e^{-2\lambda_W\tau}).
		\label{App:: Var_Y_W}
	\end{align}
	
	\myparagraph{Solution, Expected Value and Variance for $Y^S(t)$}
	To obtain a solution for $Y^S(t)$, Itô’s formula is applied with $F(t,Y^S(t)) = e^{\lambda_S t} Y^S(t)$ that gives
	\begin{align}
		\mathrm{d} F(t,Y^S(t)) = \lambda_S e^{\lambda_S t} Y^S(t) \mathrm{d} t + e^{\lambda_S t} \mathrm{d} Y^S(t) = -\lambda_S c_W^{} e^{\lambda_St} Y^W(t) \mathrm{d} t+ e^{\lambda_St} \sigma_S^{} \mathrm{d} B^S(t).
	\end{align}
	Analogously, using the integral form with initial condition $Y^S(\tl) = y_S^{}$ on $[\tl,\tu]$ and multiplying by $e^{-\lambda_S t}$ yields
	\begin{align}
		Y^S(t) = y_S^{} e^{-\lambda_S(t-\tl)} + I_{S,1}(t) +  I_{S,2}(t),
	\end{align}
	with $I_{S,1}(t) = \int_{\tl}^t \sigma_S^{} e^{-\lambda_S (t-r)} \mathrm{d} B^S(r)$ and $I_{S,2}(t)= -\lambda_S c_W^{} \int_{\tl}^t e^{-\lambda_S (t-r)} Y^W(r) \mathrm{d} r$.
	The integral expression $I_{S,2}(t)$ can further be simplified by using the solution for $Y^W(t)$, see \eqref{solution_Y_W},
	\begin{align}
		I_{S,2}(t) &= - \lambda_S c_W^{} \int_{\tl}^t e^{-\lambda_S (t-r)} \bigg( e^{-\lambda_W (r-\tl)} y_W^{} +  \int_{\tl}^r \sigma_W^{} e^{-\lambda_W(r-u)}\mathrm{d} B^W(u) \bigg) \mathrm{d} r \nonumber\\
        &=- \lambda_S c_W^{} \bigg[ I_{S,2.1}(t) + I_{S,2.2}(t) \bigg],
	\end{align}
	with $I_{S,2.1}(t) = \int_{\tl}^t e^{-\lambda_S (t-r)} e^{-\lambda_W (r-\tl)} y_W^{} \mathrm{d}r$ and $I_{S,2.2}(t) = \int_{t_n}^t \int_{\tl}^r \sigma_W^{} e^{-\lambda_S (t-r)} e^{-\lambda_W(r-u)}\mathrm{d} B^W(u)\, \mathrm{d}r$.
	Evaluating each integral individually gives
	\begin{align}
		I_{S,2.1}(t) & = e^{-\lambda_S t} e^{-\lambda_W \tl} y_W^{} \int_{\tl}^t e^{(\lambda_S -\lambda_W) r}  \mathrm{d}r = \frac{e^{-\lambda_S t} e^{-\lambda_W \tl} y_W^{}}{\lambda_S -\lambda_W} ( e^{(\lambda_S -\lambda_W) t} - e^{(\lambda_S -\lambda_W) \tl}) \nonumber\\
		& = \frac{y_W^{}}{\lambda_S -\lambda_W} ( e^{ -\lambda_W( t-\tl)} - e^{ -\lambda_S( t-\tl)}),
	\end{align}
	and by interchanging the order of integration
	\begin{align}
		I_{S,2.2}(t) &= \int_{\tl}^t \int_{u}^t \sigma_W^{} e^{-\lambda_S (t-r)} e^{-\lambda_W(r-u)} \mathrm{d}r \,\mathrm{d} B^W(u)  = \sigma_W^{} e^{-\lambda_S t} \int_{\tl}^t e^{\lambda_W u} \int_{u}^t e^{(\lambda_S -\lambda_W)r} \mathrm{d}r \,\mathrm{d} B^W(u) \nonumber\\ 
		&= \frac{\sigma_W^{} e^{-\lambda_S t}}{\lambda_S -\lambda_W} \int_{\tl}^t e^{\lambda_W u} ( e^{(\lambda_S -\lambda_W)t} - e^{(\lambda_S -\lambda_W)u}) \mathrm{d} B^W(u) \nonumber\\
        &=   \frac{\sigma_W^{} }{\lambda_S -\lambda_W} \int_{\tl}^t  ( e^{ -\lambda_W( t-u)} - e^{ -\lambda_S( t-u)}) \mathrm{d} B^W(u).
	\end{align}
	Combining the results for $I_{S,2.1}(t)$ and $I_{S,2.2}(t)$ yields the solution
	\begin{align}
		\begin{split}
			Y^S(t) &= y_S^{} e^{-\lambda_S(t-\tl)} + \int_{\tl}^t \sigma_S^{} e^{-\lambda_S (t-r)} \mathrm{d} B^S(r) \\
              & \hspace{0.3cm}- \frac{\lambda_S c_W^{}}{\lambda_S -\lambda_W} \bigg( y_W^{} ( e^{ -\lambda_W( t-\tl)} - e^{ -\lambda_S( t-\tl)}) \\
            &\hspace{0.3cm}+  \sigma_W^{}  \int_{\tl}^t  ( e^{ -\lambda_W( t-r)} - e^{ -\lambda_S( t-r)}) \mathrm{d} B^W(r)\bigg).
		\end{split}
		\label{solution_Y_S}
	\end{align}
	Analogous to the calculation for $Y^W(t)$ it is pointed out that both stochastic integrals are martingales with mean zero Gaussian processes. Furthermore, for the expected value of $Y^S(t)$, we get
	\begin{align}
		m_{Y^S}(\tau,y_W^{},y_W^{}) = \mathbb{E}[Y^S(t)] = y_S^{} e^{-\lambda_S\tau} - \tfrac{\lambda_S c_W^{}}{\lambda_S -\lambda_W} y_W^{} ( e^{ -\lambda_W\tau} - e^{ -\lambda_S\tau}),
	\end{align}
	with $\tau = t-\tl$.
	For the variance, only the variances of the stochastic integrals need to consider. Due to the independence of $B^W$ and $B^S$, this results in
	\begin{align}
		\textrm{Var}(Y^S(t)) &= \mathbb{E}\bigg[ \bigg( \int_{\tl}^t \sigma_S^{} e^{-\lambda_S(t-r)}\mathrm{d} B^S(r) \bigg)^2 \bigg] + \mathbb{E}\bigg[ \bigg(  \int_{\tl}^t \tfrac{\lambda_S c_W^{} \sigma_W^{}}{\lambda_S -\lambda_W}  ( e^{ -\lambda_W( t-r)} - e^{ -\lambda_S( t-r)}) \mathrm{d} B^W(r) \bigg)^2 \bigg] \nonumber\\
		&= \mathbb{E}[I_{S,1}(t)^2] +  \mathbb{E}[I_{S,2}(t)^2].
	\end{align}
	The latter follows from the Itô isometry. Note that the first expected value is derived the same way as for $Y^W(t)$, i.e.,
	\begin{align}
		\Sigma_{Y^S}^2(\tau) = \mathbb{E}[I_{S,1}(t)^2] = \mathbb{E}\bigg[ \int_{\tl}^t \sigma_S^2 e^{-2\lambda_S(t-r)}\mathrm{d} r  \bigg] = \frac{\sigma_S^2}{2 \lambda_S} ( 1- e^{-2\lambda_S\tau}).
		\label{App:: Var_Y_S}
	\end{align}
	The second expected value is calculated by multiplying with $( e^{ -\lambda_W( t-r)} - e^{ -\lambda_S( t-r)})^2$ such that
	\begin{align}
		\mathbb{E}[I_{S,2}(t)^2] &= \bigg( \tfrac{\lambda_S c_W^{} \sigma_W^{}}{\lambda_S -\lambda_W}\bigg)^2 \mathbb{E}\bigg[  \int_{\tl}^t \big(  e^{ -2\lambda_W( t-r)} - 2e^{ -\lambda_W( t-r)}e^{ -\lambda_S( t-r)} + e^{ -2\lambda_S( t-r)}\big) \mathrm{d} r  \bigg] \nonumber\\
		&= \bigg( \tfrac{\lambda_S c_W^{} \sigma_W^{}}{\lambda_S -\lambda_W}\bigg)^2 \bigg(\int_{\tl}^t   e^{ -2\lambda_W( t-r)} \mathrm{d} r - 2\int_{\tl}^t e^{ -(\lambda_S+\lambda_W)( t-r)} \mathrm{d} r + \int_{\tl}^t e^{ -2\lambda_S( t-r)} \mathrm{d} r \bigg) \nonumber\\
		& = \bigg( \tfrac{\lambda_S c_W^{} \sigma_W^{}}{\lambda_S -\lambda_W}\bigg)^2 \bigg( \frac{\Sigma_{Y^W}^2(\tau)}{\sigma_W^2} + \frac{\Sigma_{Y^S}^2(\tau)}{\sigma_S^2} -\frac{2}{\lambda_S + \lambda_W} \big(1- e^{-(\lambda_S + \lambda_W)\tau}\big) \bigg) \nonumber\\
		& = \bigg( \tfrac{\lambda_S c_W^{} \sigma_W^{}}{\lambda_S -\lambda_W}\bigg)^2 \bigg( \Sigma_{Y^W}^2(\tau) + \frac{\sigma_W^2}{\sigma_S^2} \Sigma_{Y^S}^2(\tau) -\frac{2\sigma_W^2}{\lambda_S + \lambda_W} \big(1- e^{-(\lambda_S + \lambda_W)\tau}\big) \bigg).
	\end{align}
	To summarize, the variance of $Y^S(t)$ is given by 
	\begin{align}
		\Sigma_S^2(\tau) &= \textrm{Var}(Y^S(t))\nonumber\\ 
        &= \Sigma_{Y^S}^2(\tau) + \bigg( \tfrac{\lambda_S c_W^{} \sigma_W^{}}{\lambda_S -\lambda_W}\bigg)^2 \bigg( \Sigma_{Y^W}^2(\tau) + \frac{\sigma_W^2}{\sigma_S^2} \Sigma_{Y^S}^2(\tau) -\frac{2\sigma_W^2}{\lambda_S + \lambda_W} (1- e^{-(\lambda_S + \lambda_W)\tau}) \bigg).
	\end{align}
	Because $Y^S(t)$ is a sum of independent normal distributed random variables, it holds that
	\begin{align}
		Y^S(t) \sim \mathcal{N}(m_{Y^S}(\tau,y_W^{},y_W^{}) ,\Sigma_S^2(\tau)).
	\end{align}
	
	\myparagraph{$(Y^W,Y^S)$ is bivariate Gaussian} 
	In order to show that $Y^W(t)$ and $Y^S(t)$ define a bivariate Gaussian random variable, one has to shown that $a Y^W(t)+ b Y^S(t)$ is normally distributed for all $a,b \in \mathbb{R}$. Let $a,b \in \mathbb{R}$.
	Let $t \in [\tl,\tu]$, then by the closed-form solutions given $Y^W(\tl) =y_W^{}$ and $Y^S(\tl) =y_S^{}$, we get the linear combination
	\begin{align}
		a Y^W(t)+ b Y^S(t) = a m_{Y^W}^{}(\tau,y_W^{}) +b m_{Y^S}^{}(\tau,y_W^{},y_S^{}) + M^W(t) + M^S(t),
	\end{align}
	with independent Gaussian martingale processes
	\begin{align}
    \begin{split}
		M^W(t) &= \sigma_W^{} a \int_{\tl}^t e^{-\lambda_W(t-r)} - b \tfrac{\lambda_S c_W^{}}{\lambda_S -\lambda_W}     \int_{\tl}^t  ( e^{ -\lambda_W( t-r)} - e^{ -\lambda_S( t-r)}) \mathrm{d} B^W(r), \\
		M^S(t) &= b \sigma_S^{} \int_{\tl}^t  e^{-\lambda_S (t-r)} \mathrm{d} B^S(r).
        \end{split}
	\end{align}
	Note that $M^W(t)$ and $M^S(t)$ are independent normal distributed random variables. Thus, $a Y^W(t)+ b Y^S(t)$ is a normal distributed random variable, implying that ($Y^W(t),Y^S(t))$ is bivariate Gaussian.
	
	\myparagraph{Covariance of $Y^W$ and $Y^S$} 
	The pair $(Y^W(t),Y^S(t))$ is bivariate Gaussian with conditional mean and variance
	\begin{align}
		m_{Y}(\tau,y_W^{},y_S^{}) = \begin{pmatrix}
			m_{Y^W}^{}(\tau,y_W^{}) \\
			m_{Y^S}^{}(\tau,y_W^{},y_S^{})
		\end{pmatrix}, ~~~~ \Sigma_Y(\tau) = \begin{pmatrix}
			\Sigma_W^2(\tau) & \Sigma_{WS}(\tau) \\
			\Sigma_{WS}(\tau) & \Sigma_S^2(\tau)
		\end{pmatrix}.
	\end{align}
	Here, $\Sigma_{WS}(\tau)$ describes the covariance between $Y^W$ and $Y^S$ and is defined by
	\begin{align}
		\Sigma_{WS}(\tau) &= \Cov(Y^W(t),Y^S(t)) = \mathbb{E}[Y^W(t)Y^S(t)] - \mathbb{E}[Y^W(t)]\mathbb{E}[Y^S(t)] \nonumber\\
		&= \mathbb{E}[Y^W(t)Y^S(t)] - m_{Y^W}^{}(\tau,y_W^{})m_{Y^S}^{}(\tau,y_W^{},y_S^{}).
		\label{Covariance proof expression 1}
	\end{align}
	As for the proof above we use the definitions of $I_W(t), I_{S,1}(t)$ and $I_{S,2}(t)$.
	Noticing that all of these expressions are martingales with mean zero. Also $B^W$ and $B^S$ are independent Brownian motions and therefore $I_W(t)$ is independent of $I_{S,1}(t)$ and $I_{S,2}(t)$.
	Using the solutions for $Y^W(t)$ and $Y^S(t)$ in \eqref{solution_Y_W} and \eqref{solution_Y_S}, respectively, leads to the following calculation for $\mathbb{E}[Y^W(t)Y^S(t)]$
	\begin{align}
		\mathbb{E}[Y^W(t)Y^S(t)] &=  \mathbb{E} [ (m_{Y^W}^{}(\tau,y_W^{}) + I_W(t) ) ( m_{Y^S}^{}(\tau,y_W^{},y_S^{}) + I_{S,1}(t) + I_{S,2}(t))] \nonumber\\
		&= \mathbb{E} [m_{Y^W}^{}(\tau,y_W^{}) m_{Y^S}^{}(\tau,y_W^{},y_S^{})] + \mathbb{E} [ m_{Y^W}^{}(\tau,y_W^{})I_{S,1}(t) ] + \mathbb{E} [ m_{Y^W}^{}(\tau,y_W^{}) I_{S,2}(t)] \nonumber\\
		&\hspace{0.3cm} + \mathbb{E} [ I_W(t) m_{Y^S}^{}(\tau,y_W^{},y_S^{})] + \mathbb{E} [ I_W(t)I_{S,1}(t)] + \mathbb{E} [ I_W(t)I_{S,2}(t)] \nonumber\\
		&= m_{Y^W}^{}(\tau,y_W^{}) m_{Y^S}^{}(\tau,y_W^{},y_S^{}) + \mathbb{E} [ I_W(t)] \mathbb{E}[I_{S,1}(t)] + \mathbb{E} [ I_W(t)I_{S,2}(t)] \nonumber\\
		&= m_{Y^W}^{}(\tau,y_W^{}) m_{Y^S}^{}(\tau,y_W^{},y_S^{}) + \mathbb{E} [ I_W(t)I_{S,2}(t)].
	\end{align}
	Substitution into \eqref{Covariance proof expression 1} and applying Itô isometry yield the following result for the covariance
	\begin{align}
		\Sigma_{WS}(\tau) &=   \mathbb{E} [ I_W(t)I_{S,2}(t)] \nonumber\\
		&= \mathbb{E}\bigg[\bigg( \int_{\tl}^t \sigma_W^{} e^{-\lambda_W(t-r)}\mathrm{d} B^W(r) \bigg) \bigg(  - \tfrac{\lambda_S c_W^{}}{\lambda_S -\lambda_W}  \sigma_W^{}  \int_{\tl}^t  ( e^{ -\lambda_W( t-r)} - e^{ -\lambda_S( t-r)}) \mathrm{d} B^W(r) \bigg) \bigg] \nonumber\\
		&= \int_{\tl}^t  - \tfrac{\lambda_S c_W^{}}{\lambda_S -\lambda_W}  \sigma_W^2 e^{-\lambda_W(t-r)}( e^{ -\lambda_W( t-r)} - e^{ -\lambda_S( t-r)}) \mathrm{d}r \nonumber\\
		&=  - \tfrac{\lambda_S c_W^{}}{\lambda_S -\lambda_W}  \sigma_W^2 \int_{\tl}^t  (e^{ -2\lambda_W( t-r)} - e^{ -(\lambda_S+\lambda_W)( t-r)}) \mathrm{d}r \nonumber\\
        & =  - \tfrac{\lambda_S c_W^{}}{\lambda_S -\lambda_W}  \sigma_W^2 \bigg( \int_{\tl}^t e^{ -2\lambda_W( t-r)} \mathrm{d}r -\int_{\tl}^t e^{ -(\lambda_S+\lambda_W)( t-r)} \mathrm{d}r \bigg) \nonumber\\
		&=- \tfrac{\lambda_S c_W^{}}{\lambda_S -\lambda_W}  \sigma_W^2 \bigg( \tfrac{\Sigma_{Y^W}^2(\tau)}{\sigma_W^2} - \tfrac{1}{(\lambda_S+\lambda_W)} \bigg(1- e^{ -(\lambda_S+\lambda_W)\tau} \bigg) \bigg) \nonumber\\
		&= - \tfrac{\lambda_S c_W^{}}{\lambda_S -\lambda_W}  \bigg( \Sigma_{Y^W}^2(\tau) - \tfrac{\sigma_W^2}{(\lambda_S+\lambda_W)} \bigg(1- e^{ -(\lambda_S+\lambda_W)\tau} \bigg) \bigg)
		\label{App:: Cov_WS}
	\end{align}
	To show that $\Sigma_{WS}(\tau) \leq 0$, two different cases must be considered.
    
	\myparagraph{Case 1: $\lambda_S > \lambda_W > 0$}
	In this case, it holds that $\lambda_S - \lambda_W > 0$ and therefore $- \tfrac{\lambda_S c_W^{}}{\lambda_S -\lambda_W} \leq  0.$
	Note that $\lambda_S + \lambda_W > 2 \lambda_W$ and using \eqref{App:: Var_Y_W} this leads to
	\begin{align}
		\Sigma_{Y^W}^2(\tau) = \frac{\sigma_W^2}{2\lambda_W}\big(1- e^{ -2\lambda_W \tau} \big) > \frac{\sigma_W^2}{\lambda_S + \lambda_W}\big(1- e^{ -(\lambda_S+\lambda_W) \tau} \big).
	\end{align}
	It follows from \eqref{App:: Cov_WS} that 
	\begin{align}
		\Sigma_{WS}(\tau) = - \tfrac{\lambda_S c_W^{}}{\lambda_S -\lambda_W}  \bigg( \Sigma_{Y^W}(\tau)^2 - \tfrac{\sigma_W^2}{\lambda_S+\lambda_W} \big(1- e^{ -(\lambda_S+\lambda_W)\tau} \big) \bigg) \leq 0,
	\end{align}
	being zero if and only if $c_W^{} = 0$.
    
	\myparagraph{Case 2: $ \lambda_W >  \lambda_S > 0 $} 
	It is $\lambda_S - \lambda_W < 0$ and gives $- \tfrac{\lambda_S c_W^{}}{\lambda_S -\lambda_W} > 0.$
	Conversely, to the first case, we get 	$\lambda_S + \lambda_W < 2 \lambda_W$, leading to
	\begin{align}
		\Sigma_{Y^W}^2(\tau) = \frac{\sigma_W^2}{2\lambda_W}\big(1- e^{ -(2\lambda_W)\tau} \big) < \frac{\sigma_W^2}{\lambda_S + \lambda_W}\big(1- e^{ -(\lambda_S+\lambda_W)\tau} \big),
	\end{align}
	and moreover this results in
	\begin{align}
		\Sigma_{WS}(\tau) = - \tfrac{\lambda_S c_W^{}}{\lambda_S -\lambda_W}  \bigg( \Sigma_{Y^W}^2(\tau) - \tfrac{\sigma_W^2}{\lambda_S+\lambda_W} \big(1- e^{ -(\lambda_S+\lambda_W)\tau} \big) \bigg) \leq 0,
	\end{align}
	being zero if and only if $c_W^{} = 0$, which completes the proof.
    
	\myparagraph{Positive Definiteness of $\Sigma(\tau)$}
	The matrix $\Sigma(\tau)$ is positive definite if $\trace \Sigma(\tau) >0$ and $\det \Sigma(\tau) >0$. Because of $\Sigma_W^2(\tau), \Sigma_S^2(\tau) > 0$, it is $\trace \Sigma(\tau) = \Sigma_W^2(\tau) + \Sigma_S^2(\tau) > 0$. Furthermore,  $\det \Sigma(\tau) =  \Sigma_W^2(\tau) \Sigma_S^2(\tau) - \Sigma_{WS}^2 > 0$  if $\Sigma_W^2(\tau) \Sigma_S^2(\tau) > \Sigma_{WS}^2(\tau)$. Let $\gamma = \bigg( \tfrac{\lambda_S c_W^{} \sigma_W^{}}{\lambda_S -\lambda_W}\bigg)^2$, then calculations yield
	\begin{align}
		\Sigma_W^2(\tau) \Sigma_S^2(\tau) &= \Sigma_W^2(\tau)\bigg[ \Sigma_{Y^S}^2(\tau) + \gamma \big( \Sigma_{Y^W}^2(\tau) + \frac{\sigma_W^2}{\sigma_S^2} \Sigma_{Y^S}^2(\tau) -\frac{2\sigma_W^2}{\lambda_S + \lambda_W} \big(1- e^{-(\lambda_S + \lambda_W)\tau}\big) \big)\bigg] \nonumber\\
		&= \Sigma_W^2(\tau)\Sigma_{Y^S}^2(\tau) + \gamma \bigg( \Sigma_{Y^W}^4(\tau) + \frac{\sigma_W^2}{\sigma_S^2} \Sigma_W^2(\tau) \Sigma_{Y^S}^2(\tau) -\frac{2 \Sigma_W^2(\tau) \sigma_W^2}{\lambda_S + \lambda_W} \big(1- e^{-(\lambda_S + \lambda_W)\tau} \big) \bigg) \nonumber\\
		&= \Sigma_W^2(\tau)\Sigma_{Y^S}^2(\tau) + \gamma \bigg( \bigg[ \Sigma_W^2(\tau) - \frac{ \sigma_W^2}{\lambda_S + \lambda_W} (1- e^{-(\lambda_S + \lambda_W)\tau}) \bigg]^2 \nonumber\\ 
        & \hspace{0.3cm}+ \frac{\sigma_W^2}{\sigma_S^2} \Sigma_W^2(\tau) \Sigma_{Y^S}^2(\tau)  - \frac{ \sigma_W^4}{(\lambda_S + \lambda_W)^2} (1- e^{-(\lambda_S + \lambda_W)\tau})^2\bigg) \nonumber\\ 
		&= \Sigma_W^2(\tau)\Sigma_{Y^S}^2(\tau) + \Sigma_{WS}^2 + \gamma A,
	\end{align}
	where $A =  \frac{\sigma_W^2}{\sigma_S^2} \Sigma_W^2(\tau) \Sigma_{Y^S}^2(\tau)  - \frac{ \sigma_W^4}{(\lambda_S + \lambda_W)^2} (1- e^{-(\lambda_S + \lambda_W)\tau})^2$. Note that $\Sigma_W^2(\tau)\Sigma_{Y^S}^2(\tau) > 0$ and $\gamma > 0$, so we need to verify that $A \geq 0$ in order to show that $\Sigma$ is positive definite. Using \eqref{App:: Var_Y_W} and \eqref{App:: Var_Y_S}, we obtain
	\begin{align}
		A &= \sigma_W^4\bigg[ \frac{1 -e^{-2\lambda_W}}{2\lambda_W} \frac{1 -e^{-2\lambda_S}}{2\lambda_S} - \bigg( \frac{1 -e^{-(\lambda_W + \lambda_S)}}{\lambda_W + \lambda_S} \bigg)^2\bigg],
	\end{align}
	and by using $\int_{0}^1 e^{-az} \, \mathrm{d}z = \frac{1 -e^{-a}}{a}$ this leads to 
	\begin{align}
		\frac{1 -e^{-(\lambda_W + \lambda_S)}}{\lambda_W + \lambda_S} = \int_{0}^1 e^{-(\lambda_W + \lambda_S)z}  \, \mathrm{dz} =  \int_{0}^1 e^{-\lambda_W z} e^{ -\lambda_Sz} \, \mathrm{dz}.
	\end{align}
	Applying the Cauchy-Schwarz inequality gives
	\begin{align}
		\bigg(\frac{1 -e^{-(\lambda_W + \lambda_S)}}{\lambda_W + \lambda_S}\bigg)^2 \leq \int_{0}^1 e^{-2\lambda_W z} \, \mathrm{dz} \, \int_{0}^1 e^{ -2\lambda_Sz} \, \mathrm{dz} = \frac{1 -e^{-2\lambda_W}}{2\lambda_W} \frac{1 -e^{-2\lambda_S}}{2\lambda_S},
	\end{align}
	which implies $A \geq 0$. Therefore, $\Sigma_W^2(\tau) \Sigma_S^2(\tau) > \Sigma_{WS}^2(\tau)$ and $\Sigma$ is positive definite.  \qed
	\\
	\subsection{Proof of Proposition \ref{prop bivariate dist}} \noindent
	Let $0\leq \tl < \tu \leq \Hor$ with $\tl = t_n$ and $\tu = t_{n+1}$ and $(Y^W_{n},Y^S_n) = (y_W,y_S)$.
	The pair $(\log W_{n+1}, S_{n+1})$ is bivariate Gaussian, as it results from a linear transformation of the bivariate Gaussian $(Y^W_{n+1},Y^S_{n+1})$. Moreover, we have $\Delta t = \tau = t_{n+1}-t_n \equiv const$ and therefore the covariance matrix is given by $\Sigma = \Sigma_Y(\Delta t) = \Sigma_Y(\tau)$. The linear transformation leads to the following mean
	\begin{align}
		m^{WS}_{n+1}(w,s) &= \begin{pmatrix}
			m^W_{n+1}(w) \\
			m^S_{n+1}(w,s)
		\end{pmatrix},
	\end{align}
	with 
	\begin{align}
    \begin{split}
		m^W_{n+1}(w) &= \mu_W^{}(t_{n+1}) +  m_{Y^W}(\Delta t, y_W ), \\
		m^S_{n+1}(w,s) &= \;\mu_S^{}(t_{n+1}) \,+  m_{Y^S}(\Delta t,y_W, y_S ),
        \end{split}
	\end{align}
	where $y_W = \log w -\mu_W^{}(t_{n})$ and $y_S =s-\mu_S^{}(t_{n})$. \qed
	
	\subsection{Proof of Corollary \ref{cor_WS-sequence}} \noindent
	The Cholesky decomposition for $\Sigma$ can be verified by some simple calculations. Furthermore, recalling Proposition \ref{prop bivariate dist}, the pair $(\log W ,S)$ is bivariate Gaussian with mean $m_{WS}$ and covariance matrix $\Sigma$. Based on this, given a standard normally distributed random vector $\noise_{n+1} = (\noise^W_{n+1},\noise^S_{n+1})^\top \sim \mathcal{N}(0_2,I_2)$, we can write $(\log W ,S)$ as a linear transformation of the form
	\begin{align}
		(\log W_{n+1} 	S_{n+1}) 
		=  m^{WS}_{n+1}(W_n,S_n)	
		+ A \noise_{n+1}, \quad \text{ with } A =\begin{pmatrix}
				\Sigma_W & 0 \\
				\rho \Sigma_S & \sqrt{1- \rho^2}\Sigma_S
			\end{pmatrix}.
	\end{align}
	Calculation of the matrix vector product yields \eqref{disccrete dynamics wind and price}. \qed
	
	\section{Calibration of Wind Speed and Energy Price Model}\label{SM:calibration} \noindent
	The wind speed and the energy price are modeled by continuous-time  Ornstein-Uhlenbeck processes on $[0,\Hor]$, see \eqref{ModelDescribtion:OrnsteinUhlenbeckProcesses}. Choosing appropriate parameters is important in order to obtain a realistic dynamical system. Therefore, real-world data is needed to calibrate the model parameters. 
	Data for calibration can be used from the website of Deutscher Wetterdienst\footnote{Deutscher Wetterdienst, Climate data center, \url{https://opendata.dwd.de/climate_environment}} and ENTSO-E\footnote{ENTSO-E, Transparency platform, \url{https://transparency.entsoe.eu/dashboard/show}}. Note that the electricity price data is provided in a resolution of every 15 minutes or every hour. The latter results by averaging the 15 minute values. In our calibration, wind speed and electricity price are provided in an hourly resolution for the entire year 2020, i.e., 8760 data points each. The model assumes that data for wind speed is given by realizations of a random variable $e^{Z(t)}$ for a Gaussian process $Z$. Therefore, the fetched data has to be transformed by applying a logarithmic transformation. Denote the individual data points by $z^\dagger_1,\ldots,z^\dagger_{8760}$ with $\dagger \in \{ W,S \}$, where $z_j^\dagger$ is the data for the $j$-th hour for wind and price respectively. Moreover, outliers need to be removed from the dataset. This is realized by applying the $3$-$\sigma$ rule and remove values outside of the $99.7$ \% confidence interval. 
	
	\subsection{Estimation of the Seasonality Function} \noindent
	For an approximation of the seasonality functions $\mu_W^{}$ and $\mu_S^{}$, an appropriate modeling function is fitted to this data by  least-square regression. 
	The following models where chosen in order capture the daily, half-daily and yearly trends in the underlying data: for wind speed 
	\begin{align}
		&\mu_W^{}(t) = k^W_0 + k^W_1  \cos \bigg( \frac{2\pi(t-t^W_1)}{8760} \bigg) + k^W_2  \cos \bigg( \frac{2\pi(t-t^W_2)}{24} \bigg),
	\end{align}
	and electricity price
	\begin{align}
		& \mu_S^{}(t) = k^S_0 + k^S_1  \cos \bigg( \frac{2\pi(t-t^S_1)}{8760} \bigg) + k^S_2  \cos \bigg( \frac{2\pi(t-t^S_2)}{24} \bigg) + k^S_3  \cos \bigg( \frac{2\pi(t-t^S_3)}{12} \bigg),
	\end{align}
	The denominator in the cosine functions reflects the daily, half-daily and yearly components of the data. The choice of these functions is motivated by the patterns in the data. Wind speeds are commonly higher in winter compared to summer. This fact is also reflected in the daily variations, where wind is slower at night and faster during the day. Electricity prices need to adjust to the changing weather conditions, especially with regard to renewable energy resulting in a yearly seasonality. This adjustments are also needed on a daily scale, arising due to day-to-night cycle of renewable energy sources. The price data reveals that there is a need to include a half-daily trend, due to the fact that most electricity peak load is typically used during two periods in a day: in the morning (people wake up) and in the afternoon (people return from work). The fitted parameters for the seasonalites can be found in Table \ref{parameters seasonalities}. 
	\begin{table}[ht]
		\centering
		\caption{Parameters of the seasonality functions $\mu_W$ and $\mu_S$}
		\begin{tabular}{cccccccc}
			\midrule
			\textbf{Parameters for $\mu_S$} & $k_0^S$  & $k_1^S$ & $k_2^S$ & $k_3^S$ & $t_1^S$ &  $t_2^S$&  $t_3^S$\\
			\midrule
			& -11.2038 & 4.2571 & -6.6642 & 30.4945 & -14782.5 & -6.7823 & -9.5016\\ 
			\midrule
			\textbf{Parameters for $\mu_W$} & $k_0^W$  & $k_1^W$ & $k_2^W$ & $t_1^W$ & $t_2^W$ & &\\
			\midrule
			& 0.1357 & -0.328 & 1.6496 & 1034.1 & 1.1707 & &  \\
		\end{tabular}
		\label{parameters seasonalities}
	\end{table}
	
	\subsection{Estimation of Ornstein-Uhlenbeck Parameters} \noindent
	Recall that $\log W(t) = \mu_W(t) + Y^W(t)$ and $S(t) = \mu(t) + Y^S(t)$, where $Y^W$ and $Y^S$ in \eqref{ModelDescription:OUY_t} are Ornstein-Uhlenbeck processes described in Lemma \ref{Lemma: E of Y^S, Y^W}.
	The discrete-time analogue of the closed-form solution \eqref{model_descibtion::closed_form_soulution_OU} is an autoregressive process of the form
	\begin{align}
    \begin{split}
		Y^W_{n+1} &= p_W Y_n^W + \Sigma_{Y^W} \noise^W_{n+1}, \\
		Y^S_{n+1} &= p_S Y_n^S + q_S Y_n^W+ \Sigma_{S} \noise^S_{n+1},
        \end{split}
		\label{Appendix:AutoregressiveY}
	\end{align}
	with a sequence $(\noise_n)_{n=1,\ldots,N}$ of  independent standard normally distributed random vectors $\noise_{n+1}=(\noise_{n+1}^W,\noise_{n+1}^S)^\top \sim \mathcal{N}(0_2,I_2)$. The parameters $p_W,p_S,q_S,\Sigma_{Y^W}^2$ and $\Sigma_{S}^2$ are identified with their associated quantities given in \eqref{model_describtion:: conditional_mean_variance}, specified as
	\begin{align}
		\begin{split}
			& p_W =  e^{-\lambda_W \Delta t}, \quad \Sigma_{Y^W}^2 = \tfrac{\sigma_W^2}{2 \lambda}(1- e^{-2\lambda_W \Delta t}), \\
			& p_S =  e^{-\lambda_S \Delta t}, \quad q_S = - \frac{\lambda_S c_W^{}}{\lambda_S - \lambda_W}(e^{-\lambda_W\Delta t} - e^{-\lambda_S\Delta t}), \\
			& \Sigma_{S}^2 = \Sigma_{Y^S}^2 + \tfrac{(\lambda_S c_W^{})^2}{(\lambda_S - \lambda_W)^2} \bigg[ \Sigma_{Y^W}^2 + \tfrac{\sigma_W^2}{\sigma_S^2}\Sigma_{Y^S}^2 - \tfrac{2 \sigma_W^2}{\lambda_S+ \lambda_W}(1- e^{-(\lambda_S+\lambda_W) \Delta t}) \bigg].
		\end{split}
		\label{Appendix:AutoregressionParameters}
	\end{align}
	Let $z^W_i,z^S_i,~ i = 1,\ldots,N = 8760$ be the actual time series data in an hourly resolution for $\log(W)$ and $S$, respectively. The realizations $y^W_i = z^W_i - \mu_W^{}(i)$ of $Y^W$ are calculated by removing the seasonality $\mu_W^{}$ from the data. The same holds for the realizations $y^S_i = z^S_i - \mu_S^{}(i)$ of $Y^S$, which are calculated by removing the seasonality $\mu_S^{}$. Regarding the autoregressive structure \eqref{Appendix:AutoregressiveY} and the fact that the random vector $(Y^W,Y^S)^\top$ is bivariate Gaussian, the conditional joint distribution is given by
	\begin{align}
		f_{n}^\theta(y^W,y^S|y_{n-1}^W,y_{n-1}^S) = \frac{1}{2 \pi \Sigma_W \Sigma_S \sqrt{1- \rho^2}} \exp \bigg\{ -\frac{1}{2} A_{n}(\theta) \bigg\},
	\end{align}
	with correlation $\rho$, parameter vector $\theta = (p_W,p_S,q_S,\Sigma_W^2,\Sigma_S^2,\rho)$ and 
	\begin{align}
		A_n(\theta) =  \tfrac{1}{(1-\rho^2)} \bigg[ \tfrac{(y^W-p_Wy_{n-1}^W)^2}{\Sigma_W^2} - 2\rho  \tfrac{(y^W-p_Wy_{n-1}^W)(y^S- p_Sy_{n-1}^S - q_Sy_{n-1}^W)}{\Sigma_W \Sigma_S} +  \tfrac{(y^S - p_Sy_{n-1}^S - q_Sy_{n-1}^W)^2}{\Sigma_S^2} \bigg ].
	\end{align}
	In \cite{roussas2003introduction} a derivation for the maximum-likelihood estimator \cite{Holy2018} of a bivariate Gaussian process can be found. Although our model differs and extends this approach, we use a similar approach to determine the estimator. 
	The maximum-likelihood estimator of the parameter $\theta$ are obtained by maximizing the log-likelihood function
	\begin{align}
		\ell(\theta) = \log \bigg(\prod\limits_{n=1}^N f_{n}^\theta(y^W,y^S|y_{n-1}^W,y_{n-1}^S) \bigg) = \sum \limits_{n=1}^N \log(f_{n}^\theta(y^W,y^S|y_{n-1}^W,y_{n-1}^S)).
	\end{align}
	Note that $\log(f_{n}^\theta(y^W,y^S|y_{n-1}^W,y_{n-1}^S)) = \log( (2 \pi \Sigma_W \Sigma_S \sqrt{1- \rho^2})^{-1}) + A_{n}(\theta)$ and thus
	\begin{align}
		\ell(\theta) = -N(\log(2\pi) + \log(\Sigma_W) + \log(\Sigma_S) + 0.5\log(1-\rho^2)) -\frac{1}{2} \sum \limits_{n=1}^N A_n(\theta).
	\end{align}
	The necessary condition $\nabla_{\theta} \ell(\theta) = 0$, yields a potential maximizer of the log-likelihood function. First, we are going to calculate $\nabla_{\theta} A_n(\theta)$. The partial derivatives are given by
	\begin{align}
    \begin{split}
		\frac{\partial\, A_n(\theta)}{\partial p_W} & = q_2 \bigg[ - \frac{y_{n-1}^W}{\Sigma_W}\frac{(y^W-p_Wy_{n-1}^W)}{\Sigma_W} + \rho \frac{y_{n-1}^W}{\Sigma_W}\frac{(y^S- p_Sy_{n-1}^S - q_Sy_{n-1}^W)}{\Sigma_S} \bigg], \\
		\frac{\partial\, A_n(\theta)}{\partial p_S} & = q_2 \bigg[ - \frac{y_{n-1}^S}{\Sigma_S}\frac{(y^S- p_Sy_{n-1}^S - q_Sy_{n-1}^W)}{\Sigma_S} + \rho \frac{y_{n-1}^S}{\Sigma_S}\frac{(y^W-p_Wy_{n-1}^W)}{\Sigma_W} \bigg], \\
		\frac{\partial\, A_n(\theta)}{\partial q_S} & = q_2 \bigg[ - \frac{y_{n-1}^W}{\Sigma_S}\frac{(y^S- p_Sy_{n-1}^S - q_Sy_{n-1}^W)}{\Sigma_S} + \rho \frac{y_{n-1}^W}{\Sigma_S}\frac{(y^W-p_Wy_{n-1}^W)}{\Sigma_W} \bigg], \\
		\frac{\partial\, A_n(\theta)}{\partial \Sigma_W^2} & = q_1 \bigg[ -\frac{(y^W-p_Wy_{n-1}^W)^2}{\Sigma_W^4} + \rho \frac{(y^W-p_Wy_{n-1}^W)(y^S- p_Sy_{n-1}^S - q_Sy_{n-1}^W)}{\Sigma_W^3 \Sigma_S} \bigg], \\
		\frac{\partial\, A_n(\theta)}{\partial \Sigma_S^2} & = q_1 \bigg[ -\frac{(y^S - p_Sy_{n-1}^S - q_Sy_{n-1}^W)^2}{\Sigma_S^4} + \rho \frac{(y^W-p_Wy_{n-1}^W)(y^S- p_Sy_{n-1}^S - q_Sy_{n-1}^W)}{\Sigma_W \Sigma_S^3} \bigg],
         \end{split}
	\end{align}
with $q_1=\frac{1}{(1-\rho^2)}, ~q_2=\frac{2}{(1-\rho^2)}$ and finally
	\begin{align}
    \begin{split}
		\frac{\partial\, A_n(\theta)}{\partial \rho}  &= \frac{2}{(1-\rho^2)^2} \bigg[\rho\bigg( \frac{(y^W-p_Wy_{n-1}^W)^2}{\Sigma_W^2} + \frac{(y^S - p_Sy_{n-1}^S - q_Sy_{n-1}^W)^2}{\Sigma_S^2} \\
		&\hspace{0.3cm}-(1+\rho^2) \frac{(y^W-p_Wy_{n-1}^W)(y^S - p_Sy_{n-1}^S - q_Sy_{n-1}^W)}{\Sigma_W\Sigma_S} \bigg) \bigg].
        \end{split}
	\end{align}
	This leads to the partial derivatives for $\ell(\theta)$
	\begin{align}
    \begin{split}
		\frac{\partial\, \ell(\theta)}{\partial p_W} & = -\frac{1}{2} \sum \limits_{n=1}^N \frac{\partial\, A_n(\theta)}{\partial p_W},\quad
		\frac{\partial\, \ell(\theta)}{\partial p_S} = -\frac{1}{2} \sum \limits_{n=1}^N \frac{\partial\, A_n(\theta)}{\partial p_S} , \quad 
		\frac{\partial\, \ell(\theta)}{\partial q_S} = -\frac{1}{2} \sum \limits_{n=1}^N \frac{\partial\, A_n(\theta)}{\partial q_S} \\
		\frac{\partial\, \ell(\theta)}{\partial \Sigma_W^2} & = -\frac{N}{2\Sigma_W^2} -\frac{1}{2} \sum \limits_{n=1}^N \frac{\partial\, A_n(\theta)}{\partial \Sigma_W^2},\quad \frac{\partial\, \ell(\theta)}{\partial \Sigma_S^2}  = -\frac{N}{2\Sigma_S^2} -\frac{1}{2} \sum \limits_{n=1}^N \frac{\partial\, A_n(\theta)}{\partial \Sigma_S^2} \\ 
		\frac{\partial\, \ell(\theta)}{\partial \rho} & = \frac{N \rho}{(1-\rho^2)} -\frac{1}{2} \sum \limits_{n=1}^N \frac{\partial\, A_n(\theta)}{\partial \rho}.
        \end{split}
	\end{align}
	We use the following notation to provide a compact derivation of the maximum-likelihood estimator of $\theta$, for this let
	\begin{align}
    \begin{split}
		A &= \sum \limits_{n=1}^N y_n^W y_{n-1}^W,~B = \sum \limits_{n=1}^N (y_{n-1}^W)^2, ~C = \sum \limits_{n=1}^N y_{n-1}^S y_n^S, ~D = \sum \limits_{n=1}^N (y_{n-1}^S)^2,~E = \sum \limits_{n=1}^N y_{n-1}^W y_n^S, \\
		F &= \sum \limits_{n=1}^N y_{n-1}^S y_n^W, ~G = \sum \limits_{n=1}^N y_{n-1}^W y_{n-1}^S, ~H = \sum \limits_{n=1}^N (y_{n}^W-p_Wy_{n-1}^W)^2,\\
		I &= \sum \limits_{n=1}^N (y_{n}^S-p_Sy_{n-1}^S-q_S y_{n-1}^W)^2,~ J = \sum \limits_{n=1}^N (y_{n}^W-p_Wy_{n-1}^W)(y_{n}^S-p_Sy_{n-1}^S-q_S y_{n-1}^W).
        \end{split}
	\end{align}
	Equating the derivatives with respect to $p_W,p_S$ and $q_S$ to zero, we get that
	\begin{align}
    \begin{split}
		\frac{\partial\, \ell(\theta)}{\partial p_W} &= \frac{1}{\Sigma_W^2}(A -p_WB) - \frac{\rho}{\Sigma_W \Sigma_S}(E -p_SG -q_S B) = 0, \\
		\frac{\partial\, \ell(\theta)}{\partial p_S} &= \frac{1}{\Sigma_S^2}(C -p_SD -q_S G) - \frac{\rho}{\Sigma_W \Sigma_S}(F -p_W G) = 0, \\
		\frac{\partial\, \ell(\theta)}{\partial q_S} &= \frac{1}{\Sigma_S^2}(E -p_SG -q_S B) - \frac{\rho}{\Sigma_W \Sigma_S}(A -p_WB) = 0.
        \end{split}
	\end{align}
	Multiplying the first equation by $\frac{\rho \Sigma_W }{\Sigma_W \Sigma_S}$ and subtracting the third, yields
	\begin{align}
    \begin{split}
		&- \frac{\rho^2}{\Sigma_S^2}(E -p_SG -q_S B) - \frac{1}{\Sigma_S^2}(E -p_SG -q_S B) = 0\\
        \Longleftrightarrow~ &\bigg( \frac{\rho^2}{\Sigma_S^2} - \frac{1}{\Sigma_S^2} \bigg) (E -p_S G -q_SB) = 0.
        \end{split}
	\end{align}
	From this, we can deduce that the first and/or second factor needs to be zero. If the first factor is zero, then $\rho^2 = 1$, which by assumption is excluded. Therefore, $E -p_S G -q_SB = 0$ holds, which gives the estimator
	\begin{align}
		\widehat{q}_S = \frac{E -p_S G}{B} = \frac{\sum \limits_{n=1}^N y_{n-1}^W y_n^S -p_S \sum \limits_{n=1}^N y_{n-1}^W y_{n-1}^S}{\sum \limits_{n=1}^N (y_{n-1}^W)^2}.
		\label{qS}
	\end{align}
	Plugging this into the equation for $\frac{\partial\, \ell(\theta)}{\partial p_W} = 0$, we get $\frac{1}{\Sigma_W^2}(A -p_WB) = 0$ and moreover
	\begin{align}
		\widehat{p}_W = \frac{A}{B} = \frac{\sum \limits_{n=1}^N y_n^W y_{n-1}^W}{\sum \limits_{n=1}^N (y_{n-1}^W)^2}.
		\label{pW}
	\end{align}
	Rearranging the second and third equation for $\frac{\rho}{\Sigma_W\Sigma_S}$ takes the form
	\begin{align}
		\frac{\rho}{\Sigma_W\Sigma_S} = \frac{1}{\Sigma_S^2} \frac{C -p_S D -\widehat{q}_SG}{F-\widehat{p}_WG} = \frac{1}{\Sigma_S^2} \frac{A-\widehat{p}_W B}{E -p_SG - \widehat{q}_SB}.
	\end{align}
	Using \eqref{pW}, leads to an expression for $p_S$ given according to 
	\begin{align}
		\widehat{p}_S = \frac{C -\widehat{q}_SG}{D} = \frac{\sum \limits_{n=1}^N y_{n-1}^S y_n^S -\widehat{q}_S\sum \limits_{n=1}^N y_{n-1}^W y_{n-1}^S}{\sum \limits_{n=1}^N (y_{n-1}^S)^2}.
		\label{pS}
	\end{align}
	Next, setting  the derivatives with respect to $\Sigma_W^2,\Sigma_S^2$ and $\rho$ to zero 
	\begin{align}
    \begin{split}
		\frac{\partial\, \ell(\theta)}{\partial \Sigma_W^2} & = -\frac{N}{2\Sigma_W^2} + \frac{H}{2(1-\rho^2)\Sigma_W^4} - \frac{\rho J}{2(1-\rho^2)\Sigma_W \Sigma_S^3}  = 0,\\
		\frac{\partial\, \ell(\theta)}{\partial \Sigma_S^2} & = -\frac{N}{2 \Sigma_S^2} + \frac{I}{2(1-\rho^2)\Sigma_S^4} - \frac{\rho J}{2(1-\rho^2)\Sigma_W^3 \Sigma_S}  = 0,\\
		\frac{\partial\, \ell(\theta)}{\partial \rho} & = \frac{N \rho}{1-\rho^2} -\frac{1}{(1-\rho^2)^2}\bigg( \rho \bigg[\frac{H}{\Sigma_W^2} +\frac{I}{\Sigma_S^2} \bigg] - \frac{1+\rho^2}{\Sigma_W \Sigma_S} J\bigg)=0,
        \end{split}
	\end{align}
	and rearranging with respect to $1-\rho^2$ leading to
	\begin{align}
		1-\rho^2 &= \bigg( \frac{H}{\Sigma_W^2} -\frac{\rho J}{\Sigma_W\Sigma_S}\bigg) \frac{1}{N} = \bigg( \frac{I}{\Sigma_S^2} -\frac{\rho J}{\Sigma_W\Sigma_S}\bigg) \frac{1}{N} = \bigg( \frac{H}{\Sigma_W^2} +\frac{I}{\Sigma_S^2} -\frac{(1+\rho^2) J}{\rho \Sigma_W\Sigma_S}\bigg) \frac{1}{N}.
	\end{align}
	From the second equality, we conclude that $\widehat{\Sigma}_W^2 = \frac{\Sigma_S^2 H}{I}$.
	Plugging this into the third equality, gives
	\begin{align}
		\bigg( \frac{I}{\Sigma_S^2} -\frac{\rho J\sqrt{I}}{\sqrt{H}\Sigma_S^2})\bigg) \frac{1}{N} = \bigg( 2\frac{I}{\Sigma_S^2} -\frac{(1+\rho^2) J\sqrt{I}}{\rho \sqrt{H}\Sigma_S^2})\bigg) \frac{1}{N}.
		\label{secound_relation}
	\end{align}
	Solving for $\rho$ leads to the estimator
	\begin{align}
		\widehat{\rho} = \frac{J}{\sqrt{H I}} = \frac{\sum\limits_{n=1}^N (y_n^W-\widehat{p}_Wy_{n-1}^W)(y_n^S -\widehat{p}_S y_{n-1}^S -\widehat{q}_S y_{n-1}^W)}{\bigg( \sum\limits_{n=1}^N (y_n^W-\widehat{p}_Wy_{n-1}^W)^2 \cdot \sum\limits_{n=1}^N (y_n^S -\widehat{p}_S y_{n-1}^S -\widehat{q}_S y_{n-1}^W)^2 \bigg)^{1/2}}.
	\end{align}
	Using $\widehat{\rho}$ in the first equality in \eqref{secound_relation} implies
	\begin{align}
		1 -\frac{J^2}{HI} = \bigg( \frac{I}{\Sigma_S^2} -\frac{J^2}{H\Sigma_S^2}\bigg) \frac{1}{N} = I\bigg( 1 -\frac{J^2}{HI}\bigg) \frac{1}{N \Sigma_S^2},
	\end{align}
	which results in 
	\begin{align}
		\widehat{\Sigma}_S^2 = \frac{I}{N} = \frac{1}{N} \sum \limits_{n=1}^N (y_n^S -\widehat{p}_Sy_{n-1}^S - \widehat{q}_Sy_{n-1}^W)^2 \quad\text{and}\quad \widehat{\Sigma}_W^2 = \frac{H}{N}= \frac{1}{N} \sum \limits_{n=1}^N (y_n^W -\widehat{p}_Wy_{n-1}^W)^2.
	\end{align}
	Moreover, we get $\widehat{\rho} = \frac{\widehat{\Sigma}_{WS}}{{\widehat{\Sigma}_{W}}{\widehat{\Sigma}_{S}}}$, with the maximum-likelihood estimate of the covariance
	\begin{align}
		\widehat{\Sigma}_{WS} = \frac{1}{N}\sum\limits_{n=1}^N (y_n^W-\widehat{p}_Wy_{n-1}^W)(y_n^S -\widehat{p}_S y_{n-1}^S -\widehat{q}_S y_{n-1}^W).
	\end{align}
	By using the estimators for $\widehat{\theta} =(\widehat{p}_W,\widehat{p}_S,\widehat{q}_S,\widehat{\Sigma}_W^2,\widehat{\Sigma}_S^2, \widehat{\rho})$ and \eqref{Appendix:AutoregressionParameters}, we obtain estimators
	\begin{align}
		\widehat{\lambda}_W = - \frac{\log(\widehat{p}_W)}{\Delta t}, \quad \widehat{\sigma}_W^2 = \frac{2 \widehat{\lambda}_W \widehat{\Sigma}_{W}^2}{1- e^{-2 \widehat{\lambda}_W \Delta t}},
	\end{align}
	for the parameters of the process $Y^W$ and for the process $Y^S$
	\begin{align}
		&\widehat{\lambda}_S = - \frac{\log(\widehat{p}_S)}{\Delta t},\quad  \widehat{c}_W = -\frac{\widehat{q}_S (\widehat{\lambda}_S - \widehat{\lambda}_W)}{\widehat{\lambda}_S(e^{-\widehat{\lambda}_W \Delta t} - e^{-\widehat{\lambda}_S \Delta t})},\quad \widehat{\sigma}_S^2 = \frac{2 \widehat{\lambda}_S (\widehat{\Sigma}_{S}^2 -\kappa)}{1- e^{-2 \widehat{\lambda}_S \Delta t}},
	\end{align}
	where $\widehat{\lambda}_W, \widehat{\sigma}_W^2$ are the maximum-likelihood of $Y^W$ and where
	\begin{align}
		\kappa = \tfrac{(\widehat{\lambda}_S \widehat{c}_W)^2}{(\widehat{\lambda}_S- \widehat{\lambda}_W)^2} \bigg [\widehat{\Sigma}_{Y^W}^2 + \tfrac{\widehat{\sigma}_W^2}{2\widehat{\lambda}_S}(1- e^{-\widehat{\lambda}_S \Delta t}) - \tfrac{2 \widehat{\sigma}_W^2}{\widehat{\lambda}_S+ \widehat{\lambda}_W}(1- e^{-(\widehat{\lambda}_S+ \widehat{\lambda}_W) \Delta t}) \bigg].
	\end{align}
    
	\subsection{Wind Turbine Power}
\begin{figure}[ht]
		\centering
		\includegraphics[scale = 0.53]{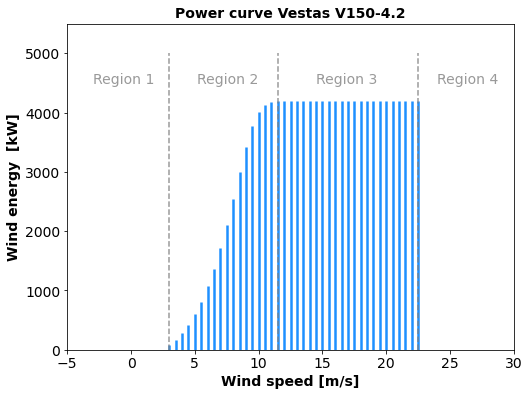}
		\caption{Power curve of a Vestas V150-4.2 with given data (blue bars) and its different operation zones. In region 1 no energy is produced until the cut-in speed of 3 m/s is reached. After that, in region 2 the energy production increases until a certain threshold speed is met. From there, the energy level stays constant at the rated power generation of 4200 \kW, ~in region 3. If the wind speed reaches the cut-out speed of 22.5 m/s, no energy is generated in region 4 in order to prevent damage on the WT.
		}
		\label{model_descibtion::power_curve_data}
	\end{figure}

\begin{figure}[ht]
		\centering
		\includegraphics[scale = 0.53]{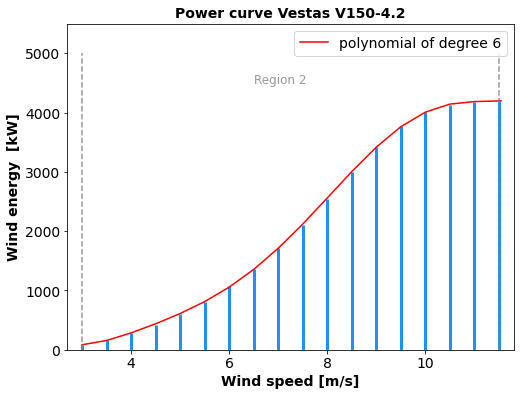}
		\caption{Power curve data of a Vestas V150-4.2 in region 2 (blue bars) and the fitted polynomial of degree 6 (red). More precisely, the interpolated power curve (red), denoted here as $\PW_{R_2}$, determines the generated wind power $\PW$ depending on the wind speed $W$, in region 2.}\label{model_descibtion::power_curve_fit_region_2}
	\end{figure}

\noindent
	\label{sec:Wind_Turbine_Power}
	The WT power output $\PW$ induced by a specific wind speed $W$ can be calculated with its corresponding power curve $P_{\text{WT}}(\cdot)$, via $\PW=P_{\text{WT}}(W)$. In this work, a Vestas V150-4.2\footnote{Vestas, \url{https://www.vestas.com/en/products/4-mw-platform/V150-4-2-MW}} is used for numerical experiments. The corresponding power curve data is shown in Figure \ref{model_descibtion::power_curve_data}, which is available in a database\footnote{The wind power, \url{www.thewindpower.net/turbine_en_1490_vestas_v150-4000-4200.php}}.
	The data of the power curve can be separated into four different regions and suggest a segmented modeling approach for the wind power curve. The regions are classified as follows: Low wind speed, where no energy is produced until the cut-in speed of $w_\text{in} = 3$ m/s is met, denoted as region 1 ($R_1 = [0,w_\text{in})$). From cut-in to a rated level of $w_\text{r} = 11.5$ m/s the WT power generation increases with the underlying wind speed in region 2 ($R_2 = [w_\text{in},w_{\text{Pmax}})$), until reaching its rated power generation of $\PW_{\max} = 4200$ \kW. In region 3 ($R_3 = [w_{\text{Pmax}},w_\text{out})$) the power output remains constant on the rated maximal output. In order to prevent damage on the WT for high speeds, the generation will turn off at a cut-out speed $w_\text{out} = 22.5$ m/s, summarized in region 4 ($R_4 = [w_\text{out},\infty$). Note that the power curve takes constant values except for region 2, where the behavior has to be approximated by a function $\PW_{R_2}(w)$. Moreover, we choose a polynomial of degree 6 \cite{WANG2019109422}, to fit the underlying data points of the power curve in region 2, see Figure \ref{model_descibtion::power_curve_fit_region_2}, i.e.,
	\begin{align}
		\PW_{R_2}(w) = a_6 w^6 + a_5 w^5 +a_3 w^4 +a_3 w^3 +a_2 w^2 +a_1 w +a_0,
		\label{PW_region2_construction}
	\end{align}
	with coefficients $a_0,a_1,a_2,a_3,a_4,a_5,a_6 \in \mathbb{R}$, which can be found in Table \ref{power_curve_region_2_parameters}. We are aware that the behavior in region 2 can be mimicked by a cubic polynomial \cite{WANG2019109422}. However, our motivation for a polynomial of degree 6 is to ensure a smooth transition from the power output in region 2 to the rated power in region 3. To obtain a function that is a representation of the power curve, the following form is used
	\begin{align}
		P_{\text{WT}}(w) = \mathbbm{1}_{R_2}(w) \PW_{R_2}(w) + \mathbbm{1}_{R_3}(w) \PW_{\max}.
	\end{align}
	
	\begin{table}
		\centering
		\caption{Parameters for polynomial function $\PW_{R_2}$}
		\begin{tabular}{cccccccc}
			\midrule
			\textbf{Parameters} & $a_0$  & $a_1$ & $a_2$ & $a_3$ & $a_4$ &  $a_5$&  $a_6$\\
			\midrule
			& 0.1959 & -8.16 & 133.46 & -1101.46 & 4918.22 & -11117.58 & 9941.94 \\
		\end{tabular}
		\label{power_curve_region_2_parameters}
	\end{table}
	
	\section{Construction of Cost Functionals} 
	\subsection{Construction of Running Costs}
	\label{Appedix:Contruction_Cost}
\begin{figure}[ht]
		\centering
		\begin{multicols}{2}
			\includegraphics[scale = 0.39]{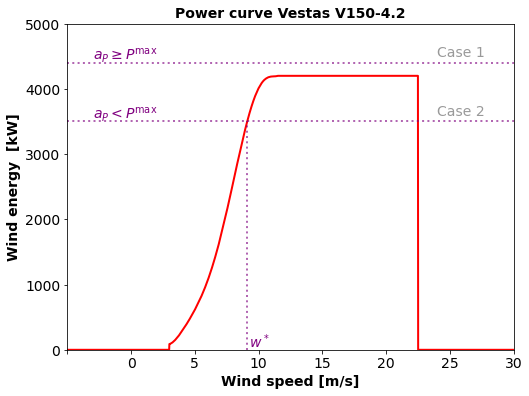} \\
			\hspace{-0.4cm} \includegraphics[scale = 0.39]{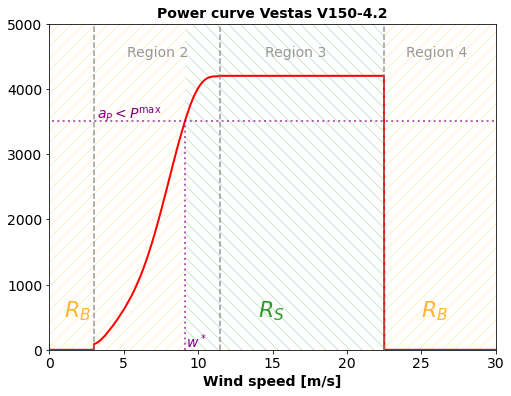}
		\end{multicols}
		\caption{Different cases of the constant HTHP energy consumption $\PH_n$. Case 1 ($\PH_n \geq P^{max}$): Due to the maximal rated power of the WT one always needs to buy additional electricity from the grid. Case 2 ($\PH_n < P^{max}$):  One can find a wind speed $w^* \in R_2$ for which the electricity consumption $\PH_n$ of the HTHP is exactly covered. For all speeds below $w^*$ electricity from the grid needs to be bought. Above this value the WT produces an energy excess, which can be sold. The right plot visualizes the buying and selling regions $R_B$ and $R_S$, respectively.}  
		\label{consumption_cases}
	\end{figure}

	During one time period $[t_n,t_{n+1})$ the HTHP consumes a constant amount of energy $\PH_n$. For all times $t \in [t_n, t_{n+1})$ the electricity consumption must be covered by the available wind energy $\PW(t) = P_{\text{WT}}(W(t))$ (which is a random variable) and the grid power $\PG(t)$, i.e.,
	\begin{align}
		\PH_n = \PW(t) + \PG(t).
	\end{align}
	The expected operational cost at each time step $t \in [t_n, t_{n+1})$ emerges from the price of the used grid power. Given state $x = (r,w,s)$ and action $a$ this yields the operational running costs $C_n(x,a)$ for the time interval $[t_n , t_{n+1})$
	\begin{align}
		C_n(x,a) &= \int_{t_n}^{t_{n+1}} \mathbb{E}_{n,x,a}[S(t)(\PHfun(a) - \PW(t))^+ - \zeta S_{sell}(t)(\PHfun(a) - \PW(t))^-] \mathrm{d} t
	\end{align}
	with $\zeta \in \{ 0,1 \}$, $S_{sell}(t) = S(t) - \eta(t)$ and $\eta:[0,\Hor] \to \mathbb{R}^+$. The amount of energy that needs to be bought is represented by the term $(\PHfun(a) - \PW(t))^+$ and the amount energy being sold by $(\PHfun(a) - \PW(t))^-$.
	Since $\PW(t)$ is bounded from above by $P^{\max}$, two cases are distinguished for the constant value of $\PH_n = \PHfun(\THTout,\THTin)$ as shown in Figure \ref{consumption_cases}. \\[2ex]
    
	\myparagraph{Case one} 
    If $\PH_n \geq P^{\max}$ we always need to buy additional electricity from the grid, resulting in 
	\begin{align}
		(\PH_n - \PW(t))^+ = \PH_n -\PW(t) \quad \text{ and } \quad (\PH_n - \PW(t))^- = 0.
	\end{align}
	Moreover, in this case the running cost is given by
	\begin{align}
		C^1_n(x,a) &= \int_{t_n}^{t_{n+1}} \mathbb{E}_{n,x,a}[S(t)(\PH_n - \PW(t))] \mathrm{d}t = \int_{t_n}^{t_{n+1}} \PH_n\mathbb{E}_{n,x,a}[S(t)]  - \mathbb{E}_{n,x,a}[S(t) \PW(t))] \mathrm{d}t \nonumber\\
		&= \int_{t_n}^{t_{n+1}} \PH_n m^S(\tau,y_W^{},y_S^{})  - \mathbb{E}_{n,x,a}[S(t) \PW(t))] \mathrm{d}t,
		\label{Cost_1}
	\end{align}
	since $S(t) \sim \mathcal{N}(m^S(\tau,y_W^{},y_S^{}) ,\Sigma_S^2(\tau))$ with $\tau = t-t_n$, see Proposition \ref{prop bivariate dist}.\\[2ex]
    
	\myparagraph{Case two} 
    For $\PH_n < P^{\max}$ there exists a $w^* \in R_2$ such that $\PW_{R_2}(w^*) = \PH_n$, which separates the wind speeds into $R_B =  R_1 \cup R^*_{-} \cup R_4 = (0,w_\text{in}) \cup [w_\text{in},w^*) \cup [w_\text{out},\infty) $, meaning to purchase power from the grid and $R_S = R^*_{+} \cup R_3 = [w^*,w_r) \cup [w_r,w_\text{out})$, where it is able to sell unused electricity. In our case, the minimal electricity consumption of the HTHP is greater than the power generated at cut-in of the WT, we get $w^* \geq w_\text{in}$. Moreover, since $\PH_n < P^{\max}$, we also have $w^* < w_r$ and $w^* \in R_2 = [w_\text{in}, w_r)$. Incorporating this to the running costs leads to  
	\begin{align}
    \begin{split}
		C^2_n(x,a) &= \int_{t_n}^{t_{n+1}} \mathbb{E}_{n,x,a}[\mathbbm{1}_{R_B}(W(t))S(t) (\PH_n - \PW(t)) \\
        &\hspace{0.3cm}+ \mathbbm{1}_{R_S}(W(t)) \zeta S_{sell}(t) (\PH_n - \PW(t))] \mathrm{d} t.
        \end{split}
		\label{Cost_2}
	\end{align}
	However, $w^*$ is a zero of $\PW_{R_2}(w) - \PH_n$, which is a polynomial of degree 6 and has at most 6 zeros. By construction, there is a zero in region 2 that needs to be found. Since no closed-form solution for $w^*$ exists, it is computed numerically.
    
	\myparagraph{Calculation of Expected Values} 
	A detailed look at the cost functions \eqref{Cost_1} and \eqref{Cost_2} shows that it is required to determine the expected values for 
	\begin{align}
		\begin{split}
			& \mathbb{E}_{n,x,a}[S(t) \PW(t))],\quad \mathbb{E}_{n,x,a}[\mathbbm{1}_{R_B}(W(t))S(t) \PW(t))], \quad\mathbb{E}_{n,x,a}[\mathbbm{1}_{R_S}(W(t))S(t) \PW(t))], \\
			& \mathbb{E}_{n,x,a}[\mathbbm{1}_{R_B}(W(t)) S(t)], \quad \mathbb{E}_{n,x,a}[\mathbbm{1}_{R_S}(W(t)) S(t)], \quad \mathbb{E}_{n,x,a}[\mathbbm{1}_{R_S}(W(t)) \PW(t)],
		\end{split}
		\label{expected_terms}
	\end{align}
	with
	\begin{align}
		S(t) \PW(t)) = \mathbbm{1}_{R_2}(W(t)) \PW_{R_2}(W(t)) S(t) + \PW_{\max} \mathbbm{1}_{R_3}(W(t)) S(t).
	\end{align}
	Observing that $R_B \cap R_2 = [w_\text{in},w^*)  = R^*_{-},~R_S \cap R_2 = [w^*,w_r) = R^*_{+},~R_B \cap R_3 = \emptyset,~R_S \cap R_3 = R_3$, $\mathbbm{1}_{R_B} = \mathbbm{1}_{R_1} +\mathbbm{1}_{R_{-}^*}  + \mathbbm{1}_{R_4}$ and $\mathbbm{1}_{R_S} = \mathbbm{1}_{R_{+}^*}  + \mathbbm{1}_{R_3}$  
    yields for the expressions in \eqref{expected_terms}
	\begin{align}
    \begin{split}
		\mathbb{E}_{n,x,a}[S(t) \PW(t))] &= \sum_{k=0}^6 a_k \mathbb{E}_{n,x,a}[\mathbbm{1}_{R_{2} }(W(t)) W(t)^k S(t)] \\
        &\hspace{0.3cm}+ P^{\max} \mathbb{E}_{n,x,a}[\mathbbm{1}_{R_3}(W(t)) S(t)], \\
		\mathbb{E}_{n,x,a}[\mathbbm{1}_{R_B}(W(t))S(t) \PW(t))] &= \sum_{k=0}^6 a_k \mathbb{E}_{n,x,a}[\mathbbm{1}_{R^*_{-}}(W(t)) W(t)^k S(t)], \\
		\mathbb{E}_{n,x,a}[\mathbbm{1}_{R_S}(W(t))S(t) \PW(t))] &= \sum_{k=0}^6 a_k \mathbb{E}_{n,x,a}[\mathbbm{1}_{R^*_{+} }(W(t)) W(t)^k S(t)] \\
        &\hspace{0.3cm}+ P^{\max} \mathbb{E}_{n,x,a}[\mathbbm{1}_{R_3}(W(t)) S(t)],\\
		\mathbb{E}_{n,x,a}[\mathbbm{1}_{R_B}(W(t)) S(t)] &=  \mathbb{E}_{n,x,a}[(\mathbbm{1}_{R_1}(W(t)) + \mathbbm{1}_{R^*_{-}}(W(t)) +\mathbbm{1}_{R_4}(W(t))) S(t)], \\
		\mathbb{E}_{n,x,a}[\mathbbm{1}_{R_S}(W(t)) S(t)] &=  \mathbb{E}_{n,x,a}[\mathbbm{1}_{R^*_{+}}(W(t))S(t)] + \mathbb{E}_{n,x,a}[\mathbbm{1}_{R_3}(W(t)) S(t)], \\
		\mathbb{E}_{n,x,a}[\mathbbm{1}_{R_S}(W(t)) \PW(t)] &= \sum_{k=0}^6 a_k \mathbb{E}_{w}[\mathbbm{1}_{R^*_{+} }(W(t)) W(t)^k ] + P^{\max} \mathbb{E}_{w}[\mathbbm{1}_{R_3}(W(t))].
        \end{split}
	\end{align}
	Note that the expressions above are expected values of the form 
	\begin{align}
		E(t,k,R) = \mathbb{E}_{n,x,a}[\mathbbm{1}_{R}(W(t)) W(t)^k S(t)] \quad \text{ and } \quad E_0(t,k,R) = \mathbb{E}_{w}[\mathbbm{1}_{R}(W(t)) W(t)^k],
	\end{align}
	with $R \in \{R_2, R_3,R_4,R^*_{-},R_{B,S},R^*_{+},R_S \}$ and $k \in \mathbb{N}$. Furthermore, $E_0(t,k,R)$ is the special case of $E(t,k,R)$, when the price is replaced with $S(t) = 1$. Hence, it is sufficient to find a closed-form expression for $E(t,k,R)$ in order to get the corresponding expected values in $C_n^1$ and $C_n^2$, respectively.
	In total, the cost function $C_n$ is given by
	\begin{align}
		C_n(x,a) = \begin{cases}
			& C^1_n(x,a),\quad\PH_n \geq P^{\max}, \\
			& C^2_n(x,a),\quad \PH_n < P^{\max},
		\end{cases}
	\end{align}
	with $C_n^1$ and $C_n^2$ specified as
	\begin{align}
		C_n^1(x,a) &= \int_{t_n}^{t_{n+1}} \Psi^1(t,x,a) \mathrm{d} t
		\label{C_n^1 cost},\\
		C^2_n(x,a) &= \int_{t_n}^{t_{n+1}} \Psi^2(t,x,a) \mathrm{d}t,
		\label{C_n^2 cost}
	\end{align}
	with the functions 
	\begin{align}
		\Psi^1(t,x,a) = &\PH_n m^S(\tau,y_W^{},y_S^{})  - \sum_{k=0}^6 a_k E(t,k,R_2) - P^{\max}E(t,0,R_3),\\
		\Psi^2(t,x,a) = &\PH_n(E(t,0,R_B) +\zeta[ E(t,0,R_S)-\eta(t) E_0(t,0,R_S)]) \nonumber\\ &- \sum_{k=0}^6 a_k \big( E(t,k,R^*_{-}) + \zeta[ E(t,k,R^*_{+}) - \eta(t) E_0(t,k,R^*_{+})]\big) \\
		&- \zeta P^{\max}[ E(t,0,R_3) - \eta(t) E_0(t,0,R_3)].\nonumber
	\end{align}
	To keep expressions compact, note that for regions $R_B$ and $R_S$, we have
	\begin{align}
    \begin{split}
		&E(t,0,R_B) =  E(t,0,R_1) + E(t,0,R^*_{-}) + E(t,0,R_4),\\
		&E(t,0,R_S) =  E(t,0,R^*_{+}) + E(t,0,R_3), \\
		&E_0(t,0,R_S) =  E_0(t,0,R^*_{+}) + E_0(t,0,R_3).
        \end{split}
	\end{align}
	
	\myparagraph{Details on $E(t,k,R)$ and $E_0(t,k,R)$} 
	The calculation of $E(t,k,R)$ with $R = [a,b),a < b \in \mathbb{R}$ requires the joint density $f_{WS}$ of $W(t)$ and $S(t)$. For the sake of simplicity, we write $m^W, m^S$ for $m^W(\tau,y_W^{}),m^S(\tau,y_W^{},y_S^{})$ and $\Sigma_W^2,\Sigma_S^2$ for $\Sigma_W^2(\tau), \Sigma_S^2(\tau)$ leading to
	\begin{align}
    \begin{split}
		E(t,k,R)& = \int_{(0,\infty)} \int_\mathbb{R}  \mathbbm{1}_{R}(w) w^k s \,f_{WS}(w,s) \mathrm{d}s\, \mathrm{d}w = \int_{a}^b \int_\mathbb{R}   w^k s \,f_{WS}(w,s) \mathrm{d}s\, \mathrm{d}w \\
		& = \int_{a}^b \int_\mathbb{R}    \, \tfrac{w^k s}{2 \pi \Sigma_W, \Sigma_S \sqrt{1- \rho^2} w} \exp\bigg\{ -\tfrac{1}{2(1-\rho^2)} \bigg[ \tfrac{(\log(w)-m^W)^2}{\Sigma_W^2} \\
        &\hspace{0.3cm}- 2\rho  \tfrac{(\log(w)-m^W)(s-m^S)}{\Sigma_W\Sigma_S}+\tfrac{(s-m^S)^2}{\Sigma_S^2} \bigg ]\bigg\} \mathrm{d}s\, \mathrm{d}w.
        \end{split}
	\end{align}
	Using the substitutions $u = u(w) = \frac{\log(w)-m^W}{\Sigma_W}, v = v(s) = \frac{s-m^S}{\Sigma_S}$ with
	\begin{align}
		\frac{\mathrm{d}u}{\mathrm{d}w} = \frac{1}{w \Sigma_W},~ \frac{\mathrm{d}v}{\mathrm{d}s} = \frac{1}{\Sigma_S} \quad \text{and} \quad  w = e^{\Sigma_W u + m^W},~ s = \Sigma_S v + m^S,
	\end{align}
	results in
	\begin{align}
		E(t,k,R) &= \int_{u(a)}^{u(b)} \int_{\mathbb{R}} \tfrac{w^ks}{2\pi \sqrt{1-\rho^2}} \exp \big\{ -\tfrac{u^2 - 2\rho  uv +  v^2}{2(1-\rho^2)} \big\} \mathrm{d}v\, \mathrm{d}u \nonumber\\
		&= \int_{u(a)}^{u(b)} e^{k(\Sigma_W u + m^W)} \int_{\mathbb{R}} \tfrac{\Sigma_S v + m^S}{2\pi \sqrt{1-\rho^2}} \exp \big\{ -\tfrac{u^2 - 2\rho  uv +  v^2}{2(1-\rho^2)} \big\} \mathrm{d}v\, \mathrm{d}u \nonumber\\
		&= \int_{u(a)}^{u(b)} e^{k(\Sigma_W u + m^W)} \int_{\mathbb{R}} \tfrac{\Sigma_S v + m^S}{2\pi \sqrt{1-\rho^2}} \exp \big\{ -\tfrac{ u^2 - \rho^2 u^2 + \rho^2 u^2 - 2\rho  uv +  v^2}{2(1-\rho^2)} \big\} \mathrm{d}v\, \mathrm{d}u \nonumber\\
		&= \int_{u(a)}^{u(b)} \tfrac{1}{\sqrt{2\pi}} e^{k(\Sigma_W u + m^W)} e^{-\tfrac{u^2}{2}}  I_1\, \mathrm{d}u,
	\end{align}
	with $I_1  = \int_{\mathbb{R}} \tfrac{\Sigma_S v + m^S}{\sqrt{2\pi} \sqrt{1-\rho^2}} \exp \big\{ -\tfrac{  (v -\rho u)^2}{2(1-\rho^2)} \big\} \mathrm{d}v$.
	For solving integral $I_1$, we note that the function in it is the density of a random variable $V \sim \mathcal{N}(\rho u, 1-\rho^2)$, therefore 
	\begin{align}
		I_1 &= \int_{\mathbb{R}} \tfrac{\Sigma_S v }{\sqrt{2\pi} \sqrt{1-\rho^2}} \exp \big\{ -\tfrac{  (v -\rho u)^2}{2(1-\rho^2)} \big\} \mathrm{d}v + \int_{\mathbb{R}} \tfrac{m^S}{\sqrt{2\pi} \sqrt{1-\rho^2}} \exp \big\{ -\tfrac{  (v -\rho u)^2}{2(1-\rho^2)} \big\} \mathrm{d}v \nonumber\\ 
		&=  \Sigma_S \mathbb{E}[V] + m^S \cdot 1 = \Sigma_S \rho u + m^S.
	\end{align}
	Going back to the calculation of $E(t,k,R)$ and substituting the above, we get
	\begin{align}
		E(t,k,R) &= \int_{u(a)}^{u(b)} \tfrac{1}{\sqrt{2\pi}} e^{k(\Sigma_W u + m^W)} e^{-\tfrac{u^2}{2}} (\Sigma_S \rho u + m^S) \mathrm{d}u \nonumber\\
		&= \int_{u(a)}^{u(b)} \tfrac{\Sigma_S \rho u + m^S }{\sqrt{2\pi}} \exp \big\{-\tfrac{u^2-2k\Sigma_W u +(k\Sigma_W)^2 -(k\Sigma_W)^2}{2}  + km^{W}\big\} \mathrm{d}u \nonumber\\
		&= \int_{u(a)}^{u(b)} \tfrac{\Sigma_S \rho u + m^S }{\sqrt{2\pi}} \exp \big\{-\tfrac{(u-k\Sigma_W)^2}{2}+ \tfrac{(k\Sigma_W)^2}{2} + km^{W}\big\} \mathrm{d}u \nonumber\\
		&= e^{\tfrac{(k\Sigma_W)^2}{2} + km^{W}} I_2, 
	\end{align}
	with $I_2 = \int_{u(a)}^{u(b)} \tfrac{\Sigma_S \rho u + m^S }{\sqrt{2\pi}} \exp \big\{-\tfrac{(u-k\Sigma_W)^2}{2}\big\} \mathrm{d}u$.
	Analogous to $I_1$, the function within $I_2$ is the density of the random variable $U \sim \mathcal{N}(k \Sigma_W, 1)$, which yields
	\begin{align}
		I_2 &= \Sigma_S \rho \int_{u(a)}^{u(b)} \tfrac{ u}{\sqrt{2\pi}} e^{-\tfrac{(u-k\Sigma_W)^2}{2}} \mathrm{d}u + m^S \int_{u(a)}^{u(b)}\tfrac{1}{\sqrt{2\pi}} e^{-\tfrac{(u-k\Sigma_W)^2}{2}} \mathrm{d}u \nonumber\\ 
		&= I_3+ m^S[\Phi(u(b)-k\Sigma_W) - \Phi(u(a)- k\Sigma_W)],
	\end{align}
	where $I_3 = \Sigma_S \rho \int_{u(a)}^{u(b)} \tfrac{ u}{\sqrt{2\pi}} e^{-\tfrac{(u-k\Sigma_W)^2}{2}} \mathrm{d}u$.
	The solution of the first integral can be derived as follows:
	\begin{align}
		I_3 &= \Sigma_S \rho \int_{u(a)}^{u(b)} \tfrac{ u - k \Sigma_W + k \Sigma_W}{\sqrt{2\pi}} e^{-\tfrac{(u-k\Sigma_W)^2}{2}} \mathrm{d}u \nonumber\\ 
		&= \Sigma_S \rho \bigg( \int_{u(a)}^{u(b)} \tfrac{ u - k \Sigma_W}{\sqrt{2\pi}} e^{-\tfrac{(u-k\Sigma_W)^2}{2}} \mathrm{d}u + k \Sigma_W \int_{u(a)}^{u(b)} \tfrac{1 }{\sqrt{2\pi}} e^{-\tfrac{(u-k\Sigma_W)^2}{2}} \mathrm{d}u \bigg) \nonumber\\ 
		&= \Sigma_S \rho \bigg( \tfrac{1}{\sqrt{2\pi}}( e^{-\tfrac{(u(a)-k\Sigma_W)^2}{2}} -e^{-\tfrac{(u(b)-k\Sigma_W)^2}{2}} )\nonumber\\
        &\hspace{0.3cm}+ k \Sigma_W[\Phi(u(b)-k\Sigma_W) - \Phi(u(a)- k\Sigma_W)] \bigg),
	\end{align}
	where $\int_{u(a)}^{u(b)} \tfrac{1 }{\sqrt{2\pi}} e^{-\tfrac{(u-k\Sigma_W)^2}{2}} \mathrm{d}u= \Phi(u(b)-k\Sigma_W) - \Phi(u(a)- k\Sigma_W)$. The last equality follows from the substitution $x = -\tfrac{(u-k\Sigma_W)^2}{2}$, as well as by solving the integral. Substituting the expression for $I_3$ into $I_2$ and back into $E(t,k,R)$ gives
	\begin{align}
    \begin{split}
		E(t,k,R) &=  e^{\tfrac{(k\Sigma_W)^2}{2} + km^{W}(t,y_W^{})}\bigg[ \tfrac{\Sigma_S \rho}{\sqrt{2\pi}}( \exp \big\{ -\tfrac{(u(a)-k\Sigma_W)^2}{2} \big\} - \exp\{ -\tfrac{(u(b)-k\Sigma_W)^2}{2} \big\} )\\
		&\hspace{0.3cm}+ (\Sigma_S \rho k \Sigma_W +m^S(t,y_W^{},y_S^{}))[\Phi(u(b)-k\Sigma_W) - \Phi(u(a)- k\Sigma_W)]\bigg],
        \end{split}
	\end{align}
	with $u(w) = \frac{\log(w)-m^W(t,y_W^{})}{\Sigma_W}$. The closed-form expression for $E_0(t,k,R)$ is derived from the expression of $E(t,k,R)$ by letting $m^S \equiv 1$ and $\Sigma_S \equiv 0$, i.e, price is constant with $S(t) = 1$ for all $t \in [t_n, t_{n+1})$.
	
\begin{figure}[ht]
		\centering
		\begin{multicols}{2}
			\hspace*{-1.5cm}\includegraphics[scale = 0.47]{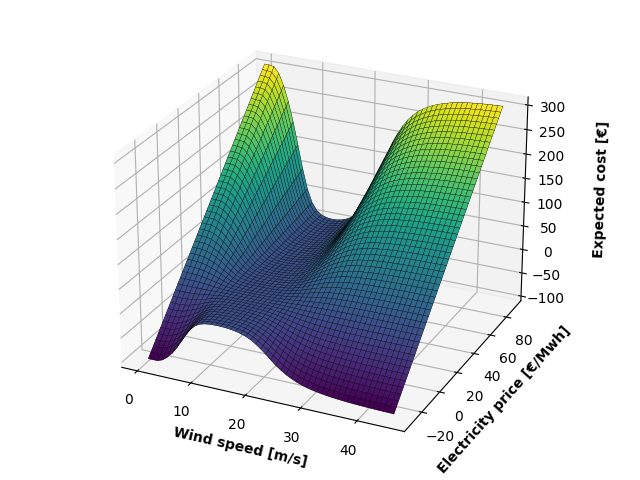}\\
			\includegraphics[scale = 0.47]{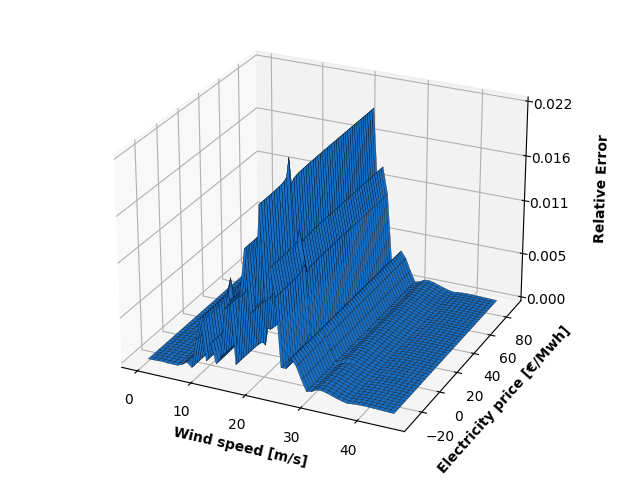}
		\end{multicols}
		\caption{Left: The cost functional $C_n^2$ calculated by an integral solver for $\PH_n = 3500$, with respect to wind speed and electricity price is shown. Right: The relative error between the Gauss-Legendre quadrature for $k^I = 2$ and integral solver is shown.} 
		\label{Surface Cost Functional}
	\end{figure}

	\myparagraph{Numerical Integration for the Cost Functional} 
	The computation of the cost functionals requires the evaluation of time dependent integrals, see \eqref{C_n^1 cost} and  \eqref{C_n^2 cost}. However, there are no closed-form expressions to these integrals and therefore approximation is needed. Built-in integral solvers could be used for this task. These solvers usually achieve high accuracy, but have a slow calculation time. Since both, BDP and Q-learning, make repeated use of the cost functionals, a fast alternative to the integral solver needs to be found, which is comparably accurate but faster in calculation. In our study, we use Gaussian quadrature \cite{Quarteroni2010} with Legendre polynomials to approximate the integral expressions for $C^{1/2}$ in \eqref{C_n^1 cost} and $\eqref{C_n^2 cost}$. For this, let $t_0,..,t_{k^{I}}$ be $k^I + 1$ given distinct interpolation points in the interval $[t_n,t_{n+1}]$. The resulting approximation is given by
	\begin{align}
		C_n^{1/2} = \int_{t_n}^{t_{n+1}} \Psi^{1/2}(t,x,a)\, \mathrm{d}t \approx \frac{\Delta t}{2}\sum_{l=0}^{k^I} \omega_l \Psi^{1/2}(t_l,x,a),
	\end{align}
	with $\Delta t = t_{n+1} -t_n$ and weights $\omega_l,~l=1,\ldots,k^I$.
	We find that $k^I = 2$ is sufficient to obtain good results. 
	Figure \ref{Surface Cost Functional}
	shows the solution to $C_n^2$ obtained by an integral solver and the relative error induced by this Gauss-Legendre quadrature rule. It can be seen that the relative error is less than $2\%$, while acceleration the calculation by a factor of $100$ for the evaluation of the running cost.
	
	\subsection{Construction of Terminal Costs with Penalization and Liquidation} \noindent\label{Appendix::terminal_cost_construction}
	In this work, we consider terminal costs as defined in \eqref{terminal_penalty}. This means, if the terminal TES temperature $R_N = r$ is below or above a critical temperature $r_\text{crit} \in [\rmin,\rmax]$, a penalty or reward with respect to the difference $r_\text{crit} - r$ is applied. The penalty can be interpreted such that the TES has to be charged to the critical value, thus additional costs are incurring. Conversely, the reward results from the liquidation of the energy saved. This is the energy that would be needed to charge the TES from $r_\text{crit}$ to the larger temperature $r$. In both cases, we assign a value to the energy associated with this difference given by a function $g_\text{Ter}$. In turn these are multiplied by predetermined fixed prices $s_\text{Pen}, s_\text{Liq} \geq 0$, respectively, to evaluate the TES filling level. Usually, the fixed prices satisfy $s_\text{Pen} \ge  s_\text{Liq}$. In the following, we will describe the penalization and liquidation process in order to construct $g_\text{Ter}$.\\[1ex]
	For the case $r_\text{crit} > r$, a penalization cost is applied, which is described by the following scenario. Since the TES temperature is below the critical threshold it needs to be charged manually. We assume that for all $t > \Hor$ 
	the produced wind power $\PW(t) = 0$ and that the HTHP has to operate at maximum compressor shaft speed $\dmax$. Therefore, the TES is charged as fast as possible and the energy consumption is covered solely by grid energy. From \eqref{consumption dependency}, the maximum energy consumption $P_\text{max}$ in this case is given by
	\begin{align}
		P_\text{max} = \npump F_2(\TSGout,\mdot,\Tc, d_\text{max}).
	\end{align}
	Moreover, the maximum compressor shaft speed $\dmax$, yields the maximum outlet temperature 
	\begin{align}
		\tHToutmax = F_1(\TSGout, \mdot, \Tc, \dmax).
	\end{align}
	The HTHP run on this continuous-time setting until time $t_\text{Ter} = \Hor + h(r,r_\text{crit})$, for which the TES temperature is raised to $r_\text{crit}$. Here, $h(r,r_\text{crit}) >0$ describes the charging period length, which depends on the current TES state $r$ and the critical value $r_\text{crit}$. From the storage dynamics \eqref{storage_dynamic}, we have
	\begin{align}
		r_\text{crit} = R(t_\text{Pen}) = r + h(r,r_\text{crit}) \frac{\npump \mdot c_\text{p,f}^{}}{m_\text{s}^{} c_\text{p,s}^{}} (\tHToutmax - \TSGin),
	\end{align}
	and determine the length of the charging period according to
	\begin{align}
		h(r,r_\text{crit}) = \frac{m_\text{s}^{} c_\text{p,s}^{} }{\npump \mdot c_\text{p,f}^{} (\tHToutmax - \TSGin)} (r_\text{crit} - r). 
	\end{align}
	The consumed electrical energy associated is given by the function
	\begin{align}
		g_\text{Ter}(r,r_\text{crit}) = h(r,r_\text{crit}) P_\text{max}.
	\end{align}
	The liquidation in the case $r_\text{crit} < r$ results from the saved energy that would be required to charge the TES from the lower critical temperature $r_\text{crit}$ to the higher TES temperature $r$. Analogously to the penalization, we need to charge for a period of length $h(r_\text{crit},r)$, with only difference that the role of $r$ and $r_\text{crit}$ have changed. Thus, the energy savings associated with the liquidation is given by
	\begin{align}
		g_\text{Ter}(r_\text{crit},r) = h(r_\text{crit},r) P_\text{max}.
	\end{align}
	Summarizing and multiplying the energies from penalization and liquidation with its corresponding predetermined prices $s_{Pen}, s_{Liq}$ yields the terminal cost function
	\begin{align}
		G_N(x) =  
		\begin{cases}
			\phantom{-}g_\text{Ter}(r,r_\text{crit})s_\text{Pen}, & r < r_\text{crit}, \\
			- g_\text{Ter}(r_\text{crit},r)s_\text{Liq}, & r \geq, r_\text{crit}
		\end{cases}.
	\end{align}
	
	\section{Details State Space Discretization} \noindent
	\label{Appendix: State Discretization}
	In order to calculate and store the value function with BDP, the state space $\mathcal{X}$ is going to be discretized into grid points. Moreover, we want to use the additive structure and closed-form solutions of the Ornstein-Uhlenbeck processes \eqref{ModelDescribtion:OrnsteinUhlenbeckProcesses}, to construct a family of discretizations $\mathcal{X}_n, ~n=0,\ldots,N$. By this construction we make sure, that the value function is calculated in regions of interested, i.e., states that are located around the seasonalities. In the first step of the construction, we define the reference sets $\mathcal{X}^\text{Ref}_n \subset \mathcal{X}$
	\begin{align}
		\mathcal{X}^\text{Ref}_n = [\rmin,\rmax] \times [\underline{w}_n, \overline{w}_n] \times [\underline{s}_n, \overline{s}_n],
	\end{align}
	with boundary values $\underline{w}_n, \overline{w}_n \in (0,\infty)$ and $\underline{s}_n, \overline{s}_n \in (-\infty,\infty)$. For our specific case we construct these boundaries with the $k^\text{Ref}$-$\sigma$-rule motivated by the normal distribution of the closed-form solutions \eqref{disccrete dynamics wind and price}, leading to
	\begin{align}
    \begin{split}
		\underline{w}_n &= \exp(\mu_W^{}(t_n) - k^\text{Ref} \Sigma_W),~ \overline{w}_n = \exp(\mu_W^{}(t_n) + k^{Ref} \Sigma_W),\\
		\underline{s}_n & = \mu_S^{}(t_n) - k^\text{Ref} \Sigma_S,~ \overline{s}_n = \mu_S^{}(t_n) + k^\text{Ref} \Sigma_S.
        \end{split}
	\end{align}  
	For states in $\mathcal{X}^\text{Ref}_n$, we want to achieve a good approximation, meaning that the errors introduced by the simplifications, namely discretization of states and actions as well as quantization of the expected value, should be small. The error of the state and action space discretization can be reduced by choosing a finer grid, i.e., increasing the number of grid points. However, handling the error introduced by the approximation of the expected value is more difficult. A main cause comes from the quantizer $\mathcal{Z}$. For $x \in \mathcal{X}^\text{Ref}_n$, the quantization points $\mathcal{T}_n(x,a,z)$ with $z \in \mathcal{Z}$ are unlikely to be contained in the discrete set $\widetilde{\mathcal{X}}_{n+1}$. It is therefore necessary to interpolate or extrapolate the value function for these states, whereby not only interpolation errors but also extrapolation errors must be taken into account. In order to reduce extrapolation of the value function for states in the reference set, we enlarge them by choosing boundary values $w_{\min,n} \leq \underline{w}_n,~\overline{w}_n \leq w_{\max,n}$ and $s_{\min,n} \leq \underline{s}_n,~\overline{s}_n \leq s_{\max,n}$. This means, we construct an enlarged set $\mathcal{X}_{n}$ with $\mathcal{X}^\text{Ref}_{n} \subset \mathcal{X}_{n}$ given by
	\begin{align}
		\mathcal{X}_{n} = [\rmin,\rmax] \times [w_{\min,n},w_{\max,n}] \times [s_{\min,n},s_{\max,n}],
	\end{align}
	in the sense that for (almost) every $x \in \mathcal{X}^\text{Ref}_{n-1}$, the transition $\mathcal{T}_{n-1}(x,a,z)$ will be contained in $\mathcal{X}_{n}$. This helps to avoid or reduce the need of extrapolation for states in the reference set. The construction of the values for $w_{\min,n},w_{\max,n},s_{\min,n}$ and $s_{\max,n}$ is again motivated by the normal distribution of the closed-form solution \eqref{disccrete dynamics wind and price}. For $w_{\min,n}$ and $w_{\max,n}$, we apply the $k^\text{Ext}$-$\sigma$-rule leading to
	\begin{align}
		w_{\min,n} = \exp(m^W_n(\underline{y}_{W,n-1}) - k^\text{Ext} \Sigma_W),\quad w_{\max,n} = \exp(m^W_n(\overline{y}_{W,n-1}) + k^\text{Ext} \Sigma_W),
	\end{align}
	for $\underline{y}_{W,n-1} = \log(\underline{w}_{n-1}) - \mu_W^{}(t_{n-1})$ and $\overline{y}_{W,n-1} = \log(\overline{w}_{n-1}) - \mu_W^{}(t_{n-1})$. In order to construct $s_{\min,n}$ and $s_{\max,n}$ we proceed analogously. Considering that the electricity price and the wind speed are negatively correlated, i.e., prices are high if wind speeds are low and visa versa, we end up with the following values
	\begin{align}
		s_{\min,n} = m^S_n(\overline{y}_{W,n-1},\underline{y}_{S,n-1}) - k^\text{Ext} \Sigma_S,\quad s_{\max,n} = m^S_n(\underline{y}_{W,n-1},\overline{y}_{S,n-1}) + k^\text{Ext} \Sigma_S,
	\end{align}
	where $\underline{y}_{W,n-1}$ and $\overline{y}_{W,n-1}$ are as above and $\underline{y}_{S,n-1} = \underline{s}_{n-1} - \mu_S^{}(t_{n-1})$ and $\overline{y}_{S,n-1} = \overline{s}_{n-1} - \mu_S^{}(t_{n-1})$. Note that we set $\mathcal{X}_{0} = \mathcal{X}^\text{Ref}_0$. Given the family of time-varying state spaces $\mathcal{X}_n$ we chose equidistant grid points in each of the intervals $[\rmin,\rmax],~[w_{\min,n},w_{\max,n}]$ and $[s_{\min,n},s_{\max,n}]$. This leads to the discretized state spaces $\widetilde{\mathcal{X}}_n \subset \mathbb{R}^3$ used for the calculations in BDP. In our specific case we set $k^\text{Ref} = 3$ and $k^\text{Ext} = 4$.
	\begin{figure}[ht]
		\centering
		\hspace*{-0.7cm}\includegraphics[scale = 0.35]{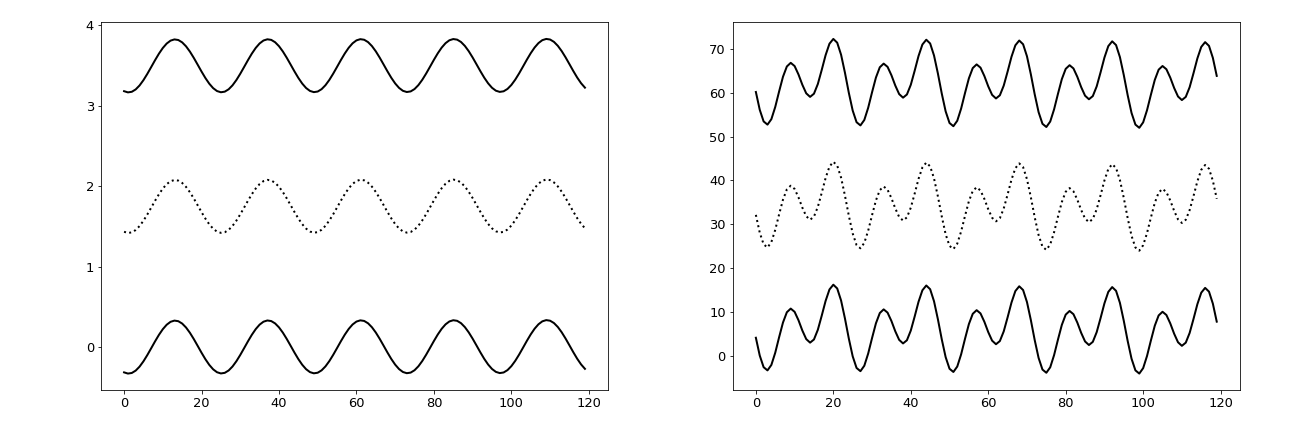}
		\caption{Visualization of the state space boundaries for $\log W$ (left) and $S$ (right). The solid black lines show the boundaries value for $w_{\min,n},w_{\max,n},s_{\min,n}$ and $s_{\max,n}$ over time and the dotted line indicates the seasonality functions $\mu_W$ and $\mu_S$ respetively.}
		\label{state_space_bounds}
	\end{figure}


	\bigskip
	\begin{footnotesize}
		{\neu				
			
			\smallskip\noindent\textbf{Acknowledgments~}
			The authors thank    Ibrahim Mbouandi Njiasse   (BTU Cottbus--Senftenberg)	for valuable discussions that improved this paper.

			\smallskip\noindent
			\textbf{Funding~}	
			E.~Pilling and R.~Wunderlich gratefully acknowledges the  support by the Federal Ministry of Education and Research (BMBF),  award number	05M2022.

		}
	\end{footnotesize}

	\addcontentsline{toc}{section}{References}
	\bibliographystyle{acm}  

\end{document}